# The Hydra-$k$ partial fields

*R.A. Pendavingh[1] and S.H.M. van Zwam[2,3]*

[1]Technische Universiteit Eindhoven, P.O. Box 513, 5600MB Eindhoven, The Netherlands. E-mail: rudi@win.tue.nl

[2]University of Waterloo, 200 University Ave West, Waterloo, ON, Canada, N2L 3G1.

[3]Centrum Wiskunde en Informatica, Amsterdam, The Netherlands. E-mail: Stefan.van.Zwam@cwi.nl

This research was supported by NWO, grants 613.000.561 and 613.000.805.

## Introduction

This report is an addendum to the following two papers, and the following thesis:

**[lift]** R.A. Pendavingh and S.H.M. van Zwam, *Lifts of matroid representations over partial fields* (2010), J. Comb. Theory Ser. B, Vol. 100(1), pp. 36--67.
**[conf]** R.A. Pendavingh and S.H.M. van Zwam, *Confinement of matroid representations to subsets of partial fields* (2010), J. Comb. Theory Ser. B, accepted. Preprint at arXiv:0806.4487
**[thesis]** S.H.M. van Zwam, *Partial Fields in Matroid Theory* (2009), Ph.D. thesis, Eindhoven University of Technology. Available online through http://www.tue.nl/bib/

In these works we introduced the Hydra-$k$ partial fields to study inequivalent representations of quinary matroids. Some of the results in those papers depend on somewhat lengthy calculations; those were left out of the proofs. In this report we carry out these computations, using the *Mathematica* computer algebra system. The code is not involved; in an appendix we summarize some of the less transparent language features that we use, for the benefit of readers not familiar with *Mathematica*.

We assume that the reader is familiar with either the first two or the third of the three works cited above. We will start in the next section with the relevant theoretical results. After that we define a few basic functions that will be of use in the computations. Then we will compute, for each $k \in \{2, 3, 4, 5\}$, the fundamental elements of the partial field $\mathbb{H}_k$, and the automorphism group of $\mathbb{H}_k$. Finally we will verify that $\mathbb{H}_k$ is indeed the lift partial field for the appropriate subset of $GF(5) \times GF(5) \times ... \times GF(5)$-matrices, thus establishing the following result:

**Theorem 1.** Let $M$ be a 3-connected matroid. For each $k \in \{1, 2, 3, 4, 5, 6\}$ the following hold:
1. If $M$ has at least $k$ inequivalent representations over GF(5) then $M$ is representable over $\mathbb{H}_k$.
2. If $M$ has a $U_{2,5}$- or $U_{3,5}$-minor and is representable over $\mathbb{H}_k$, then $M$ has at least $k$ inequivalent representations over GF(5).

# Theoretical results

In this section we list the main theoretical results on which our computations are based. We do not elaborate on the definitions, nor do we give proofs: consult the references above for a detailed treatment. Similarly, the references above will give proper credit to the results being discussed.

**Definition 2.** A partial field is a pair $(R, G)$ of a commutative ring $R$ and a subgroup $G$ of the group of units of $R$ such that $-1 \in G$. We write $p \in \mathbb{P}$ if $p \in G \bigcup \{0\}$.

**Definition 3.** A fundamental element of a partial field $\mathbb{P}$ is an element p such that $1 - p \in \mathbb{P}$. The set of fundamental elements of a partial field is denoted by $\mathcal{F}(\mathbb{P})$. The associates of a fundamental $p$ are $\mathrm{Asc}(p) := \left\{p,\ 1-p,\ \frac{1}{1-p},\ \frac{p}{p-1},\ \frac{p-1}{p},\ \frac{1}{p}\right\}$.

Note that the associates are all fundamental elements.

**Definition 4.** A $\mathbb{P}$-matrix is a matrix such that the determinant of each square submatrix is in $\mathbb{P}$.

**Definition 5.** A partial-field homomorphism $\phi : \mathbb{P} \to \mathbb{P}'$ is a function satisfying $\phi(1) = 1$, $\phi(p)\,\phi(q) = \phi(p\,q)$, and, if $p + q \in \mathbb{P}$ then $\phi(p) + \phi(q) = \phi(p + q)$.

**Lemma 6.** Let $\mathbb{P}_1 = (R_1, G_1)$, $\mathbb{P}_2 = (R_2, G_2)$ be partial fields. Then $\mathbb{P}_1 \times \mathbb{P}_2 := (R_1 \times R_2, G_1 \times G_2)$ is a partial field with fundamental elements $(p, q)$, where $p$ is a fundamental element of $\mathbb{P}_1$ and $q$ is a fundamental element of $\mathbb{P}_2$.

**Lemma 7.** Let $\mathbb{P}_1$, $\mathbb{P}_2$ be partial fields, and $\phi : \mathbb{P}_1 \to \mathbb{P}_2$ a partial-field homomorphism. If $p \in \mathbb{P}_1$ is a fundamental element, then $\phi(p)$ is a fundamental element of $\mathbb{P}_2$.

The following is an application of Whittle's Stabilizer Theorem:

**Lemma 8.** Let $M$ be a 3-connected quinary matroid, and $N$ a minor of $M$ isomorphic to $U_{2,5}$ or $U_{3,5}$. Then each representation over GF(5) of $N$ extends to at most one representation of $M$.

A minor of a matrix $A$ is a matrix obtained from $A$ through a sequence of pivoting, row and column scaling, and deleting of rows and columns, permuting of rows and columns.

**Definition 9.** Let $A$ be a $\mathbb{P}$-matrix. The cross ratios of $A$ are

$$\mathrm{Cr}(A) := \left\{ p : \begin{pmatrix} 1 & 1 \\ p & 1 \end{pmatrix} \text{is minor of } A \right\}.$$

**Definition 10.** Let $\mathbb{P}$, $\hat{\mathbb{P}}$ be partial fields, $\phi : \hat{\mathbb{P}} \to \mathbb{P}$ a partial-field homomorphism, and $A$ a $\mathbb{P}$-matrix. A lifting function for $\phi$ is a function $^{\uparrow} : \mathrm{Cr}(A) \to \hat{\mathbb{P}}$ satisfying, for all $p, q \in \mathrm{Cr}(A)$,

1. $\phi\left(p^{\uparrow}\right) = p$
2. if $p + q = 1$ in $\mathbb{P}$ then $p^{\uparrow} + q^{\uparrow} = 1$ in $\hat{\mathbb{P}}$
3. if $p\,q = 1$ in $\mathbb{P}$ then $p^{\uparrow}\,q^{\uparrow} = 1$ in $\hat{\mathbb{P}}$.

Without going into details (for which we refer to **[lift]**), a local lift of $A$ is a $\hat{\mathbb{P}}$-matrix $\hat{A}$ that is obtained from $A$ through the lifting function, satisfying $\phi(\hat{A})$ is scaling-equivalent to $A$; it is a global lift if the lifting construction commutes with pivoting in $A$. The following is a special case:

**Lemma 11.** Let $A = \begin{pmatrix} 1 & 1 & 1 \\ 1 & p & q \end{pmatrix}$ be a $\mathbb{P}$-matrix. Then $A$ has a local lift if and only if $q^{\uparrow} / p^{\uparrow} = (q/p)^{\uparrow} \in \mathcal{F}(\hat{\mathbb{P}})$. In that case the local lift is $\begin{pmatrix} 1 & 1 & 1 \\ 1 & p^{\uparrow} & q^{\uparrow} \end{pmatrix}$, which is also a global lift.

**Theorem 12 (Lift Theorem).** Let $\mathbb{P}$, $\hat{\mathbb{P}}$ be partial fields, $\phi : \hat{\mathbb{P}} \to \mathbb{P}$ a partial-field homomorphism, let $A$ be a $\mathbb{P}$-matrix, and let ${}^{\uparrow} : \mathrm{Cr}(A) \to \mathcal{F}(\hat{\mathbb{P}})$ be a lifting function for $\phi$. Exactly one of the following is true:

1. $A$ has a global lift $\hat{A}$;
2. $A$ has a minor $D$ such that there is no local lift of $D$, and $D$ or $D^T$ is equal to $\begin{pmatrix} 1 & 1 & 0 & 1 \\ 1 & 0 & 1 & 1 \\ 0 & 1 & 1 & 1 \end{pmatrix}$ or $\begin{pmatrix} 1 & 1 & 1 \\ 1 & p & q \end{pmatrix}$, for some $p, q \in \mathbb{P}$.

We will apply this theorem to a class $\mathcal{A}$ of $\mathbb{P}$-matrices, which is closed under minors and duality.

Next up is the Confinement Theorem, which we will use to limit the number of cross ratios that need lifting.

**Definition 13.** Let $\mathbb{P}$ be a partial field and $\mathbb{P}'$ a sub-partial field of $\mathbb{P}$. Then $\mathbb{P}'$ is induced if, for all $p, q \in \mathbb{P}'$, $p + q \in \mathbb{P}$ implies $p + q \in \mathbb{P}'$.

In particular, if $\mathbb{P}' = \mathbb{P} \cap R'$ for a subring $R'$ of the ring containing $\mathbb{P}$, then certainly $\mathbb{P}'$ is induced.

**Theorem 14 (Confinement Theorem).** Let $\mathbb{P}$ be a partial field, $\mathbb{P}'$ an induced sub-partial field, and $D$ a 3-connected scaled $\mathbb{P}'$-matrix. Let $A$ be a 3-connected $\mathbb{P}$-matrix with a $D$-minor. Exactly one of the following is true:

1. $A$ is a scaled $\mathbb{P}'$-matrix;
2. $A$ has a 3-connected minor $A'$ such that
   a. $A'$ is not a scaled $\mathbb{P}'$-matrix;
   b. $D$ can be obtained from $A'$ by deleting at most one row and at most one column;
   c. if $A'$ has both a row and a column more than $D$, then for at least one of them its deletion is 3-connected.

In this report we will only apply the following corollary:

**Corollary 15.** Let $\mathbb{P} = (R, G)$ be a partial field, and $p \in \mathcal{F}(\mathbb{P})$. Let $D = \begin{pmatrix} 1 & 1 \\ p & 1 \end{pmatrix}$ be a $\mathbb{P}$-matrix. Let $R'$ be the subring of $R$ generated by $\mathrm{Cr}(D)$, and $\mathbb{P}' := \mathbb{P} \cap R'$. Let $A$ be a $\mathbb{P}$-matrix with a $D$-minor. If there is no $q$ such that $\begin{pmatrix} 1 & 1 & 1 \\ 1 & p & q \end{pmatrix}$ is a minor of $A$ or $A^T$, then $A$ is a scaled $\mathbb{P}'$-matrix.

*Proof.* Clearly $\mathbb{P}'$ is an induced sub-partial field. Since $A$ has no minor that is a 3-connected one-element extension of $D$, Case (2) of the Confinement Theorem does not hold, and therefore Case (1) must hold. □

Starting with the next section we will interleave the Mathematica commands that carry out computations, their output, and discussion of the results being computed.

# Basic functions

The variables we will use. Make sure nothing gets assigned to them by accident:

```
Protect[α, β, γ]
```

```
{α, β, γ}
```

The associates of a fundamental element or list of fundamental elements:

```
associates[1] := {0, 1};
associates[0] := {0, 1};
associates[S_List] := Union @@ (associates /@ S);
associates[p_] := Union[Simplify[{p, 1 - p, 1/(1 - p), p/(p - 1), (p - 1)/p, 1/p}]];
```

Example:

```
associates[{1, 2}]
```

$$\left\{-1, 0, \frac{1}{2}, 1, 2\right\}$$

We will be working with rational functions of polynomials. Testing if two such expressions are equal is not straightforward, since some terms may have moved around. We use Mathematica's built-in arithmetic to deal with this.

```
equality[x_, y_] := Simplify[x - y] == 0;
```

```
occursInList[L_List, x_] := Or @@ (equality[x, #] & /@ L);
```

Example, comparing Mathematica's default function, which only does pattern matching, with our implementation:

```
MemberQ[{1, α, -1/(α - 1), -α/(α^3 + 1)}, -(-1/(1 - α))]
```

```
False
```

This can be explained by the way Mathematica stores these expressions.

```
InputForm[-1/(α - 1)]
```

```
-(-1 + α)^(-1)
```

```
InputForm[-(-1/(1 - α))]
```

```
(1 - α)^(-1)
```

Checking if the difference equals zero is more reliable:

```
occursInList[{1, α, -1/(α - 1), -α/(α^3 + 1)}, -(-1/(1 - α))]
```

```
True
```

Note that occursInList is not particularly efficient. This is no problem until we get to Hydra-5. For that partial field the various computations we do take up a number of hours. Since these computations are only one-time checks, we did not make much effort to optimize them.

Converting a rational number to GF(5):

```
toGF5[n_Integer] := Mod[n, 5];
```

```
GF5inv[0] := 0;
GF5inv[1] := 1;
GF5inv[2] := 3;
GF5inv[3] := 2;
GF5inv[4] := 4;
```

```
toGF5[Rational[n_Integer, m_Integer]] := Mod[n GF5inv[Mod[m, 5]], 5];
```

```
toGF5[S_List] := toGF5 /@ S;
```

Example:

```
toGF5[14/72]
```

```
2
```

Note that this conversion is somewhat risky, as common factors of 5 may have been removed from fractions, causing undefined fractions to be mapped to GF(5). However, in the instances where we use this function this is of no concern.

```
toGF5[30/20]
```

```
4
```

# Hydra-2

## Definition and homomorphisms

The Hydra-2 partial field is $\mathbb{H}_2 := \left(\mathbb{Z}\left[\frac{1}{2}, i\right], \langle i, 1-i\rangle\right)$, where $i$ is a root of $x^2+1$. Elements are of the form $\pm 2^k i^l (1-i)^m$. This partial field embeds in the complex plane in a natural way.

```
ListPlot[{Re[#], Im[#]} & /@
  Flatten[Table[(-1)^x 2^y i^z (1 - i)^v, {x, 0, 1}, {y, -2, 2}, {z, -2, 2}, {v, -2, 2}]],
 {AspectRatio → 1, PlotStyle → PointSize[Large]}]
```

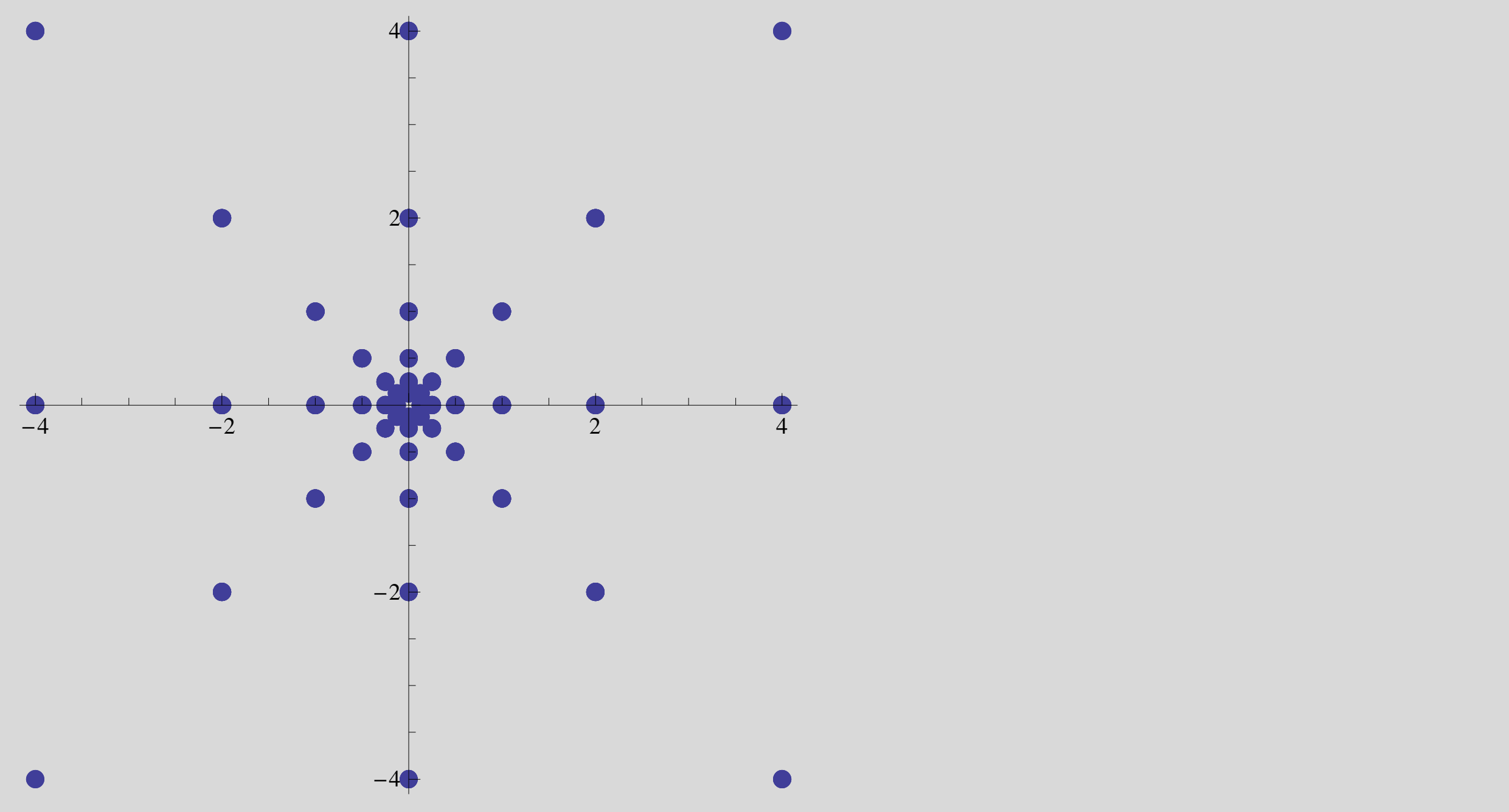

The fundamental elements of $\mathbb{H}_2$ are

```
H2funs = Together[associates[{0, 1, 2, i}]]
```

```
{-1, 0, -i, i, 1/2, 1/2 - i/2, 1/2 + i/2, 1, 1 - i, 1 + i, 2}
```

There are two homomorphisms $\mathbb{H}_2 \to \mathrm{GF}(5)$. We collect them in the partial-field homomorphism $\phi 2 : \mathbb{H}_2 \to \mathrm{GF}(5) \times \mathrm{GF}(5)$ defined by $\phi 2(i) = (2, 3)$.

```
ϕ2[-1] = {4, 4}; ϕ2[0] = {0, 0}; ϕ2[-i] = {3, 2}; ϕ2[i] = {2, 3};
ϕ2[1/2] = {3, 3}; ϕ2[1/2 - i/2] = {2, 4}; ϕ2[1/2 + i/2] = {4, 2};
ϕ2[1] = {1, 1}; ϕ2[1 - i] = {4, 3}; ϕ2[1 + i] = {3, 4}; ϕ2[2] = {2, 2};
```

## Fundamental elements

The complex argument of every element is a multiple of $\pi/4$. For an element $p$ to be fundamental, $p-1$ must also be in the partial field. Suppose, first, that $\mathrm{im}(p) > 0$. There are three choices for the complex argument. Since $\arg(p-1) > \arg(p)$, it follows that we must have $\{\arg(p), \arg(p-1)\} \in \{(\pi/4, \pi/2), (\pi/4, 3\pi/4), (\pi/2, 3\pi/4)\}$. The corresponding fundamental elements are, respectively, $\left\{1+i, \frac{1+i}{2}, i\right\}$. A similar argument holds if $\mathrm{im}(p) < 0$.

If $p$ is real, then $m$, the exponent of $(1-i)$, has to be a multiple of 4. But $(1-i)^4 = -4$. From this it follows that each such $p$ must be a power of two, and hence the real fundamental elements of $\mathbb{H}_2$ are those of the Dyadic partial field, namely $\left\{2, 1, \frac{1}{2}, 0, -1\right\}$. Hence we have proven

**Lemma 16.** The fundamental elements of $\mathbb{H}_2$ are
$\left\{0, 1, -1, 2, 1/2, i, i+1, \frac{i+1}{2}, 1-i, \frac{1-i}{2}, -i\right\}$.

For the next sections, note that the maximum of the norms of these elements is 2, and the minimum norm of the nonzero elements is 1/2:

```
Max[Abs[H2funs]]
```

```
2
```

```
Min[Abs[Complement[H2funs, {0}]]]
```

```
1
─
2
```

## The automorphism group

The partial-field automorphism group of $\mathbb{H}_2$ has order 2: it contains the identity and the automorphism that exchanges $i$ and $-i$. This is easy to see since each automorphism must keep 0 and 1 fixed, and therefore the real number line must be fixed.

In fact, the partial-field homomorphism $\phi 2$ defined above induces a $1-1$ correspondence between the automorphisms of $\mathbb{H}_2$ and coordinate permutations of $\mathrm{GF}(5) \times \mathrm{GF}(5)$.

## Representations of $U_{2,5}$

We prove the following fact:

**Lemma 17.** Let $M$ be a 3-connected, $\mathbb{H}_2$-representable matroid. If $M$ has a $U_{2,5}$- or $U_{3,5}$-minor then $M$ has at least two inequivalent representations over GF(5).

To prove this it suffices to show that each $\mathbb{H}_2$-representation of $U_{2,5}$ gives rise to two inequivalent representations over GF(5). Since any $\mathbb{H}_2$-representation matrix of $M$ must contain one of these as a minor, the result follows. First we enumerate all such representations. By normalizing and suppressing the identity matrix at the front, we see

```
A = {{1, 1, 1}, {1, p, q}}; MatrixForm[A]
```

$$\begin{pmatrix} 1 & 1 & 1 \\ 1 & p & q \end{pmatrix}$$

This is an $\mathbb{H}_2$-matrix representing $U_{2,5}$ if and only if $p$, $q$, $\frac{p}{q}$ are fundamental elements, $p$ and $q$ are not equal to 0 or 1, and $p \neq q$. Moreover, two such matrices are equivalent if and only if they are equal.

```
nonzerooneH2funs = Complement[H2funs, {0, 1}]
```

$$\left\{-1, -i, i, \frac{1}{2}, \frac{1}{2} - \frac{i}{2}, \frac{1}{2} + \frac{i}{2}, 1 - i, 1 + i, 2\right\}$$

Create a table of all candidate pairs $p$, $q$:

```
Flatten[Table[{nonzerooneH2funs[[x]], nonzerooneH2funs[[y]]},
   {x, 1, Length[nonzerooneH2funs]}, {y, 1, x - 1}], 1] ⋃
 Flatten[Table[{nonzerooneH2funs[[x]], nonzerooneH2funs[[y]]},
   {x, 1, Length[nonzerooneH2funs]}, {y, x + 1, Length[nonzerooneH2funs]}], 1]
```

$$\Big\{\{-1, -i\}, \{-1, i\}, \left\{-1, \frac{1}{2}\right\}, \left\{-1, \frac{1}{2} - \frac{i}{2}\right\}, \left\{-1, \frac{1}{2} + \frac{i}{2}\right\}, \{-1, 1 - i\}, \{-1, 1 + i\}, \{-1, 2\},$$
$$\{-i, -1\}, \{-i, i\}, \left\{-i, \frac{1}{2}\right\}, \left\{-i, \frac{1}{2} - \frac{i}{2}\right\}, \left\{-i, \frac{1}{2} + \frac{i}{2}\right\}, \{-i, 1 - i\}, \{-i, 1 + i\}, \{-i, 2\},$$
$$\{i, -1\}, \{i, -i\}, \left\{i, \frac{1}{2}\right\}, \left\{i, \frac{1}{2} - \frac{i}{2}\right\}, \left\{i, \frac{1}{2} + \frac{i}{2}\right\}, \{i, 1 - i\}, \{i, 1 + i\}, \{i, 2\},$$
$$\left\{\frac{1}{2}, -1\right\}, \left\{\frac{1}{2}, -i\right\}, \left\{\frac{1}{2}, i\right\}, \left\{\frac{1}{2}, \frac{1}{2} - \frac{i}{2}\right\}, \left\{\frac{1}{2}, \frac{1}{2} + \frac{i}{2}\right\}, \left\{\frac{1}{2}, 1 - i\right\}, \left\{\frac{1}{2}, 1 + i\right\}, \left\{\frac{1}{2}, 2\right\},$$
$$\left\{\frac{1}{2} - \frac{i}{2}, -1\right\}, \left\{\frac{1}{2} - \frac{i}{2}, -i\right\}, \left\{\frac{1}{2} - \frac{i}{2}, i\right\}, \left\{\frac{1}{2} - \frac{i}{2}, \frac{1}{2}\right\}, \left\{\frac{1}{2} - \frac{i}{2}, \frac{1}{2} + \frac{i}{2}\right\}, \left\{\frac{1}{2} - \frac{i}{2}, 1 - i\right\},$$
$$\left\{\frac{1}{2} - \frac{i}{2}, 1 + i\right\}, \left\{\frac{1}{2} - \frac{i}{2}, 2\right\}, \left\{\frac{1}{2} + \frac{i}{2}, -1\right\}, \left\{\frac{1}{2} + \frac{i}{2}, -i\right\}, \left\{\frac{1}{2} + \frac{i}{2}, i\right\}, \left\{\frac{1}{2} + \frac{i}{2}, \frac{1}{2}\right\},$$
$$\left\{\frac{1}{2} + \frac{i}{2}, \frac{1}{2} - \frac{i}{2}\right\}, \left\{\frac{1}{2} + \frac{i}{2}, 1 - i\right\}, \left\{\frac{1}{2} + \frac{i}{2}, 1 + i\right\}, \left\{\frac{1}{2} + \frac{i}{2}, 2\right\}, \{1 - i, -1\}, \{1 - i, -i\},$$
$$\{1 - i, i\}, \left\{1 - i, \frac{1}{2}\right\}, \left\{1 - i, \frac{1}{2} - \frac{i}{2}\right\}, \left\{1 - i, \frac{1}{2} + \frac{i}{2}\right\}, \{1 - i, 1 + i\}, \{1 - i, 2\}, \{1 + i, -1\},$$
$$\{1 + i, -i\}, \{1 + i, i\}, \left\{1 + i, \frac{1}{2}\right\}, \left\{1 + i, \frac{1}{2} - \frac{i}{2}\right\}, \left\{1 + i, \frac{1}{2} + \frac{i}{2}\right\}, \{1 + i, 1 - i\}, \{1 + i, 2\},$$
$$\{2, -1\}, \{2, -i\}, \{2, i\}, \left\{2, \frac{1}{2}\right\}, \left\{2, \frac{1}{2} - \frac{i}{2}\right\}, \left\{2, \frac{1}{2} + \frac{i}{2}\right\}, \{2, 1 - i\}, \{2, 1 + i\}\Big\}$$

Filter out those for which $p/q$ is not fundamental or equal to 1:

```
pqPairs = Select[%, MemberQ[nonzerooneH2funs, #[[1]] / #[[2]]] &]
```

$$\Big\{\{-1, -i\}, \{-1, i\}, \{-i, -1\}, \{-i, i\}, \left\{-i, \frac{1}{2} - \frac{i}{2}\right\}, \{-i, 1 - i\}, \{i, -1\},$$
$$\{i, -i\}, \left\{i, \frac{1}{2} + \frac{i}{2}\right\}, \{i, 1 + i\}, \left\{\frac{1}{2}, \frac{1}{2} - \frac{i}{2}\right\}, \left\{\frac{1}{2}, \frac{1}{2} + \frac{i}{2}\right\}, \left\{\frac{1}{2} - \frac{i}{2}, -i\right\},$$
$$\left\{\frac{1}{2} - \frac{i}{2}, \frac{1}{2}\right\}, \left\{\frac{1}{2} - \frac{i}{2}, \frac{1}{2} + \frac{i}{2}\right\}, \left\{\frac{1}{2} - \frac{i}{2}, 1 - i\right\}, \left\{\frac{1}{2} + \frac{i}{2}, i\right\}, \left\{\frac{1}{2} + \frac{i}{2}, \frac{1}{2}\right\},$$
$$\left\{\frac{1}{2} + \frac{i}{2}, \frac{1}{2} - \frac{i}{2}\right\}, \left\{\frac{1}{2} + \frac{i}{2}, 1 + i\right\}, \{1 - i, -i\}, \left\{1 - i, \frac{1}{2} - \frac{i}{2}\right\}, \{1 - i, 1 + i\},$$
$$\{1 - i, 2\}, \{1 + i, i\}, \left\{1 + i, \frac{1}{2} + \frac{i}{2}\right\}, \{1 + i, 1 - i\}, \{1 + i, 2\}, \{2, 1 - i\}, \{2, 1 + i\}\Big\}$$

```
Length[%]
```

```
30
```

To show that each of these gives rise to two inequivalent representations over GF(5) we consider the partial-field homomorphism

$\phi 2 : \mathbb{H}_2 \to \mathrm{GF}(5) \times \mathrm{GF}(5)$ defined in the first section.

Let $\xi_i : \mathrm{GF}(5) \times \mathrm{GF}(5) \to \mathrm{GF}(5)$ be the projection onto the $i$th coordinate. We have to show that $\xi_1(\phi 2(A))$ is not equivalent to $\xi_2(\phi 2(A))$ for all matrices $A$ computed above. We only need to look at the pairs $p$, $q$. These are their images under $\phi 2$:

```
Map[ϕ2, pqPairs, {2}]
```

```
{{{4, 4}, {3, 2}}, {{4, 4}, {2, 3}}, {{3, 2}, {4, 4}}, {{3, 2}, {2, 3}}, {{3, 2}, {2, 4}},
 {{3, 2}, {4, 3}}, {{2, 3}, {4, 4}}, {{2, 3}, {3, 2}}, {{2, 3}, {4, 2}}, {{2, 3}, {3, 4}},
 {{3, 3}, {2, 4}}, {{3, 3}, {4, 2}}, {{2, 4}, {3, 2}}, {{2, 4}, {3, 3}}, {{2, 4}, {4, 2}},
 {{2, 4}, {4, 3}}, {{4, 2}, {2, 3}}, {{4, 2}, {3, 3}}, {{4, 2}, {2, 4}}, {{4, 2}, {3, 4}},
 {{4, 3}, {3, 2}}, {{4, 3}, {2, 4}}, {{4, 3}, {3, 4}}, {{4, 3}, {2, 2}}, {{3, 4}, {2, 3}},
 {{3, 4}, {4, 2}}, {{3, 4}, {4, 3}}, {{3, 4}, {2, 2}}, {{2, 2}, {4, 3}}, {{2, 2}, {3, 4}}}
```

Indeed, each pair of representations is inequivalent:

```
Map[#[[1]] - #[[2]] &, %, {2}]
```

```
{{0, 1}, {0, -1}, {1, 0}, {1, -1}, {1, -2}, {1, 1}, {-1, 0}, {-1, 1}, {-1, 2}, {-1, -1},
 {0, -2}, {0, 2}, {-2, 1}, {-2, 0}, {-2, 2}, {-2, 1}, {2, -1}, {2, 0}, {2, -2}, {2, -1},
 {1, 1}, {1, -2}, {1, -1}, {1, 0}, {-1, -1}, {-1, 2}, {-1, 1}, {-1, 0}, {0, 1}, {0, -1}}
```

```
Position[%, {0, 0}]
```

```
{}
```

This completes the proof of Lemma 17.

## Lifting

We prove the following fact:

**Lemma 18.** Let $M$ be a 3-connected matroid with at least two inequivalent representations over GF(5). Then $M$ is representable over $\mathbb{H}_2$.

Consider a $\mathrm{GF}(5) \times \mathrm{GF}(5)$-matrix $A$ representing $M$ such that $\xi_1(A)$ and $\xi_2(A)$ are inequivalent. We want to apply the Lift Theorem, so we start by constructing a lifting function for $\mathrm{Cr}(A)$. Consider the homomorphism $\phi 2$ from the previous section. The restriction of $\phi 2$ to $\mathcal{F}(\mathbb{H}_2)$ is a bijection between $\mathcal{F}(\mathbb{H}_2)$ and $\mathcal{F}(\mathrm{GF}(5) \times \mathrm{GF}(5))$:

```
ϕ2 /@ H2funs
```

```
{{4, 4}, {0, 0}, {3, 2}, {2, 3}, {3, 3}, {2, 4}, {4, 2}, {1, 1}, {4, 3}, {3, 4}, {2, 2}}
```

It follows immediately that the function $\mathcal{F}(\mathrm{GF}(5) \times \mathrm{GF}(5)) \to \mathbb{H}_2$ that is the inverse of $\phi 2$ is a lifting function.

```
up2[x_, y_] := H2funs[[Position[ϕ2 /@ H2funs, {x, y}][[1, 1]]]];
up2[{x_, y_}] := up2[x, y];
```

Some examples:

```
up2[1, 1]
```

```
1
```

```
up2[2, 4]
```

```
1   i
— - —
2   2
```

Now suppose that $A$ has no global lift. Then the Lift Theorem (Theorem 12) tells us that $A$ or $A^T$ has a minor of the form

$$\begin{pmatrix} 1 & 1 & 0 & 1 \\ 1 & 0 & 1 & 1 \\ 0 & 1 & 1 & 1 \end{pmatrix} \text{ or } \begin{pmatrix} 1 & 1 & 1 \\ 1 & p & q \end{pmatrix}, \text{ for some } p, q \in \mathrm{Cr}(A).$$

without a local lift. The first of these certainly does have a local lift, since the dyadic partial field is contained in $\mathbb{H}_2$, and the non-Fano matroid is dyadic. Therefore we have to ensure that all $\mathrm{GF}(5) \times \mathrm{GF}(5)$-representations of $U_{2,5}$ that can occur in $A$ have a local lift. Let $D$ be such a minor. Suppose that $\xi_1(D)$ and $\xi_2(D)$ are equivalent. By the Stabilizer Theorem, and in particular Lemma 8, this implies that $\xi_1(A)$ and $\xi_2(A)$ are equivalent, a contradiction to our choice of $A$. Hence $\xi_1(D)$ and $\xi_2(D)$ are inequivalent. This gives us 30 representations of $U_{2,5}$. Again we only list the pairs $(p, q)$ in

$$D = \begin{pmatrix} 1 & 1 & 1 \\ 1 & p & q \end{pmatrix}.$$

```
U25reps = { {2, 3}, {2, 4}, {3, 2}, {3, 4}, {4, 2}, {4, 3}}
```

```
{{2, 3}, {2, 4}, {3, 2}, {3, 4}, {4, 2}, {4, 3}}
```

```
U25reppairs = Transpose /@ Permutations[U25reps, {2}]
```

```
{{{2, 2}, {3, 4}}, {{2, 3}, {3, 2}}, {{2, 3}, {3, 4}}, {{2, 4}, {3, 2}}, {{2, 4}, {3, 3}},
 {{2, 2}, {4, 3}}, {{2, 3}, {4, 2}}, {{2, 3}, {4, 4}}, {{2, 4}, {4, 2}}, {{2, 4}, {4, 3}},
 {{3, 2}, {2, 3}}, {{3, 2}, {2, 4}}, {{3, 3}, {2, 4}}, {{3, 4}, {2, 2}}, {{3, 4}, {2, 3}},
 {{3, 2}, {4, 3}}, {{3, 2}, {4, 4}}, {{3, 3}, {4, 2}}, {{3, 4}, {4, 2}}, {{3, 4}, {4, 3}},
 {{4, 2}, {2, 3}}, {{4, 2}, {2, 4}}, {{4, 3}, {2, 2}}, {{4, 3}, {2, 4}}, {{4, 4}, {2, 3}},
 {{4, 2}, {3, 3}}, {{4, 2}, {3, 4}}, {{4, 3}, {3, 2}}, {{4, 3}, {3, 4}}, {{4, 4}, {3, 2}}}
```

```
Length[%]
```

```
30
```

The lifts look as follows:

```
Map[up2, U25reppairs, {2}]
```

$\{\{2, 1+i\}, \{i, -i\}, \{i, 1+i\}, \{\frac{1}{2}-\frac{i}{2}, -i\}, \{\frac{1}{2}-\frac{i}{2}, \frac{1}{2}\}, \{2, 1-i\},$
$\{i, \frac{1}{2}+\frac{i}{2}\}, \{i, -1\}, \{\frac{1}{2}-\frac{i}{2}, \frac{1}{2}+\frac{i}{2}\}, \{\frac{1}{2}-\frac{i}{2}, 1-i\}, \{-i, i\}, \{-i, \frac{1}{2}-\frac{i}{2}\},$
$\{\frac{1}{2}, \frac{1}{2}-\frac{i}{2}\}, \{1+i, 2\}, \{1+i, i\}, \{-i, 1-i\}, \{-i, -1\}, \{\frac{1}{2}, \frac{1}{2}+\frac{i}{2}\},$
$\{1+i, \frac{1}{2}+\frac{i}{2}\}, \{1+i, 1-i\}, \{\frac{1}{2}+\frac{i}{2}, i\}, \{\frac{1}{2}+\frac{i}{2}, \frac{1}{2}-\frac{i}{2}\}, \{1-i, 2\}, \{1-i, \frac{1}{2}-\frac{i}{2}\},$
$\{-1, i\}, \{\frac{1}{2}+\frac{i}{2}, \frac{1}{2}\}, \{\frac{1}{2}+\frac{i}{2}, 1+i\}, \{1-i, -i\}, \{1-i, 1+i\}, \{-1, -i\}\}$

And indeed, for each of these $p^{\uparrow} / q^{\uparrow}$ is a fundamental element:

```
Select[%, MemberQ[nonzerooneH2funs, #[[1]] / #[[2]]] &]
```

$\{\{2, 1+i\}, \{i, -i\}, \{i, 1+i\}, \{\frac{1}{2}-\frac{i}{2}, -i\}, \{\frac{1}{2}-\frac{i}{2}, \frac{1}{2}\}, \{2, 1-i\},$
$\{i, \frac{1}{2}+\frac{i}{2}\}, \{i, -1\}, \{\frac{1}{2}-\frac{i}{2}, \frac{1}{2}+\frac{i}{2}\}, \{\frac{1}{2}-\frac{i}{2}, 1-i\}, \{-i, i\}, \{-i, \frac{1}{2}-\frac{i}{2}\},$
$\{\frac{1}{2}, \frac{1}{2}-\frac{i}{2}\}, \{1+i, 2\}, \{1+i, i\}, \{-i, 1-i\}, \{-i, -1\}, \{\frac{1}{2}, \frac{1}{2}+\frac{i}{2}\},$
$\{1+i, \frac{1}{2}+\frac{i}{2}\}, \{1+i, 1-i\}, \{\frac{1}{2}+\frac{i}{2}, i\}, \{\frac{1}{2}+\frac{i}{2}, \frac{1}{2}-\frac{i}{2}\}, \{1-i, 2\}, \{1-i, \frac{1}{2}-\frac{i}{2}\},$
$\{-1, i\}, \{\frac{1}{2}+\frac{i}{2}, \frac{1}{2}\}, \{\frac{1}{2}+\frac{i}{2}, 1+i\}, \{1-i, -i\}, \{1-i, 1+i\}, \{-1, -i\}\}$

```
Length[%]
```

```
30
```

It follows that $A$ does have a global lift. This completes the proof of Lemma 18.

# Hydra-3

## Definition and homomorphisms

The Hydra-3 partial field is $\mathbb{H}_3 = \left(\mathbb{Z}\left[\alpha, \frac{1}{1-\alpha}, \frac{1}{\alpha^2-\alpha+1}\right], \left\langle -1, \alpha, 1-\alpha, \alpha^2-\alpha+1 \right\rangle\right)$, where $\alpha$ is an indeterminate.

The generators of $\mathbb{H}_3$ are

```
H3gens = {-1, α, 1 - α, α^2 - α + 1};
```

The fundamental elements of $\mathbb{H}_3$ are

```
H3funs = Together[associates[{1, α, α^2 - α + 1, α^2/(α - 1), -α/(α - 1)^2}]]
```

$$\left\{0,\ 1,\ \frac{1-\alpha+\alpha^2}{\alpha^2},\ \frac{1}{1-\alpha},\ 1-\alpha,\ \frac{-1+\alpha}{\alpha^2},\ \frac{1}{\alpha},\ \frac{-1+\alpha}{\alpha},\ -\frac{(-1+\alpha)^2}{\alpha},\ \alpha,\ -\frac{\alpha}{(-1+\alpha)^2},\right.$$
$$\frac{\alpha}{-1+\alpha},\ -(-1+\alpha)\,\alpha,\ \frac{\alpha^2}{-1+\alpha},\ \frac{1-\alpha+\alpha^2}{\alpha},\ \frac{1-2\,\alpha+\alpha^2}{1-\alpha+\alpha^2},\ \frac{1-\alpha+\alpha^2}{(-1+\alpha)^2},\ \frac{1}{(1-\alpha)\,\alpha},\ \frac{1}{1-\alpha+\alpha^2},$$
$$\left.\frac{\alpha}{1-\alpha+\alpha^2},\ \frac{(-1+\alpha)\,\alpha}{1-\alpha+\alpha^2},\ \frac{\alpha^2}{1-\alpha+\alpha^2},\ 1-\alpha+\alpha^2,\ \frac{1-\alpha+\alpha^2}{(-1+\alpha)\,\alpha},\ \frac{1-\alpha}{1-\alpha+\alpha^2},\ \frac{-1+\alpha-\alpha^2}{-1+\alpha}\right\}$$

The fundamental elements other than 0, 1 play a special role:

```
nonzerooneH3funs = Complement[H3funs, {0, 1}];
```

There are three homomorphisms $\mathbb{H}_3 \to \mathrm{GF}(5)$. We collect them in the partial-field homomorphism $\phi 3 : \mathbb{H}_3 \to \mathrm{GF}(5) \times \mathrm{GF}(5) \times \mathrm{GF}(5)$ defined by $\phi 3(\alpha) = (2, 3, 4)$.

```
ϕ3[0] = {0, 0, 0};
ϕ3[1] = {1, 1, 1};
ϕ3[-1] = {4, 4, 4};
ϕ3[x_] := toGF5[x /. α → {2, 3, 4}];
```

This is indeed a partial-field homomorphism:

```
ϕ3 /@ H3gens
```

```
{{4, 4, 4}, {2, 3, 4}, {4, 3, 2}, {3, 2, 3}}
```

## Fundamental elements

Our aim in this section is to verify the following statement:

**Lemma 19.** The fundamental elements of $\mathbb{H}_3$ are
$\mathrm{Asc}\left\{1,\ \alpha,\ \alpha^2-\alpha+1,\ \frac{\alpha^2}{\alpha-1},\ \frac{-\alpha}{(\alpha-1)^2}\right\}$.

Elements of this partial field are of the form $\pm\alpha^x(1-\alpha)^y\left(\alpha^2-\alpha+1\right)^z$. First we have to show that the elements in the Lemma are indeed fundamental. We inspect $1-p$ for all $p$ in this list:

```
Together[1 - {1, α, α^2 - α + 1, α^2/(α - 1), -α/(α - 1)^2}]
```

```
{0, 1 - α, α - α^2, (-1 + α - α^2)/(-1 + α), (1 - α + α^2)/(-1 + α)^2}
```

These are all powers of the generators, so all elements of $\left\{1, \alpha, \alpha^2 - \alpha + 1, \frac{\alpha^2}{\alpha-1}, \frac{-\alpha}{(\alpha-1)^2}\right\}$ are indeed fundamental. It follows immediately that all their associates are fundamental.

To show that this list is indeed complete, we first try to bound the exponents of fundamental elements. Consider the homomorphism $\phi : \mathbb{H}_3 \to \mathbb{H}_2$ defined by $\phi(\alpha) = i$.

```
α^x (1 - α)^y (α^2 - α + 1)^z /. α → ⅈ
```

```
(-ⅈ)^z ⅈ^x (1 - ⅈ)^y
```

The only contribution to the norm of such an element is by $(1-i)^y$, from which it follows that $-2 \le y \le 2$, since each fundamental element is mapped to a fundamental element by a homomorphism. Next, consider the homomorphism $\psi$ such that $\psi(\alpha) = 1 - i$.

```
α^x (1 - α)^y (α^2 - α + 1)^z /. α → (1 - ⅈ)
```

```
(-ⅈ)^z ⅈ^y (1 - ⅈ)^x
```

The only contribution to the norm of such an element is by $(1-i)^x$, from which it follows that $-2 \le x \le 2$. Next, consider the homomorphism $\rho$ such that $\rho(\alpha) = \frac{1-i}{2}$.

```
α^x (1 - α)^y (α^2 - α + 1)^z /. α → (1 - ⅈ)/2
```

```
(1 - ⅈ)^x (1 + ⅈ)^y 2^(-x-y-z)
```

The norm of this expression equals $2^{-1/2\,x - 1/2\,y - z}$. Since this is at least $1/2$ and at most 2, and both $x$ and $y$ have known bounds, we deduce that $-3 \le z \le 3$. Now we have reduced the set of possible exponents of fundamental elements to a finite list.

```
candidateH3funExps =
  Flatten[Table[{s, x, y, z}, {s, 0, 1}, {x, -2, 2}, {y, -2, 2}, {z, -3, 3}], 3];
```

These are the corresponding elements of the partial field:

```
candidateH3funs = {0} ⋃ Together[(Times @@ (H3gens^#)) & /@ candidateH3funExps]
```

$$\left\{-1,\ 0,\ 1,\ \frac{1}{1-\alpha},\ 1-\alpha,\ -\frac{1}{(-1+\alpha)^2},\ \frac{1}{(-1+\alpha)^2},\ \frac{1}{-1+\alpha},\ -1+\alpha,\ -(-1+\alpha)^2,\ (-1+\alpha)^2,\ -\frac{1}{\alpha^2},\ \frac{1}{\alpha^2},\right.$$

$$\frac{1-\alpha}{\alpha^2},\ -\frac{1}{(-1+\alpha)^2\,\alpha^2},\ \frac{1}{(-1+\alpha)^2\,\alpha^2},\ -\frac{1}{(-1+\alpha)\,\alpha^2},\ \frac{1}{(-1+\alpha)\,\alpha^2},\ \frac{-1+\alpha}{\alpha^2},\ -\frac{(-1+\alpha)^2}{\alpha^2},\ \frac{(-1+\alpha)^2}{\alpha^2},$$

$$-\frac{1}{\alpha},\ \frac{1}{\alpha},\ \frac{1-\alpha}{\alpha},\ -\frac{1}{(-1+\alpha)^2\,\alpha},\ \frac{1}{(-1+\alpha)^2\,\alpha},\ -\frac{1}{(-1+\alpha)\,\alpha},\ \frac{1}{(-1+\alpha)\,\alpha},\ \frac{-1+\alpha}{\alpha},\ -\frac{(-1+\alpha)^2}{\alpha},$$

$$\frac{(-1+\alpha)^2}{\alpha},\ -\alpha,\ \alpha,\ -\frac{\alpha}{(-1+\alpha)^2},\ \frac{\alpha}{(-1+\alpha)^2},\ -\frac{\alpha}{-1+\alpha},\ \frac{\alpha}{-1+\alpha},\ -(-1+\alpha)\,\alpha,\ (-1+\alpha)\,\alpha,\ -(-1+\alpha)^2\,\alpha,$$

$$(-1+\alpha)^2\,\alpha,\ -\alpha^2,\ \alpha^2,\ -\frac{\alpha^2}{(-1+\alpha)^2},\ \frac{\alpha^2}{(-1+\alpha)^2},\ -\frac{\alpha^2}{-1+\alpha},\ \frac{\alpha^2}{-1+\alpha},\ -(-1+\alpha)\,\alpha^2,\ (-1+\alpha)\,\alpha^2,$$

$$-(-1+\alpha)^2\,\alpha^2,\ (-1+\alpha)^2\,\alpha^2,\ \frac{1}{-1+\alpha-\alpha^2},\ -1+\alpha-\alpha^2,\ \frac{-1+\alpha-\alpha^2}{(-1+\alpha)^2},\ \frac{-1+\alpha-\alpha^2}{-1+\alpha},\ \frac{-1+\alpha-\alpha^2}{\alpha^2},$$

$$\frac{-1+\alpha-\alpha^2}{(-1+\alpha)^2\,\alpha^2},\ \frac{-1+\alpha-\alpha^2}{(-1+\alpha)\,\alpha^2},\ \frac{-1+\alpha-\alpha^2}{\alpha},\ \frac{-1+\alpha-\alpha^2}{(-1+\alpha)^2\,\alpha},\ \frac{-1+\alpha-\alpha^2}{(-1+\alpha)\,\alpha},\ -\frac{1}{\left(1-\alpha+\alpha^2\right)^3},\ \frac{1}{\left(1-\alpha+\alpha^2\right)^3},$$

$$\frac{1-\alpha}{\left(1-\alpha+\alpha^2\right)^3},\ -\frac{1}{(-1+\alpha)^2\left(1-\alpha+\alpha^2\right)^3},\ \frac{1}{(-1+\alpha)^2\left(1-\alpha+\alpha^2\right)^3},\ -\frac{1}{(-1+\alpha)\left(1-\alpha+\alpha^2\right)^3},$$

$$\frac{1}{(-1+\alpha)\left(1-\alpha+\alpha^2\right)^3},\ \frac{-1+\alpha}{\left(1-\alpha+\alpha^2\right)^3},\ -\frac{(-1+\alpha)^2}{\left(1-\alpha+\alpha^2\right)^3},\ \frac{(-1+\alpha)^2}{\left(1-\alpha+\alpha^2\right)^3},\ -\frac{1}{\alpha^2\left(1-\alpha+\alpha^2\right)^3},$$

$$\frac{1}{\alpha^2\left(1-\alpha+\alpha^2\right)^3},\ \frac{1-\alpha}{\alpha^2\left(1-\alpha+\alpha^2\right)^3},\ -\frac{1}{(-1+\alpha)^2\,\alpha^2\left(1-\alpha+\alpha^2\right)^3},\ \frac{1}{(-1+\alpha)^2\,\alpha^2\left(1-\alpha+\alpha^2\right)^3},$$

$$-\frac{1}{(-1+\alpha)\,\alpha^2\left(1-\alpha+\alpha^2\right)^3},\ \frac{1}{(-1+\alpha)\,\alpha^2\left(1-\alpha+\alpha^2\right)^3},\ \frac{-1+\alpha}{\alpha^2\left(1-\alpha+\alpha^2\right)^3},\ -\frac{(-1+\alpha)^2}{\alpha^2\left(1-\alpha+\alpha^2\right)^3},$$

$$\frac{(-1+\alpha)^2}{\alpha^2\left(1-\alpha+\alpha^2\right)^3},\ -\frac{1}{\alpha\left(1-\alpha+\alpha^2\right)^3},\ \frac{1}{\alpha\left(1-\alpha+\alpha^2\right)^3},\ \frac{1-\alpha}{\alpha\left(1-\alpha+\alpha^2\right)^3},\ -\frac{1}{(-1+\alpha)^2\,\alpha\left(1-\alpha+\alpha^2\right)^3},$$

$$\frac{1}{(-1+\alpha)^2\,\alpha\left(1-\alpha+\alpha^2\right)^3},\ -\frac{1}{(-1+\alpha)\,\alpha\left(1-\alpha+\alpha^2\right)^3},\ \frac{1}{(-1+\alpha)\,\alpha\left(1-\alpha+\alpha^2\right)^3},\ \frac{-1+\alpha}{\alpha\left(1-\alpha+\alpha^2\right)^3},$$

$$-\frac{(-1+\alpha)^2}{\alpha\left(1-\alpha+\alpha^2\right)^3},\ \frac{(-1+\alpha)^2}{\alpha\left(1-\alpha+\alpha^2\right)^3},\ -\frac{\alpha}{\left(1-\alpha+\alpha^2\right)^3},\ \frac{\alpha}{\left(1-\alpha+\alpha^2\right)^3},\ -\frac{\alpha}{(-1+\alpha)^2\left(1-\alpha+\alpha^2\right)^3},$$

$$\frac{\alpha}{(-1+\alpha)^2\left(1-\alpha+\alpha^2\right)^3},\ -\frac{\alpha}{(-1+\alpha)\left(1-\alpha+\alpha^2\right)^3},\ \frac{\alpha}{(-1+\alpha)\left(1-\alpha+\alpha^2\right)^3},\ -\frac{(-1+\alpha)\,\alpha}{\left(1-\alpha+\alpha^2\right)^3},$$

$$\frac{(-1+\alpha)\,\alpha}{\left(1-\alpha+\alpha^2\right)^3},\ -\frac{(-1+\alpha)^2\,\alpha}{\left(1-\alpha+\alpha^2\right)^3},\ \frac{(-1+\alpha)^2\,\alpha}{\left(1-\alpha+\alpha^2\right)^3},\ -\frac{\alpha^2}{\left(1-\alpha+\alpha^2\right)^3},\ \frac{\alpha^2}{\left(1-\alpha+\alpha^2\right)^3},\ -\frac{\alpha^2}{(-1+\alpha)^2\left(1-\alpha+\alpha^2\right)^3},$$

$$\frac{\alpha^2}{(-1+\alpha)^2\left(1-\alpha+\alpha^2\right)^3},\ -\frac{\alpha^2}{(-1+\alpha)\left(1-\alpha+\alpha^2\right)^3},\ \frac{\alpha^2}{(-1+\alpha)\left(1-\alpha+\alpha^2\right)^3},\ -\frac{(-1+\alpha)\,\alpha^2}{\left(1-\alpha+\alpha^2\right)^3},\ \frac{(-1+\alpha)\,\alpha^2}{\left(1-\alpha+\alpha^2\right)^3},$$

$$-\frac{(-1+\alpha)^2\,\alpha^2}{\left(1-\alpha+\alpha^2\right)^3},\ \frac{(-1+\alpha)^2\,\alpha^2}{\left(1-\alpha+\alpha^2\right)^3},\ -\frac{1}{\left(1-\alpha+\alpha^2\right)^2},\ \frac{1}{\left(1-\alpha+\alpha^2\right)^2},\ \frac{1-\alpha}{\left(1-\alpha+\alpha^2\right)^2},\ -\frac{1}{(-1+\alpha)\left(1-\alpha+\alpha^2\right)^2},$$

$$\frac{1}{(-1+\alpha)\left(1-\alpha+\alpha^2\right)^2},\ \frac{-1+\alpha}{\left(1-\alpha+\alpha^2\right)^2},\ -\frac{(-1+\alpha)^2}{\left(1-\alpha+\alpha^2\right)^2},\ \frac{(-1+\alpha)^2}{\left(1-\alpha+\alpha^2\right)^2},\ -\frac{1}{\alpha^2\left(1-\alpha+\alpha^2\right)^2},$$

$$\frac{1}{\alpha^2\left(1-\alpha+\alpha^2\right)^2},\ \frac{1-\alpha}{\alpha^2\left(1-\alpha+\alpha^2\right)^2},\ -\frac{1}{(-1+\alpha)\,\alpha^2\left(1-\alpha+\alpha^2\right)^2},\ \frac{1}{(-1+\alpha)\,\alpha^2\left(1-\alpha+\alpha^2\right)^2},$$

$$\frac{-1+\alpha}{\alpha^2\left(1-\alpha+\alpha^2\right)^2},\ -\frac{(-1+\alpha)^2}{\alpha^2\left(1-\alpha+\alpha^2\right)^2},\ \frac{(-1+\alpha)^2}{\alpha^2\left(1-\alpha+\alpha^2\right)^2},\ -\frac{1}{\alpha\left(1-\alpha+\alpha^2\right)^2},\ \frac{1}{\alpha\left(1-\alpha+\alpha^2\right)^2},$$

$$\frac{1-\alpha}{\alpha\left(1-\alpha+\alpha^2\right)^2},\ -\frac{1}{(-1+\alpha)\,\alpha\left(1-\alpha+\alpha^2\right)^2},\ \frac{1}{(-1+\alpha)\,\alpha\left(1-\alpha+\alpha^2\right)^2},\ \frac{-1+\alpha}{\alpha\left(1-\alpha+\alpha^2\right)^2},\ -\frac{(-1+\alpha)^2}{\alpha\left(1-\alpha+\alpha^2\right)^2},$$

$$\frac{(-1+\alpha)^2}{\alpha\left(1-\alpha+\alpha^2\right)^2},\ -\frac{\alpha}{\left(1-\alpha+\alpha^2\right)^2},\ \frac{\alpha}{\left(1-\alpha+\alpha^2\right)^2},\ -\frac{\alpha}{(-1+\alpha)\left(1-\alpha+\alpha^2\right)^2},\ \frac{\alpha}{(-1+\alpha)\left(1-\alpha+\alpha^2\right)^2},$$

$$-\frac{(-1+\alpha)\,\alpha}{\left(1-\alpha+\alpha^2\right)^2},\ \frac{(-1+\alpha)\,\alpha}{\left(1-\alpha+\alpha^2\right)^2},\ -\frac{(-1+\alpha)^2\,\alpha}{\left(1-\alpha+\alpha^2\right)^2},\ \frac{(-1+\alpha)^2\,\alpha}{\left(1-\alpha+\alpha^2\right)^2},\ -\frac{\alpha^2}{\left(1-\alpha+\alpha^2\right)^2},\ \frac{\alpha^2}{\left(1-\alpha+\alpha^2\right)^2},$$

$$-\frac{\alpha^2}{(-1+\alpha)\left(1-\alpha+\alpha^2\right)^2},\ \frac{\alpha^2}{(-1+\alpha)\left(1-\alpha+\alpha^2\right)^2},\ -\frac{(-1+\alpha)\,\alpha^2}{\left(1-\alpha+\alpha^2\right)^2},\ \frac{(-1+\alpha)\,\alpha^2}{\left(1-\alpha+\alpha^2\right)^2},\ -\frac{(-1+\alpha)^2\,\alpha^2}{\left(1-\alpha+\alpha^2\right)^2},$$

$$\frac{(-1+\alpha)^2\,\alpha^2}{\left(1-\alpha+\alpha^2\right)^2},\ \frac{1}{1-\alpha+\alpha^2},\ \frac{1-\alpha}{1-\alpha+\alpha^2},\ -\frac{1}{(-1+\alpha)^2\left(1-\alpha+\alpha^2\right)},\ \frac{1}{(-1+\alpha)^2\left(1-\alpha+\alpha^2\right)},$$

$$\frac{-1+\alpha}{1-\alpha+\alpha^2},\ -\frac{(-1+\alpha)^2}{1-\alpha+\alpha^2},\ \frac{(-1+\alpha)^2}{1-\alpha+\alpha^2},\ -\frac{1}{\alpha^2\left(1-\alpha+\alpha^2\right)},\ \frac{1}{\alpha^2\left(1-\alpha+\alpha^2\right)},\ \frac{1-\alpha}{\alpha^2\left(1-\alpha+\alpha^2\right)},$$

$$-\frac{1}{(-1+\alpha)^2\,\alpha^2\left(1-\alpha+\alpha^2\right)},\ \frac{1}{(-1+\alpha)^2\,\alpha^2\left(1-\alpha+\alpha^2\right)},\ \frac{-1+\alpha}{\alpha^2\left(1-\alpha+\alpha^2\right)},\ -\frac{(-1+\alpha)^2}{\alpha^2\left(1-\alpha+\alpha^2\right)},$$

$$\frac{(-1+\alpha)^2}{\alpha^2\left(1-\alpha+\alpha^2\right)},\ -\frac{1}{\alpha\left(1-\alpha+\alpha^2\right)},\ \frac{1}{\alpha\left(1-\alpha+\alpha^2\right)},\ \frac{1-\alpha}{\alpha\left(1-\alpha+\alpha^2\right)},\ -\frac{1}{(-1+\alpha)^2\,\alpha\left(1-\alpha+\alpha^2\right)},$$

$$\frac{1}{(-1+\alpha)^2\,\alpha\left(1-\alpha+\alpha^2\right)},\ \frac{-1+\alpha}{\alpha\left(1-\alpha+\alpha^2\right)},\ -\frac{(-1+\alpha)^2}{\alpha\left(1-\alpha+\alpha^2\right)},\ \frac{(-1+\alpha)^2}{\alpha\left(1-\alpha+\alpha^2\right)},\ -\frac{\alpha}{1-\alpha+\alpha^2},\ \frac{\alpha}{1-\alpha+\alpha^2},$$

$$-\frac{\alpha}{(-1+\alpha)^2\left(1-\alpha+\alpha^2\right)},\ \frac{\alpha}{(-1+\alpha)^2\left(1-\alpha+\alpha^2\right)},\ -\frac{(-1+\alpha)\,\alpha}{1-\alpha+\alpha^2},\ \frac{(-1+\alpha)\,\alpha}{1-\alpha+\alpha^2},\ -\frac{(-1+\alpha)^2\,\alpha}{1-\alpha+\alpha^2},$$

$$\frac{(-1+\alpha)^2\,\alpha}{1-\alpha+\alpha^2},\ -\frac{\alpha^2}{1-\alpha+\alpha^2},\ \frac{\alpha^2}{1-\alpha+\alpha^2},\ -\frac{\alpha^2}{(-1+\alpha)^2\left(1-\alpha+\alpha^2\right)},\ \frac{\alpha^2}{(-1+\alpha)^2\left(1-\alpha+\alpha^2\right)},$$

$$-\frac{(-1+\alpha)\,\alpha^2}{1-\alpha+\alpha^2},\ \frac{(-1+\alpha)\,\alpha^2}{1-\alpha+\alpha^2},\ -\frac{(-1+\alpha)^2\,\alpha^2}{1-\alpha+\alpha^2},\ \frac{(-1+\alpha)^2\,\alpha^2}{1-\alpha+\alpha^2},\ 1-\alpha+\alpha^2,\ \frac{1-\alpha+\alpha^2}{(-1+\alpha)^2},\ \frac{1-\alpha+\alpha^2}{-1+\alpha},$$

$$-(-1+\alpha)\left(1-\alpha+\alpha^2\right),\ (-1+\alpha)\left(1-\alpha+\alpha^2\right),\ -(-1+\alpha)^2\left(1-\alpha+\alpha^2\right),\ (-1+\alpha)^2\left(1-\alpha+\alpha^2\right),$$

$$\frac{1-\alpha+\alpha^2}{\alpha^2},\ \frac{1-\alpha+\alpha^2}{(-1+\alpha)^2\,\alpha^2},\ \frac{1-\alpha+\alpha^2}{(-1+\alpha)\,\alpha^2},\ -\frac{(-1+\alpha)\left(1-\alpha+\alpha^2\right)}{\alpha^2},\ \frac{(-1+\alpha)\left(1-\alpha+\alpha^2\right)}{\alpha^2},$$

$$-\frac{(-1+\alpha)^2\left(1-\alpha+\alpha^2\right)}{\alpha^2},\ \frac{(-1+\alpha)^2\left(1-\alpha+\alpha^2\right)}{\alpha^2},\ \frac{1-\alpha+\alpha^2}{\alpha},\ \frac{1-\alpha+\alpha^2}{(-1+\alpha)^2\,\alpha},\ \frac{1-\alpha+\alpha^2}{(-1+\alpha)\,\alpha},$$

$$-\frac{(-1+\alpha)\left(1-\alpha+\alpha^2\right)}{\alpha}, \frac{(-1+\alpha)\left(1-\alpha+\alpha^2\right)}{\alpha}, -\frac{(-1+\alpha)^2\left(1-\alpha+\alpha^2\right)}{\alpha}, \frac{(-1+\alpha)^2\left(1-\alpha+\alpha^2\right)}{\alpha},$$

$$-\alpha\left(1-\alpha+\alpha^2\right), \alpha\left(1-\alpha+\alpha^2\right), -\frac{\alpha\left(1-\alpha+\alpha^2\right)}{(-1+\alpha)^2}, \frac{\alpha\left(1-\alpha+\alpha^2\right)}{(-1+\alpha)^2}, -\frac{\alpha\left(1-\alpha+\alpha^2\right)}{-1+\alpha}, \frac{\alpha\left(1-\alpha+\alpha^2\right)}{-1+\alpha},$$

$$-(-1+\alpha)\,\alpha\left(1-\alpha+\alpha^2\right), (-1+\alpha)\,\alpha\left(1-\alpha+\alpha^2\right), -(-1+\alpha)^2\alpha\left(1-\alpha+\alpha^2\right), (-1+\alpha)^2\alpha\left(1-\alpha+\alpha^2\right),$$

$$-\alpha^2\left(1-\alpha+\alpha^2\right), \alpha^2\left(1-\alpha+\alpha^2\right), -\frac{\alpha^2\left(1-\alpha+\alpha^2\right)}{(-1+\alpha)^2}, \frac{\alpha^2\left(1-\alpha+\alpha^2\right)}{(-1+\alpha)^2}, -\frac{\alpha^2\left(1-\alpha+\alpha^2\right)}{-1+\alpha}, \frac{\alpha^2\left(1-\alpha+\alpha^2\right)}{-1+\alpha},$$

$$-(-1+\alpha)\,\alpha^2\left(1-\alpha+\alpha^2\right), (-1+\alpha)\,\alpha^2\left(1-\alpha+\alpha^2\right), -(-1+\alpha)^2\alpha^2\left(1-\alpha+\alpha^2\right), (-1+\alpha)^2\alpha^2\left(1-\alpha+\alpha^2\right),$$

$$-\left(1-\alpha+\alpha^2\right)^2, \left(1-\alpha+\alpha^2\right)^2, -\frac{\left(1-\alpha+\alpha^2\right)^2}{(-1+\alpha)^2}, \frac{\left(1-\alpha+\alpha^2\right)^2}{(-1+\alpha)^2}, -\frac{\left(1-\alpha+\alpha^2\right)^2}{-1+\alpha}, \frac{\left(1-\alpha+\alpha^2\right)^2}{-1+\alpha},$$

$$-(-1+\alpha)\left(1-\alpha+\alpha^2\right)^2, (-1+\alpha)\left(1-\alpha+\alpha^2\right)^2, -(-1+\alpha)^2\left(1-\alpha+\alpha^2\right)^2, (-1+\alpha)^2\left(1-\alpha+\alpha^2\right)^2,$$

$$-\frac{\left(1-\alpha+\alpha^2\right)^2}{\alpha^2}, \frac{\left(1-\alpha+\alpha^2\right)^2}{\alpha^2}, -\frac{\left(1-\alpha+\alpha^2\right)^2}{(-1+\alpha)^2\alpha^2}, \frac{\left(1-\alpha+\alpha^2\right)^2}{(-1+\alpha)^2\alpha^2}, -\frac{\left(1-\alpha+\alpha^2\right)^2}{(-1+\alpha)\,\alpha^2}, \frac{\left(1-\alpha+\alpha^2\right)^2}{(-1+\alpha)\,\alpha^2},$$

$$-\frac{(-1+\alpha)\left(1-\alpha+\alpha^2\right)^2}{\alpha^2}, \frac{(-1+\alpha)\left(1-\alpha+\alpha^2\right)^2}{\alpha^2}, -\frac{(-1+\alpha)^2\left(1-\alpha+\alpha^2\right)^2}{\alpha^2}, \frac{(-1+\alpha)^2\left(1-\alpha+\alpha^2\right)^2}{\alpha^2},$$

$$-\frac{\left(1-\alpha+\alpha^2\right)^2}{\alpha}, \frac{\left(1-\alpha+\alpha^2\right)^2}{\alpha}, -\frac{\left(1-\alpha+\alpha^2\right)^2}{(-1+\alpha)^2\alpha}, \frac{\left(1-\alpha+\alpha^2\right)^2}{(-1+\alpha)^2\alpha}, -\frac{\left(1-\alpha+\alpha^2\right)^2}{(-1+\alpha)\,\alpha}, \frac{\left(1-\alpha+\alpha^2\right)^2}{(-1+\alpha)\,\alpha},$$

$$-\frac{(-1+\alpha)\left(1-\alpha+\alpha^2\right)^2}{\alpha}, \frac{(-1+\alpha)\left(1-\alpha+\alpha^2\right)^2}{\alpha}, -\frac{(-1+\alpha)^2\left(1-\alpha+\alpha^2\right)^2}{\alpha}, \frac{(-1+\alpha)^2\left(1-\alpha+\alpha^2\right)^2}{\alpha},$$

$$-\alpha\left(1-\alpha+\alpha^2\right)^2, \alpha\left(1-\alpha+\alpha^2\right)^2, -\frac{\alpha\left(1-\alpha+\alpha^2\right)^2}{(-1+\alpha)^2}, \frac{\alpha\left(1-\alpha+\alpha^2\right)^2}{(-1+\alpha)^2}, -\frac{\alpha\left(1-\alpha+\alpha^2\right)^2}{-1+\alpha},$$

$$\frac{\alpha\left(1-\alpha+\alpha^2\right)^2}{-1+\alpha}, -(-1+\alpha)\,\alpha\left(1-\alpha+\alpha^2\right)^2, (-1+\alpha)\,\alpha\left(1-\alpha+\alpha^2\right)^2, -(-1+\alpha)^2\alpha\left(1-\alpha+\alpha^2\right)^2,$$

$$(-1+\alpha)^2\alpha\left(1-\alpha+\alpha^2\right)^2, -\alpha^2\left(1-\alpha+\alpha^2\right)^2, \alpha^2\left(1-\alpha+\alpha^2\right)^2, -\frac{\alpha^2\left(1-\alpha+\alpha^2\right)^2}{(-1+\alpha)^2}, \frac{\alpha^2\left(1-\alpha+\alpha^2\right)^2}{(-1+\alpha)^2},$$

$$-\frac{\alpha^2\left(1-\alpha+\alpha^2\right)^2}{-1+\alpha}, \frac{\alpha^2\left(1-\alpha+\alpha^2\right)^2}{-1+\alpha}, -(-1+\alpha)\,\alpha^2\left(1-\alpha+\alpha^2\right)^2, (-1+\alpha)\,\alpha^2\left(1-\alpha+\alpha^2\right)^2,$$

$$-(-1+\alpha)^2\alpha^2\left(1-\alpha+\alpha^2\right)^2, (-1+\alpha)^2\alpha^2\left(1-\alpha+\alpha^2\right)^2, -\left(1-\alpha+\alpha^2\right)^3, \left(1-\alpha+\alpha^2\right)^3, -\frac{\left(1-\alpha+\alpha^2\right)^3}{(-1+\alpha)^2},$$

$$\frac{\left(1-\alpha+\alpha^2\right)^3}{(-1+\alpha)^2}, -\frac{\left(1-\alpha+\alpha^2\right)^3}{-1+\alpha}, \frac{\left(1-\alpha+\alpha^2\right)^3}{-1+\alpha}, -(-1+\alpha)\left(1-\alpha+\alpha^2\right)^3, (-1+\alpha)\left(1-\alpha+\alpha^2\right)^3,$$

$$-(-1+\alpha)^2\left(1-\alpha+\alpha^2\right)^3, (-1+\alpha)^2\left(1-\alpha+\alpha^2\right)^3, -\frac{\left(1-\alpha+\alpha^2\right)^3}{\alpha^2}, \frac{\left(1-\alpha+\alpha^2\right)^3}{\alpha^2}, -\frac{\left(1-\alpha+\alpha^2\right)^3}{(-1+\alpha)^2\alpha^2},$$

$$\frac{\left(1-\alpha+\alpha^2\right)^3}{(-1+\alpha)^2\alpha^2}, -\frac{\left(1-\alpha+\alpha^2\right)^3}{(-1+\alpha)\,\alpha^2}, \frac{\left(1-\alpha+\alpha^2\right)^3}{(-1+\alpha)\,\alpha^2}, -\frac{(-1+\alpha)\left(1-\alpha+\alpha^2\right)^3}{\alpha^2}, \frac{(-1+\alpha)\left(1-\alpha+\alpha^2\right)^3}{\alpha^2},$$

$$-\frac{(-1+\alpha)^2\left(1-\alpha+\alpha^2\right)^3}{\alpha^2},\ \frac{(-1+\alpha)^2\left(1-\alpha+\alpha^2\right)^3}{\alpha^2},\ -\frac{\left(1-\alpha+\alpha^2\right)^3}{\alpha},\ \frac{\left(1-\alpha+\alpha^2\right)^3}{\alpha},\ -\frac{\left(1-\alpha+\alpha^2\right)^3}{(-1+\alpha)^2\alpha},$$

$$\frac{\left(1-\alpha+\alpha^2\right)^3}{(-1+\alpha)^2\alpha},\ -\frac{\left(1-\alpha+\alpha^2\right)^3}{(-1+\alpha)\alpha},\ \frac{\left(1-\alpha+\alpha^2\right)^3}{(-1+\alpha)\alpha},\ -\frac{(-1+\alpha)\left(1-\alpha+\alpha^2\right)^3}{\alpha},\ \frac{(-1+\alpha)\left(1-\alpha+\alpha^2\right)^3}{\alpha},$$

$$-\frac{(-1+\alpha)^2\left(1-\alpha+\alpha^2\right)^3}{\alpha},\ \frac{(-1+\alpha)^2\left(1-\alpha+\alpha^2\right)^3}{\alpha},\ -\alpha\left(1-\alpha+\alpha^2\right)^3,\ \alpha\left(1-\alpha+\alpha^2\right)^3,\ -\frac{\alpha\left(1-\alpha+\alpha^2\right)^3}{(-1+\alpha)^2},$$

$$\frac{\alpha\left(1-\alpha+\alpha^2\right)^3}{(-1+\alpha)^2},\ -\frac{\alpha\left(1-\alpha+\alpha^2\right)^3}{-1+\alpha},\ \frac{\alpha\left(1-\alpha+\alpha^2\right)^3}{-1+\alpha},\ -(-1+\alpha)\alpha\left(1-\alpha+\alpha^2\right)^3,\ (-1+\alpha)\alpha\left(1-\alpha+\alpha^2\right)^3,$$

$$-(-1+\alpha)^2\alpha\left(1-\alpha+\alpha^2\right)^3,\ (-1+\alpha)^2\alpha\left(1-\alpha+\alpha^2\right)^3,\ -\alpha^2\left(1-\alpha+\alpha^2\right)^3,\ \alpha^2\left(1-\alpha+\alpha^2\right)^3,$$

$$-\frac{\alpha^2\left(1-\alpha+\alpha^2\right)^3}{(-1+\alpha)^2},\ \frac{\alpha^2\left(1-\alpha+\alpha^2\right)^3}{(-1+\alpha)^2},\ -\frac{\alpha^2\left(1-\alpha+\alpha^2\right)^3}{-1+\alpha},\ \frac{\alpha^2\left(1-\alpha+\alpha^2\right)^3}{-1+\alpha},\ -(-1+\alpha)\alpha^2\left(1-\alpha+\alpha^2\right)^3,$$

$$(-1+\alpha)\alpha^2\left(1-\alpha+\alpha^2\right)^3,\ -(-1+\alpha)^2\alpha^2\left(1-\alpha+\alpha^2\right)^3,\ (-1+\alpha)^2\alpha^2\left(1-\alpha+\alpha^2\right)^3,$$

$$-\frac{1}{\left(-1+2\alpha-2\alpha^2+\alpha^3\right)^2},\ \frac{1}{\left(-1+2\alpha-2\alpha^2+\alpha^3\right)^2},\ -\frac{1}{\alpha^2\left(-1+2\alpha-2\alpha^2+\alpha^3\right)^2},\ \frac{1}{\alpha^2\left(-1+2\alpha-2\alpha^2+\alpha^3\right)^2},$$

$$-\frac{1}{\alpha\left(-1+2\alpha-2\alpha^2+\alpha^3\right)^2},\ \frac{1}{\alpha\left(-1+2\alpha-2\alpha^2+\alpha^3\right)^2},\ -\frac{\alpha}{\left(-1+2\alpha-2\alpha^2+\alpha^3\right)^2},\ \frac{\alpha}{\left(-1+2\alpha-2\alpha^2+\alpha^3\right)^2},$$

$$-\frac{\alpha^2}{\left(-1+2\alpha-2\alpha^2+\alpha^3\right)^2},\ \frac{\alpha^2}{\left(-1+2\alpha-2\alpha^2+\alpha^3\right)^2},\ -\frac{1}{-1+2\alpha-2\alpha^2+\alpha^3},\ \frac{1}{-1+2\alpha-2\alpha^2+\alpha^3},$$

$$-\frac{1}{\alpha^2\left(-1+2\alpha-2\alpha^2+\alpha^3\right)},\ \frac{1}{\alpha^2\left(-1+2\alpha-2\alpha^2+\alpha^3\right)},\ -\frac{1}{\alpha\left(-1+2\alpha-2\alpha^2+\alpha^3\right)},\ \frac{1}{\alpha\left(-1+2\alpha-2\alpha^2+\alpha^3\right)},$$

$$-\frac{\alpha}{-1+2\alpha-2\alpha^2+\alpha^3},\ \frac{\alpha}{-1+2\alpha-2\alpha^2+\alpha^3},\ -\frac{\alpha^2}{-1+2\alpha-2\alpha^2+\alpha^3},\ \frac{\alpha^2}{-1+2\alpha-2\alpha^2+\alpha^3}\Big\}$$

```
Length[%]
```

```
351
```

To find out which elements are fundamental we map the candidates to GF($p$) for some large $p$, and choose the image of $\alpha$ so that this map is injective. Then we select the elements $x$ for which $1-x$ is also in the list. To avoid dealing with fractions we reconstruct the list using the PowerMod function:

```
p = Prime[100 000]
```

```
1 299 709
```

```
aIm = 5;
```

```
H3genIms = Mod[H3gens /. α → aIm, p]
```

```
{1 299 708, 5, 1 299 705, 21}
```

```
candidateH3funIms =
  {0} ⋃ (Mod[Times @@ PowerMod[H3genIms, #, p], p] & /@ candidateH3funExps);
```

This map is indeed injective (the union symbol removes duplicates)

```
Length[%]
```

```
351
```

Select the elements $x$ for which $1 - x$ is also in the list:

```
H3funIms = candidateH3funIms ⋂ Mod[1 - candidateH3funIms, p]
```

```
{0, 1, 5, 21, 123 783, 259 942, 259 946, 311 931, 324 922, 324 927,
 495 128, 568 624, 584 869, 618 910, 680 800, 714 841, 731 086, 804 582,
 974 783, 974 788, 987 779, 1 039 764, 1 039 768, 1 175 927, 1 299 689, 1 299 705}
```

```
Length[%]
```

```
26
```

This list of remaining candidates necessarily contains the images of all known fundamental elements. Since the number of remaining elements equals the number of known fundamental elements, we conclude that those are, indeed, all. This completes the proof of Lemma 19.

## The automorphism group

In this section we verify the following statement:

**Lemma 20.** The group of partial-field automorphisms of $\mathbb{H}_3$ is isomorphic to $S_3$, the symmetric group on three elements.

The partial-field automorphism group of $\mathbb{H}_3$ is determined uniquely by the image of $\alpha$. This is easy to see since each automorphism maps 0 to 0, 1 to 1, and $-1$ to $-1$, so the real number line is fixed. Note that $\alpha$ has to be mapped to a fundamental element, so the group will be finite.

Note that we only have to check if $\left\{1, \alpha, \alpha^2 - \alpha + 1, \frac{\alpha^2}{\alpha-1}, \frac{-\alpha}{(\alpha-1)^2}\right\}$ map to fundamental elements: the remainder follows automatically.

```
isH3fun[x_] := occursInList[H3funs, x];
areH3funs[x_List] := And @@ (isH3fun /@ x);
```

```
candidateAuts = Table[Together[{α, α^2 - α + 1, α^2/(α - 1), -α/(α - 1)^2} /. α → nonzerooneH3funs[[x]]],
   {x, 1, Length[nonzerooneH3funs]}]
```

$$\Big\{\Big\{\frac{1}{1-\alpha}, \frac{1-\alpha+\alpha^2}{(-1+\alpha)^2}, -\frac{1}{(-1+\alpha)\,\alpha}, \frac{-1+\alpha}{\alpha^2}\Big\}, \Big\{1-\alpha, 1-\alpha+\alpha^2, -\frac{(-1+\alpha)^2}{\alpha}, \frac{-1+\alpha}{\alpha^2}\Big\},$$

$$\Big\{\frac{-1+\alpha}{\alpha^2}, \frac{1-2\,\alpha+2\,\alpha^2-\alpha^3+\alpha^4}{\alpha^4}, -\frac{(-1+\alpha)^2}{\alpha^2\,(1-\alpha+\alpha^2)}, -\frac{(-1+\alpha)\,\alpha^2}{(1-\alpha+\alpha^2)^2}\Big\},$$

$$\Big\{\frac{1}{\alpha}, \frac{1-\alpha+\alpha^2}{\alpha^2}, -\frac{1}{(-1+\alpha)\,\alpha}, -\frac{\alpha}{(-1+\alpha)^2}\Big\},$$

$$\Big\{-\frac{1}{(-1+\alpha)\,\alpha}, \frac{1-\alpha+2\,\alpha^2-2\,\alpha^3+\alpha^4}{(-1+\alpha)^2\,\alpha^2}, -\frac{1}{\alpha\,(-1+2\,\alpha-2\,\alpha^2+\alpha^3)}, \frac{(-1+\alpha)\,\alpha}{(1-\alpha+\alpha^2)^2}\Big\},$$

$$\Big\{\frac{-1+\alpha}{\alpha}, \frac{1-\alpha+\alpha^2}{\alpha^2}, -\frac{(-1+\alpha)^2}{\alpha}, -(-1+\alpha)\,\alpha\Big\},$$

$$\Big\{-\frac{(-1+\alpha)^2}{\alpha}, \frac{1-3\,\alpha+5\,\alpha^2-3\,\alpha^3+\alpha^4}{\alpha^2}, -\frac{(-1+\alpha)^4}{\alpha\,(1-\alpha+\alpha^2)}, \frac{(-1+\alpha)^2\,\alpha}{(1-\alpha+\alpha^2)^2}\Big\},$$

$$\Big\{\alpha, 1-\alpha+\alpha^2, \frac{\alpha^2}{-1+\alpha}, -\frac{\alpha}{(-1+\alpha)^2}\Big\},$$

$$\Big\{-\frac{\alpha}{(-1+\alpha)^2}, \frac{1-3\,\alpha+5\,\alpha^2-3\,\alpha^3+\alpha^4}{(-1+\alpha)^4}, -\frac{\alpha^2}{(-1+\alpha)^2\,(1-\alpha+\alpha^2)}, \frac{(-1+\alpha)^2\,\alpha}{(1-\alpha+\alpha^2)^2}\Big\},$$

$$\Big\{\frac{\alpha}{-1+\alpha}, \frac{1-\alpha+\alpha^2}{(-1+\alpha)^2}, \frac{\alpha^2}{-1+\alpha}, -(-1+\alpha)\,\alpha\Big\},$$

$$\Big\{-(-1+\alpha)\,\alpha, 1-\alpha+2\,\alpha^2-2\,\alpha^3+\alpha^4, -\frac{(-1+\alpha)^2\,\alpha^2}{1-\alpha+\alpha^2}, \frac{(-1+\alpha)\,\alpha}{(1-\alpha+\alpha^2)^2}\Big\},$$

$$\Big\{\frac{\alpha^2}{-1+\alpha}, \frac{1-2\,\alpha+2\,\alpha^2-\alpha^3+\alpha^4}{(-1+\alpha)^2}, \frac{\alpha^4}{(-1+\alpha)\,(1-\alpha+\alpha^2)}, -\frac{(-1+\alpha)\,\alpha^2}{(1-\alpha+\alpha^2)^2}\Big\},$$

$$\Big\{\frac{-1+\alpha-\alpha^2}{-1+\alpha}, \frac{1-2\,\alpha+2\,\alpha^2-\alpha^3+\alpha^4}{(-1+\alpha)^2}, -\frac{(1-\alpha+\alpha^2)^2}{(-1+\alpha)\,\alpha^2}, \frac{(-1+\alpha)\,(1-\alpha+\alpha^2)}{\alpha^4}\Big\},$$

$$\Big\{\frac{1}{1-\alpha+\alpha^2}, \frac{1-\alpha+2\,\alpha^2-2\,\alpha^3+\alpha^4}{(1-\alpha+\alpha^2)^2}, -\frac{1}{\alpha\,(-1+2\,\alpha-2\,\alpha^2+\alpha^3)}, \frac{-1+\alpha-\alpha^2}{(-1+\alpha)^2\,\alpha^2}\Big\},$$

$$\Big\{\frac{1-\alpha}{1-\alpha+\alpha^2}, \frac{1-2\,\alpha+2\,\alpha^2-\alpha^3+\alpha^4}{(1-\alpha+\alpha^2)^2}, -\frac{(-1+\alpha)^2}{\alpha^2\,(1-\alpha+\alpha^2)}, \frac{(-1+\alpha)\,(1-\alpha+\alpha^2)}{\alpha^4}\Big\},$$

$$\left\{\frac{\alpha}{1-\alpha+\alpha^2}, \frac{1-3\alpha+5\alpha^2-3\alpha^3+\alpha^4}{(1-\alpha+\alpha^2)^2}, -\frac{\alpha^2}{(-1+\alpha)^2(1-\alpha+\alpha^2)}, -\frac{\alpha(1-\alpha+\alpha^2)}{(-1+\alpha)^4}\right\},$$

$$\left\{\frac{(-1+\alpha)\alpha}{1-\alpha+\alpha^2}, \frac{1-\alpha+2\alpha^2-2\alpha^3+\alpha^4}{(1-\alpha+\alpha^2)^2}, -\frac{(-1+\alpha)^2\alpha^2}{1-\alpha+\alpha^2}, -(-1+\alpha)\alpha(1-\alpha+\alpha^2)\right\},$$

$$\left\{\frac{\alpha^2}{1-\alpha+\alpha^2}, \frac{1-2\alpha+2\alpha^2-\alpha^3+\alpha^4}{(1-\alpha+\alpha^2)^2}, \frac{\alpha^4}{(-1+\alpha)(1-\alpha+\alpha^2)}, -\frac{\alpha^2(1-\alpha+\alpha^2)}{(-1+\alpha)^2}\right\},$$

$$\left\{\frac{1-2\alpha+\alpha^2}{1-\alpha+\alpha^2}, \frac{1-3\alpha+5\alpha^2-3\alpha^3+\alpha^4}{(1-\alpha+\alpha^2)^2}, -\frac{(1-2\alpha+\alpha^2)^2}{\alpha(1-\alpha+\alpha^2)}, -\frac{(1-2\alpha+\alpha^2)(1-\alpha+\alpha^2)}{\alpha^2}\right\},$$

$$\left\{1-\alpha+\alpha^2, 1-\alpha+2\alpha^2-2\alpha^3+\alpha^4, \frac{(1-\alpha+\alpha^2)^2}{(-1+\alpha)\alpha}, \frac{-1+\alpha-\alpha^2}{(-1+\alpha)^2\alpha^2}\right\},$$

$$\left\{\frac{1-\alpha+\alpha^2}{(-1+\alpha)^2}, \frac{1-3\alpha+5\alpha^2-3\alpha^3+\alpha^4}{(-1+\alpha)^4}, \frac{(1-\alpha+\alpha^2)^2}{(-1+\alpha)^2\alpha}, -\frac{(-1+\alpha)^2(1-\alpha+\alpha^2)}{\alpha^2}\right\},$$

$$\left\{\frac{1-\alpha+\alpha^2}{\alpha^2}, \frac{1-2\alpha+2\alpha^2-\alpha^3+\alpha^4}{\alpha^4}, -\frac{(1-\alpha+\alpha^2)^2}{(-1+\alpha)\alpha^2}, -\frac{\alpha^2(1-\alpha+\alpha^2)}{(-1+\alpha)^2}\right\},$$

$$\left\{\frac{1-\alpha+\alpha^2}{\alpha}, \frac{1-3\alpha+5\alpha^2-3\alpha^3+\alpha^4}{\alpha^2}, \frac{(1-\alpha+\alpha^2)^2}{(-1+\alpha)^2\alpha}, -\frac{\alpha(1-\alpha+\alpha^2)}{(-1+\alpha)^4}\right\},$$

$$\left\{\frac{1-\alpha+\alpha^2}{(-1+\alpha)\alpha}, \frac{1-\alpha+2\alpha^2-2\alpha^3+\alpha^4}{(-1+\alpha)^2\alpha^2}, \frac{(1-\alpha+\alpha^2)^2}{(-1+\alpha)\alpha}, -(-1+\alpha)\alpha(1-\alpha+\alpha^2)\right\}\}$$

```
Select[candidateAuts, areH3funs]
```

$$\left\{\left\{\frac{1}{1-\alpha}, \frac{1-\alpha+\alpha^2}{(-1+\alpha)^2}, -\frac{1}{(-1+\alpha)\alpha}, \frac{-1+\alpha}{\alpha^2}\right\}, \left\{1-\alpha, 1-\alpha+\alpha^2, -\frac{(-1+\alpha)^2}{\alpha}, \frac{-1+\alpha}{\alpha^2}\right\},\right.$$

$$\left\{\frac{1}{\alpha}, \frac{1-\alpha+\alpha^2}{\alpha^2}, -\frac{1}{(-1+\alpha)\alpha}, -\frac{\alpha}{(-1+\alpha)^2}\right\}, \left\{\frac{-1+\alpha}{\alpha}, \frac{1-\alpha+\alpha^2}{\alpha^2}, -\frac{(-1+\alpha)^2}{\alpha}, -(-1+\alpha)\alpha\right\},$$

$$\left.\left\{\alpha, 1-\alpha+\alpha^2, \frac{\alpha^2}{-1+\alpha}, -\frac{\alpha}{(-1+\alpha)^2}\right\}, \left\{\frac{\alpha}{-1+\alpha}, \frac{1-\alpha+\alpha^2}{(-1+\alpha)^2}, \frac{\alpha^2}{-1+\alpha}, -(-1+\alpha)\alpha\right\}\right\}$$

Conclusion: the automorphism group of $\mathbb{H}_3$ has 6 elements and is isomorphic to $S_3$, the symmetric group on 3 symbols. Indeed, the partial-field homomorphism $\phi 3$ defined above induces a $1-1$ correspondence between the automorphisms of $\mathbb{H}_3$ and coordinate permutations of $\mathrm{GF}(5)\times\mathrm{GF}(5)\times\mathrm{GF}(5)$. Each automorphism maps $\alpha$ to one of its associates. This completes the proof of Lemma 20.

## Representations of $U_{2,5}$

We prove the following fact:

**Lemma 21.** Let $M$ be a 3-connected, $\mathbb{H}_3$-representable matroid. If $M$ has a $U_{2,5}$- or $U_{3,5}$-minor then $M$ has at least three inequivalent representations over GF(5).

To prove this it suffices to show that each $\mathbb{H}_3$-representation of $U_{2,5}$ gives rise to three inequivalent representations over GF(5). Since any $\mathbb{H}_3$-representation matrix of $M$ must contain one of these as a minor, the result then follows. First we enumerate all such representations. By normalizing and suppressing the identity matrix at the front, we see

```
Clear[p];
```

```
A = {{1, 1, 1}, {1, p, q}}; MatrixForm[A]
```

$$\begin{pmatrix} 1 & 1 & 1 \\ 1 & p & q \end{pmatrix}$$

This is an $\mathbb{H}_3$-matrix representing $U_{2,5}$ if and only if $p$, $q$, $\frac{p}{q}$ are fundamental elements, $p$ and $q$ are not equal to 0 or 1, and $p \neq q$. Moreover, two such matrices are equivalent if and only if they are equal.

Create a table of all candidate pairs $p$, $q$. We suppress the output:

```
candidatepqPairs = Flatten[Table[{nonzerooneH3funs[[x]], nonzerooneH3funs[[y]]},
     {x, 1, Length[nonzerooneH3funs]}, {y, 1, x - 1}], 1] ⋃
   Flatten[Table[{nonzerooneH3funs[[x]], nonzerooneH3funs[[y]]},
     {x, 1, Length[nonzerooneH3funs]}, {y, x + 1, Length[nonzerooneH3funs]}], 1];
```

Filter out those for which $p/q$ is not fundamental or equal to 1:

```
pqPairs = Select[candidatepqPairs, occursInList[nonzerooneH3funs, #[[1]]/#[[2]]] &]
```

$$\Big\{\Big\{\frac{1}{1-\alpha}, \frac{1}{(1-\alpha)\alpha}\Big\}, \Big\{\frac{1}{1-\alpha}, \frac{-1+\alpha}{\alpha}\Big\}, \Big\{\frac{1}{1-\alpha}, \alpha\Big\}, \Big\{\frac{1}{1-\alpha}, -\frac{\alpha}{(-1+\alpha)^2}\Big\}, \Big\{\frac{1}{1-\alpha}, \frac{-1+\alpha-\alpha^2}{-1+\alpha}\Big\},$$
$$\Big\{\frac{1}{1-\alpha}, \frac{1}{1-\alpha+\alpha^2}\Big\}, \Big\{\frac{1}{1-\alpha}, \frac{1-\alpha}{1-\alpha+\alpha^2}\Big\}, \Big\{\frac{1}{1-\alpha}, \frac{1-\alpha+\alpha^2}{(-1+\alpha)^2}\Big\}, \Big\{1-\alpha, \frac{1}{\alpha}\Big\}, \Big\{1-\alpha, -\frac{(-1+\alpha)^2}{\alpha}\Big\},$$
$$\Big\{1-\alpha, \frac{\alpha}{-1+\alpha}\Big\}, \{1-\alpha, -(-1+\alpha)\alpha\}, \Big\{1-\alpha, \frac{-1+\alpha-\alpha^2}{-1+\alpha}\Big\}, \Big\{1-\alpha, \frac{1-\alpha}{1-\alpha+\alpha^2}\Big\}, \Big\{1-\alpha, \frac{1-2\alpha+\alpha^2}{1-\alpha+\alpha^2}\Big\},$$
$$\{1-\alpha, 1-\alpha+\alpha^2\}, \Big\{\frac{-1+\alpha}{\alpha^2}, \frac{1}{\alpha}\Big\}, \Big\{\frac{-1+\alpha}{\alpha^2}, \frac{1}{(1-\alpha)\alpha}\Big\}, \Big\{\frac{-1+\alpha}{\alpha^2}, \frac{-1+\alpha}{\alpha}\Big\}, \Big\{\frac{-1+\alpha}{\alpha^2}, -\frac{(-1+\alpha)^2}{\alpha}\Big\},$$
$$\Big\{\frac{1}{\alpha}, 1-\alpha\Big\}, \Big\{\frac{1}{\alpha}, \frac{-1+\alpha}{\alpha^2}\Big\}, \Big\{\frac{1}{\alpha}, \frac{1}{(1-\alpha)\alpha}\Big\}, \Big\{\frac{1}{\alpha}, \frac{\alpha}{-1+\alpha}\Big\}, \Big\{\frac{1}{\alpha}, \frac{1}{1-\alpha+\alpha^2}\Big\}, \Big\{\frac{1}{\alpha}, \frac{\alpha}{1-\alpha+\alpha^2}\Big\},$$
$$\Big\{\frac{1}{\alpha}, \frac{1-\alpha+\alpha^2}{\alpha^2}\Big\}, \Big\{\frac{1}{\alpha}, \frac{1-\alpha+\alpha^2}{\alpha}\Big\}, \Big\{\frac{1}{(1-\alpha)\alpha}, \frac{1}{1-\alpha}\Big\}, \Big\{\frac{1}{(1-\alpha)\alpha}, \frac{-1+\alpha}{\alpha^2}\Big\}, \Big\{\frac{1}{(1-\alpha)\alpha}, \frac{1}{\alpha}\Big\},$$

$$\left\{\frac{1}{(1-\alpha)\,\alpha}, -\frac{\alpha}{(-1+\alpha)^2}\right\}, \left\{\frac{-1+\alpha}{\alpha}, \frac{1}{1-\alpha}\right\}, \left\{\frac{-1+\alpha}{\alpha}, \frac{-1+\alpha}{\alpha^2}\right\}, \left\{\frac{-1+\alpha}{\alpha}, -\frac{(-1+\alpha)^2}{\alpha}\right\},$$

$$\left\{\frac{-1+\alpha}{\alpha}, \alpha\right\}, \left\{\frac{-1+\alpha}{\alpha}, \frac{(-1+\alpha)\,\alpha}{1-\alpha+\alpha^2}\right\}, \left\{\frac{-1+\alpha}{\alpha}, \frac{1-2\,\alpha+\alpha^2}{1-\alpha+\alpha^2}\right\}, \left\{\frac{-1+\alpha}{\alpha}, \frac{1-\alpha+\alpha^2}{\alpha^2}\right\},$$

$$\left\{\frac{-1+\alpha}{\alpha}, \frac{1-\alpha+\alpha^2}{(-1+\alpha)\,\alpha}\right\}, \left\{-\frac{(-1+\alpha)^2}{\alpha}, 1-\alpha\right\}, \left\{-\frac{(-1+\alpha)^2}{\alpha}, \frac{-1+\alpha}{\alpha^2}\right\}, \left\{-\frac{(-1+\alpha)^2}{\alpha}, \frac{-1+\alpha}{\alpha}\right\},$$

$$\left\{-\frac{(-1+\alpha)^2}{\alpha}, -(-1+\alpha)\,\alpha\right\}, \left\{\alpha, \frac{1}{1-\alpha}\right\}, \left\{\alpha, \frac{-1+\alpha}{\alpha}\right\}, \{\alpha, -(-1+\alpha)\,\alpha\}, \left\{\alpha, \frac{\alpha^2}{-1+\alpha}\right\},$$

$$\left\{\alpha, \frac{\alpha}{1-\alpha+\alpha^2}\right\}, \left\{\alpha, \frac{\alpha^2}{1-\alpha+\alpha^2}\right\}, \left\{\alpha, 1-\alpha+\alpha^2\right\}, \left\{\alpha, \frac{1-\alpha+\alpha^2}{\alpha}\right\}, \left\{-\frac{\alpha}{(-1+\alpha)^2}, \frac{1}{1-\alpha}\right\},$$

$$\left\{-\frac{\alpha}{(-1+\alpha)^2}, \frac{1}{(1-\alpha)\,\alpha}\right\}, \left\{-\frac{\alpha}{(-1+\alpha)^2}, \frac{\alpha}{-1+\alpha}\right\}, \left\{-\frac{\alpha}{(-1+\alpha)^2}, \frac{\alpha^2}{-1+\alpha}\right\}, \left\{\frac{\alpha}{-1+\alpha}, 1-\alpha\right\},$$

$$\left\{\frac{\alpha}{-1+\alpha}, \frac{1}{\alpha}\right\}, \left\{\frac{\alpha}{-1+\alpha}, -\frac{\alpha}{(-1+\alpha)^2}\right\}, \left\{\frac{\alpha}{-1+\alpha}, \frac{\alpha^2}{-1+\alpha}\right\}, \left\{\frac{\alpha}{-1+\alpha}, \frac{(-1+\alpha)\,\alpha}{1-\alpha+\alpha^2}\right\},$$

$$\left\{\frac{\alpha}{-1+\alpha}, \frac{\alpha^2}{1-\alpha+\alpha^2}\right\}, \left\{\frac{\alpha}{-1+\alpha}, \frac{1-\alpha+\alpha^2}{(-1+\alpha)^2}\right\}, \left\{\frac{\alpha}{-1+\alpha}, \frac{1-\alpha+\alpha^2}{(-1+\alpha)\,\alpha}\right\}, \{-(-1+\alpha)\,\alpha, 1-\alpha\},$$

$$\left\{-(-1+\alpha)\,\alpha, -\frac{(-1+\alpha)^2}{\alpha}\right\}, \{-(-1+\alpha)\,\alpha, \alpha\}, \left\{-(-1+\alpha)\,\alpha, \frac{\alpha^2}{-1+\alpha}\right\}, \left\{\frac{\alpha^2}{-1+\alpha}, \alpha\right\},$$

$$\left\{\frac{\alpha^2}{-1+\alpha}, -\frac{\alpha}{(-1+\alpha)^2}\right\}, \left\{\frac{\alpha^2}{-1+\alpha}, \frac{\alpha}{-1+\alpha}\right\}, \left\{\frac{\alpha^2}{-1+\alpha}, -(-1+\alpha)\,\alpha\right\}, \left\{\frac{-1+\alpha-\alpha^2}{-1+\alpha}, \frac{1}{1-\alpha}\right\},$$

$$\left\{\frac{-1+\alpha-\alpha^2}{-1+\alpha}, 1-\alpha\right\}, \left\{\frac{-1+\alpha-\alpha^2}{-1+\alpha}, 1-\alpha+\alpha^2\right\}, \left\{\frac{-1+\alpha-\alpha^2}{-1+\alpha}, \frac{1-\alpha+\alpha^2}{(-1+\alpha)^2}\right\}, \left\{\frac{1}{1-\alpha+\alpha^2}, \frac{1}{1-\alpha}\right\},$$

$$\left\{\frac{1}{1-\alpha+\alpha^2}, \frac{1}{\alpha}\right\}, \left\{\frac{1}{1-\alpha+\alpha^2}, \frac{1-\alpha}{1-\alpha+\alpha^2}\right\}, \left\{\frac{1}{1-\alpha+\alpha^2}, \frac{\alpha}{1-\alpha+\alpha^2}\right\}, \left\{\frac{1-\alpha}{1-\alpha+\alpha^2}, \frac{1}{1-\alpha}\right\},$$

$$\left\{\frac{1-\alpha}{1-\alpha+\alpha^2}, 1-\alpha\right\}, \left\{\frac{1-\alpha}{1-\alpha+\alpha^2}, \frac{1}{1-\alpha+\alpha^2}\right\}, \left\{\frac{1-\alpha}{1-\alpha+\alpha^2}, \frac{1-2\,\alpha+\alpha^2}{1-\alpha+\alpha^2}\right\}, \left\{\frac{\alpha}{1-\alpha+\alpha^2}, \frac{1}{\alpha}\right\},$$

$$\left\{\frac{\alpha}{1-\alpha+\alpha^2}, \alpha\right\}, \left\{\frac{\alpha}{1-\alpha+\alpha^2}, \frac{1}{1-\alpha+\alpha^2}\right\}, \left\{\frac{\alpha}{1-\alpha+\alpha^2}, \frac{\alpha^2}{1-\alpha+\alpha^2}\right\}, \left\{\frac{(-1+\alpha)\,\alpha}{1-\alpha+\alpha^2}, \frac{-1+\alpha}{\alpha}\right\},$$

$$\left\{\frac{(-1+\alpha)\,\alpha}{1-\alpha+\alpha^2}, \frac{\alpha}{-1+\alpha}\right\}, \left\{\frac{(-1+\alpha)\,\alpha}{1-\alpha+\alpha^2}, \frac{\alpha^2}{1-\alpha+\alpha^2}\right\}, \left\{\frac{(-1+\alpha)\,\alpha}{1-\alpha+\alpha^2}, \frac{1-2\,\alpha+\alpha^2}{1-\alpha+\alpha^2}\right\}, \left\{\frac{\alpha^2}{1-\alpha+\alpha^2}, \alpha\right\},$$

$$\left\{\frac{\alpha^2}{1-\alpha+\alpha^2}, \frac{\alpha}{-1+\alpha}\right\}, \left\{\frac{\alpha^2}{1-\alpha+\alpha^2}, \frac{\alpha}{1-\alpha+\alpha^2}\right\}, \left\{\frac{\alpha^2}{1-\alpha+\alpha^2}, \frac{(-1+\alpha)\,\alpha}{1-\alpha+\alpha^2}\right\}, \left\{\frac{1-2\,\alpha+\alpha^2}{1-\alpha+\alpha^2}, 1-\alpha\right\},$$

$$\left\{\frac{1-2\,\alpha+\alpha^2}{1-\alpha+\alpha^2}, \frac{-1+\alpha}{\alpha}\right\}, \left\{\frac{1-2\,\alpha+\alpha^2}{1-\alpha+\alpha^2}, \frac{1-\alpha}{1-\alpha+\alpha^2}\right\}, \left\{\frac{1-2\,\alpha+\alpha^2}{1-\alpha+\alpha^2}, \frac{(-1+\alpha)\,\alpha}{1-\alpha+\alpha^2}\right\}, \left\{1-\alpha+\alpha^2, 1-\alpha\right\},$$

$$\left\{1-\alpha+\alpha^2, \alpha\right\}, \left\{1-\alpha+\alpha^2, \frac{-1+\alpha-\alpha^2}{-1+\alpha}\right\}, \left\{1-\alpha+\alpha^2, \frac{1-\alpha+\alpha^2}{\alpha}\right\}, \left\{\frac{1-\alpha+\alpha^2}{(-1+\alpha)^2}, \frac{1}{1-\alpha}\right\},$$

$$\left\{\frac{1-\alpha+\alpha^2}{(-1+\alpha)^2}, \frac{\alpha}{-1+\alpha}\right\}, \left\{\frac{1-\alpha+\alpha^2}{(-1+\alpha)^2}, \frac{-1+\alpha-\alpha^2}{-1+\alpha}\right\}, \left\{\frac{1-\alpha+\alpha^2}{(-1+\alpha)^2}, \frac{1-\alpha+\alpha^2}{(-1+\alpha)\,\alpha}\right\}, \left\{\frac{1-\alpha+\alpha^2}{\alpha^2}, \frac{1}{\alpha}\right\},$$

$$\left\{\frac{1-\alpha+\alpha^2}{\alpha^2}, \frac{-1+\alpha}{\alpha}\right\}, \left\{\frac{1-\alpha+\alpha^2}{\alpha^2}, \frac{1-\alpha+\alpha^2}{\alpha}\right\}, \left\{\frac{1-\alpha+\alpha^2}{\alpha^2}, \frac{1-\alpha+\alpha^2}{(-1+\alpha)\,\alpha}\right\}, \left\{\frac{1-\alpha+\alpha^2}{\alpha}, \frac{1}{\alpha}\right\},$$

$$\left\{\frac{1-\alpha+\alpha^2}{\alpha}, \alpha\right\}, \left\{\frac{1-\alpha+\alpha^2}{\alpha}, 1-\alpha+\alpha^2\right\}, \left\{\frac{1-\alpha+\alpha^2}{\alpha}, \frac{1-\alpha+\alpha^2}{\alpha^2}\right\}, \left\{\frac{1-\alpha+\alpha^2}{(-1+\alpha)\,\alpha}, \frac{-1+\alpha}{\alpha}\right\},$$

```
Length[%]
```

```
120
```

To show that each of these gives rise to three inequivalent representations over GF(5) we consider the partial-field homomorphism $\phi 3$.

Let $\xi_i : \mathrm{GF}(5) \times \mathrm{GF}(5) \times \mathrm{GF}(5) \to \mathrm{GF}(5)$ be the projection onto the $i$th coordinate. We have to show that $\xi_i(\phi(A))$ is not equivalent to $\xi_j(\phi(A))$ for all matrices $A$ computed above. We only need to look at the pairs $p$, $q$. These are their images under $\phi$:

```
Map[ϕ3, pqPairs, {2}]
```

```
{{{4, 2, 3}, {2, 4, 2}}, {{4, 2, 3}, {3, 4, 2}}, {{4, 2, 3}, {2, 3, 4}}, {{4, 2, 3}, {3, 3, 4}},
 {{4, 2, 3}, {2, 4, 4}}, {{4, 2, 3}, {2, 3, 2}}, {{4, 2, 3}, {3, 4, 4}}, {{4, 2, 3}, {3, 3, 2}},
 {{4, 3, 2}, {3, 2, 4}}, {{4, 3, 2}, {2, 2, 4}}, {{4, 3, 2}, {2, 4, 3}}, {{4, 3, 2}, {3, 4, 3}},
 {{4, 3, 2}, {2, 4, 4}}, {{4, 3, 2}, {3, 4, 4}}, {{4, 3, 2}, {2, 2, 3}}, {{4, 3, 2}, {3, 2, 3}},
 {{4, 3, 3}, {3, 2, 4}}, {{4, 3, 3}, {2, 4, 2}}, {{4, 3, 3}, {3, 4, 2}}, {{4, 3, 3}, {2, 2, 4}},
 {{3, 2, 4}, {4, 3, 2}}, {{3, 2, 4}, {4, 3, 3}}, {{3, 2, 4}, {2, 4, 2}}, {{3, 2, 4}, {2, 4, 3}},
 {{3, 2, 4}, {2, 3, 2}}, {{3, 2, 4}, {4, 4, 3}}, {{3, 2, 4}, {2, 3, 3}}, {{3, 2, 4}, {4, 4, 2}},
 {{2, 4, 2}, {4, 2, 3}}, {{2, 4, 2}, {4, 3, 3}}, {{2, 4, 2}, {3, 2, 4}}, {{2, 4, 2}, {3, 3, 4}},
 {{3, 4, 2}, {4, 2, 3}}, {{3, 4, 2}, {4, 3, 3}}, {{3, 4, 2}, {2, 2, 4}}, {{3, 4, 2}, {2, 3, 4}},
 {{3, 4, 2}, {4, 3, 4}}, {{3, 4, 2}, {2, 2, 3}}, {{3, 4, 2}, {2, 3, 3}}, {{3, 4, 2}, {4, 2, 4}},
 {{2, 2, 4}, {4, 3, 2}}, {{2, 2, 4}, {4, 3, 3}}, {{2, 2, 4}, {3, 4, 2}}, {{2, 2, 4}, {3, 4, 3}},
 {{2, 3, 4}, {4, 2, 3}}, {{2, 3, 4}, {3, 4, 2}}, {{2, 3, 4}, {3, 4, 3}}, {{2, 3, 4}, {4, 2, 2}},
 {{2, 3, 4}, {4, 4, 3}}, {{2, 3, 4}, {3, 2, 2}}, {{2, 3, 4}, {3, 2, 3}}, {{2, 3, 4}, {4, 4, 2}},
 {{3, 3, 4}, {4, 2, 3}}, {{3, 3, 4}, {2, 4, 2}}, {{3, 3, 4}, {2, 4, 3}}, {{3, 3, 4}, {4, 2, 2}},
 {{2, 4, 3}, {4, 3, 2}}, {{2, 4, 3}, {3, 2, 4}}, {{2, 4, 3}, {3, 3, 4}}, {{2, 4, 3}, {4, 2, 2}},
 {{2, 4, 3}, {4, 3, 4}}, {{2, 4, 3}, {3, 2, 2}}, {{2, 4, 3}, {3, 3, 2}}, {{2, 4, 3}, {4, 2, 4}},
 {{3, 4, 3}, {4, 3, 2}}, {{3, 4, 3}, {2, 2, 4}}, {{3, 4, 3}, {2, 3, 4}}, {{3, 4, 3}, {4, 2, 2}},
 {{4, 2, 2}, {2, 3, 4}}, {{4, 2, 2}, {3, 3, 4}}, {{4, 2, 2}, {2, 4, 3}}, {{4, 2, 2}, {3, 4, 3}},
 {{2, 4, 4}, {4, 2, 3}}, {{2, 4, 4}, {4, 3, 2}}, {{2, 4, 4}, {3, 2, 3}}, {{2, 4, 4}, {3, 3, 2}},
 {{2, 3, 2}, {4, 2, 3}}, {{2, 3, 2}, {3, 2, 4}}, {{2, 3, 2}, {3, 4, 4}}, {{2, 3, 2}, {4, 4, 3}},
 {{3, 4, 4}, {4, 2, 3}}, {{3, 4, 4}, {4, 3, 2}}, {{3, 4, 4}, {2, 3, 2}}, {{3, 4, 4}, {2, 2, 3}},
 {{4, 4, 3}, {3, 2, 4}}, {{4, 4, 3}, {2, 3, 4}}, {{4, 4, 3}, {2, 3, 2}}, {{4, 4, 3}, {3, 2, 2}},
 {{4, 3, 4}, {3, 4, 2}}, {{4, 3, 4}, {2, 4, 3}}, {{4, 3, 4}, {3, 2, 2}}, {{4, 3, 4}, {2, 2, 3}},
 {{3, 2, 2}, {2, 3, 4}}, {{3, 2, 2}, {2, 4, 3}}, {{3, 2, 2}, {4, 4, 3}}, {{3, 2, 2}, {4, 3, 4}},
 {{2, 2, 3}, {4, 3, 2}}, {{2, 2, 3}, {3, 4, 2}}, {{2, 2, 3}, {3, 4, 4}}, {{2, 2, 3}, {4, 3, 4}},
 {{3, 2, 3}, {4, 3, 2}}, {{3, 2, 3}, {2, 3, 4}}, {{3, 2, 3}, {2, 4, 4}}, {{3, 2, 3}, {4, 4, 2}},
 {{3, 3, 2}, {4, 2, 3}}, {{3, 3, 2}, {2, 4, 3}}, {{3, 3, 2}, {2, 4, 4}}, {{3, 3, 2}, {4, 2, 4}},
 {{2, 3, 3}, {3, 2, 4}}, {{2, 3, 3}, {3, 4, 2}}, {{2, 3, 3}, {4, 4, 2}}, {{2, 3, 3}, {4, 2, 4}},
 {{4, 4, 2}, {3, 2, 4}}, {{4, 4, 2}, {2, 3, 4}}, {{4, 4, 2}, {3, 2, 3}}, {{4, 4, 2}, {2, 3, 3}},
 {{4, 2, 4}, {3, 4, 2}}, {{4, 2, 4}, {2, 4, 3}}, {{4, 2, 4}, {3, 3, 2}}, {{4, 2, 4}, {2, 3, 3}}}
```

Indeed, each pair of representations is inequivalent. A concise way of checking this is to show that for all $1 \leq i < j \leq 3$, $\left| \xi_i(p) - \xi_j(p) \right| + \left| \xi_i(q) - \xi_j(q) \right| > 0$.

```
Map[{Abs[#[[1, 1]] - #[[1, 2]]] + Abs[#[[2, 1]] - #[[2, 2]]],
   Abs[#[[1, 1]] - #[[1, 3]]] + Abs[#[[2, 1]] - #[[2, 3]]],
   Abs[#[[1, 2]] - #[[1, 3]]] + Abs[#[[2, 2]] - #[[2, 3]]]} &, %, {1}]
```

```
{{4, 1, 3}, {3, 2, 3}, {3, 3, 2}, {2, 2, 2}, {4, 3, 1}, {3, 1, 2}, {3, 2, 1}, {2, 2, 2},
 {2, 3, 3}, {1, 4, 3}, {3, 3, 2}, {2, 2, 2}, {3, 4, 1}, {2, 3, 1}, {1, 3, 2}, {2, 2, 2},
 {2, 2, 2}, {3, 1, 2}, {2, 2, 2}, {1, 3, 2}, {2, 3, 3}, {2, 2, 2}, {3, 1, 4}, {3, 2, 3},
 {2, 1, 3}, {1, 2, 3}, {2, 2, 2}, {1, 3, 4}, {4, 1, 3}, {3, 1, 2}, {3, 1, 4}, {2, 1, 3},
 {3, 2, 3}, {2, 2, 2}, {1, 3, 4}, {2, 3, 3}, {2, 1, 3}, {1, 2, 3}, {2, 2, 2}, {3, 1, 4},
 {1, 4, 3}, {1, 3, 2}, {1, 3, 4}, {1, 2, 3}, {3, 3, 2}, {2, 3, 3}, {2, 2, 2}, {3, 4, 1},
 {1, 3, 2}, {2, 3, 1}, {2, 2, 2}, {1, 4, 3}, {2, 2, 2}, {2, 1, 3}, {2, 2, 2}, {2, 3, 1},
 {3, 3, 2}, {3, 2, 3}, {2, 2, 2}, {4, 3, 1}, {3, 1, 2}, {3, 2, 1}, {2, 2, 2}, {4, 1, 3},
 {2, 2, 2}, {1, 2, 3}, {2, 2, 2}, {3, 2, 1}, {3, 4, 1}, {2, 3, 1}, {4, 3, 1}, {3, 2, 1},
 {4, 3, 1}, {3, 4, 1}, {3, 2, 1}, {2, 3, 1}, {3, 1, 2}, {2, 1, 3}, {2, 1, 1}, {1, 1, 2},
 {3, 2, 1}, {2, 3, 1}, {2, 1, 1}, {1, 2, 1}, {1, 2, 3}, {1, 3, 2}, {1, 1, 2}, {1, 2, 1},
 {2, 1, 3}, {3, 1, 2}, {2, 1, 1}, {1, 1, 2}, {2, 3, 1}, {3, 2, 1}, {1, 2, 1}, {2, 1, 1},
 {1, 3, 2}, {1, 2, 3}, {1, 2, 1}, {1, 1, 2}, {2, 2, 2}, {2, 2, 2}, {3, 2, 1}, {1, 2, 3},
 {2, 2, 2}, {2, 2, 2}, {2, 3, 1}, {2, 1, 3}, {2, 2, 2}, {2, 2, 2}, {1, 3, 2}, {3, 1, 2},
 {1, 3, 4}, {1, 4, 3}, {1, 2, 3}, {1, 3, 2}, {3, 1, 4}, {4, 1, 3}, {2, 1, 3}, {3, 1, 2}}
```

```
Position[%, 0]
```

```
{}
```

This completes the proof of Lemma 21.

## Lifting

We prove the following fact:

**Lemma 22.** Let $M$ be a 3-connected matroid with at least three inequivalent representations over GF(5). Then $M$ is representable over $\mathbb{H}_3$.

Consider a $\mathrm{GF}(5)\times\mathrm{GF}(5)\times\mathrm{GF}(5)$-matrix $A$ representing $M$ such that the three projections $\xi_i(A)$ are pairwise inequivalent. We want to apply the Lift Theorem, so we start by constructing a lifting function for $\mathrm{Cr}(A)$. First we need to find out what $\mathrm{Cr}(A)$ is.

**Claim 22.1.** No element of $\mathrm{Cr}(A)$ is of the form $(x, x, x)$ for $x \notin \{0, 1\}$.

*Proof.* Assume that $(x, x, x) \in \mathrm{Cr}(A)$ for some $x \in \{2, 3, 4\}$. Let $\mathbb{P}'$ be the sub-partial field of $\mathrm{GF}(5)\times\mathrm{GF}(5)\times\mathrm{GF}(5)$ such that the three coordinates of all elements are equal. Since $A$ is not a $\mathbb{P}'$-matrix, the Confinement Theorem, in particular Corollary 15, implies that either $A$ or $A^T$ contains a minor

$$D = \begin{pmatrix} (1,1,1) & (1,1,1) & (1,1,1) \\ (1,1,1) & (x,x,x) & (p,q,r) \end{pmatrix},$$

with not all three of $p$, $q$, $r$ equal. There are only two elements $u$ of GF(5) such that $\begin{pmatrix} 1 & 1 & 1 \\ 1 & x & u \end{pmatrix}$ represents $U_{2,5}$. Hence two of the three projections of $D$ must be equivalent. But then the Stabilizer Theorem, in particular Lemma 8, implies that two of the three projections of $A$ are equivalent, a contradiction. □

We define $F' := \mathcal{F}(\mathrm{GF}(5)\times\mathrm{GF}(5)\times\mathrm{GF}(5)) \setminus \{(2,2,2), (3,3,3), (4,4,4)\}$.

```
Fp = Complement[
  {{0, 0, 0}, {1, 1, 1}} ⋃ Flatten[Table[{x, y, z}, {x, 2, 4}, {y, 2, 4}, {z, 2, 4}], 2],
  {{2, 2, 2}, {3, 3, 3}, {4, 4, 4}}]
```

```
{{0, 0, 0}, {1, 1, 1}, {2, 2, 3}, {2, 2, 4}, {2, 3, 2}, {2, 3, 3},
 {2, 3, 4}, {2, 4, 2}, {2, 4, 3}, {2, 4, 4}, {3, 2, 2}, {3, 2, 3},
 {3, 2, 4}, {3, 3, 2}, {3, 3, 4}, {3, 4, 2}, {3, 4, 3}, {3, 4, 4}, {4, 2, 2},
 {4, 2, 3}, {4, 2, 4}, {4, 3, 2}, {4, 3, 3}, {4, 3, 4}, {4, 4, 2}, {4, 4, 3}}
```

```
Length[%]
```

```
26
```

Consider the homomorphism $\phi 3$. The restriction of $\phi 3$ to $\mathcal{F}(\mathbb{H}_3)$ is a bijection between $\mathcal{F}(\mathbb{H}_3)$ and $F'$:

```
Sort[ϕ3 /@ H3funs]
```

```
{{0, 0, 0}, {1, 1, 1}, {2, 2, 3}, {2, 2, 4}, {2, 3, 2}, {2, 3, 3},
 {2, 3, 4}, {2, 4, 2}, {2, 4, 3}, {2, 4, 4}, {3, 2, 2}, {3, 2, 3},
 {3, 2, 4}, {3, 3, 2}, {3, 3, 4}, {3, 4, 2}, {3, 4, 3}, {3, 4, 4}, {4, 2, 2},
 {4, 2, 3}, {4, 2, 4}, {4, 3, 2}, {4, 3, 3}, {4, 3, 4}, {4, 4, 2}, {4, 4, 3}}
```

It follows immediately that the function $F' \to \mathbb{H}_3$ that is the inverse of $\phi$ is a lifting function.

```
up3[x_, y_, z_] := H3funs[[Position[ϕ3 /@ H3funs, {x, y, z}][[1, 1]]]];
up3[{x_, y_, z_}] := up3[x, y, z];
```

Some examples:

```
up3[1, 1, 1]
```

```
1
```

```
up3[2, 4, 4]
```

$$\frac{-1+\alpha-\alpha^2}{-1+\alpha}$$

Now suppose that $A$ has no global lift. Then the Lift Theorem (Theorem 12) tells us that $A$ or $A^T$ has a minor of the form $\begin{pmatrix} 1 & 1 & 0 & 1 \\ 1 & 0 & 1 & 1 \\ 0 & 1 & 1 & 1 \end{pmatrix}$ or $\begin{pmatrix} 1 & 1 & 1 \\ 1 & p & q \end{pmatrix}$, for some $p, q \in \text{Cr}(A)$. without a local lift. The first of these has a cross ratio (2, 2, 2), so it does not occur as a minor of $A$, by Claim 22.1. Therefore we have to ensure that all $\text{GF}(5) \times \text{GF}(5) \times \text{GF}(5)$-representations of $U_{2,5}$ that can occur in $A$ have a local lift. Let $D$ be such a minor. Suppose that $\xi_i(D)$ and $\xi_j(D)$ are equivalent. By the Stabilizer Theorem, and in particular Lemma 8, this implies that $\xi_i(A)$ and $\xi_j(A)$ are equivalent, a contradiction to our choice of $A$. Hence $\xi_i(D)$ and $\xi_j(D)$ are inequivalent. This gives us $6 \times 5 \times 4 = 120$ representations of $U_{2,5}$. Again we only list the pairs $(p, q)$ in

$$D = \begin{pmatrix} 1 & 1 & 1 \\ 1 & p & q \end{pmatrix}.$$

```
U25reps = { {2, 3}, {2, 4}, {3, 2}, {3, 4}, {4, 2}, {4, 3}}
```

```
{{2, 3}, {2, 4}, {3, 2}, {3, 4}, {4, 2}, {4, 3}}
```

```
U25reptriples = Transpose /@ Permutations[U25reps, {3}]
```

```
{{{2, 2, 3}, {3, 4, 2}}, {{2, 2, 3}, {3, 4, 4}}, {{2, 2, 4}, {3, 4, 2}}, {{2, 2, 4}, {3, 4, 3}},
 {{2, 3, 2}, {3, 2, 4}}, {{2, 3, 3}, {3, 2, 4}}, {{2, 3, 4}, {3, 2, 2}}, {{2, 3, 4}, {3, 2, 3}},
 {{2, 3, 2}, {3, 4, 4}}, {{2, 3, 3}, {3, 4, 2}}, {{2, 3, 4}, {3, 4, 2}}, {{2, 3, 4}, {3, 4, 3}},
 {{2, 4, 2}, {3, 2, 4}}, {{2, 4, 3}, {3, 2, 2}}, {{2, 4, 3}, {3, 2, 4}}, {{2, 4, 4}, {3, 2, 3}},
 {{2, 4, 2}, {3, 3, 4}}, {{2, 4, 3}, {3, 3, 2}}, {{2, 4, 3}, {3, 3, 4}}, {{2, 4, 4}, {3, 3, 2}},
 {{2, 2, 3}, {4, 3, 2}}, {{2, 2, 3}, {4, 3, 4}}, {{2, 2, 4}, {4, 3, 2}}, {{2, 2, 4}, {4, 3, 3}},
 {{2, 3, 2}, {4, 2, 3}}, {{2, 3, 3}, {4, 2, 4}}, {{2, 3, 4}, {4, 2, 2}}, {{2, 3, 4}, {4, 2, 3}},
 {{2, 3, 2}, {4, 4, 3}}, {{2, 3, 3}, {4, 4, 2}}, {{2, 3, 4}, {4, 4, 2}}, {{2, 3, 4}, {4, 4, 3}},
 {{2, 4, 2}, {4, 2, 3}}, {{2, 4, 3}, {4, 2, 2}}, {{2, 4, 3}, {4, 2, 4}}, {{2, 4, 4}, {4, 2, 3}},
 {{2, 4, 2}, {4, 3, 3}}, {{2, 4, 3}, {4, 3, 2}}, {{2, 4, 3}, {4, 3, 4}}, {{2, 4, 4}, {4, 3, 2}},
 {{3, 2, 2}, {2, 3, 4}}, {{3, 2, 3}, {2, 3, 4}}, {{3, 2, 4}, {2, 3, 2}}, {{3, 2, 4}, {2, 3, 3}},
 {{3, 2, 2}, {2, 4, 3}}, {{3, 2, 3}, {2, 4, 4}}, {{3, 2, 4}, {2, 4, 2}}, {{3, 2, 4}, {2, 4, 3}},
 {{3, 3, 2}, {2, 4, 3}}, {{3, 3, 2}, {2, 4, 4}}, {{3, 3, 4}, {2, 4, 2}}, {{3, 3, 4}, {2, 4, 3}},
 {{3, 4, 2}, {2, 2, 3}}, {{3, 4, 2}, {2, 2, 4}}, {{3, 4, 3}, {2, 2, 4}}, {{3, 4, 4}, {2, 2, 3}},
 {{3, 4, 2}, {2, 3, 3}}, {{3, 4, 2}, {2, 3, 4}}, {{3, 4, 3}, {2, 3, 4}}, {{3, 4, 4}, {2, 3, 2}},
 {{3, 2, 2}, {4, 3, 4}}, {{3, 2, 3}, {4, 3, 2}}, {{3, 2, 4}, {4, 3, 2}}, {{3, 2, 4}, {4, 3, 3}},
 {{3, 2, 2}, {4, 4, 3}}, {{3, 2, 3}, {4, 4, 2}}, {{3, 2, 4}, {4, 4, 2}}, {{3, 2, 4}, {4, 4, 3}},
 {{3, 3, 2}, {4, 2, 3}}, {{3, 3, 2}, {4, 2, 4}}, {{3, 3, 4}, {4, 2, 2}}, {{3, 3, 4}, {4, 2, 3}},
 {{3, 4, 2}, {4, 2, 3}}, {{3, 4, 2}, {4, 2, 4}}, {{3, 4, 3}, {4, 2, 2}}, {{3, 4, 4}, {4, 2, 3}},
 {{3, 4, 2}, {4, 3, 3}}, {{3, 4, 2}, {4, 3, 4}}, {{3, 4, 3}, {4, 3, 2}}, {{3, 4, 4}, {4, 3, 2}},
 {{4, 2, 2}, {2, 3, 4}}, {{4, 2, 3}, {2, 3, 2}}, {{4, 2, 3}, {2, 3, 4}}, {{4, 2, 4}, {2, 3, 3}},
 {{4, 2, 2}, {2, 4, 3}}, {{4, 2, 3}, {2, 4, 2}}, {{4, 2, 3}, {2, 4, 4}}, {{4, 2, 4}, {2, 4, 3}},
 {{4, 3, 2}, {2, 2, 3}}, {{4, 3, 2}, {2, 2, 4}}, {{4, 3, 3}, {2, 2, 4}}, {{4, 3, 4}, {2, 2, 3}},
 {{4, 3, 2}, {2, 4, 3}}, {{4, 3, 2}, {2, 4, 4}}, {{4, 3, 3}, {2, 4, 2}}, {{4, 3, 4}, {2, 4, 3}},
 {{4, 4, 2}, {2, 3, 3}}, {{4, 4, 2}, {2, 3, 4}}, {{4, 4, 3}, {2, 3, 2}}, {{4, 4, 3}, {2, 3, 4}},
 {{4, 2, 2}, {3, 3, 4}}, {{4, 2, 3}, {3, 3, 2}}, {{4, 2, 3}, {3, 3, 4}}, {{4, 2, 4}, {3, 3, 2}},
 {{4, 2, 2}, {3, 4, 3}}, {{4, 2, 3}, {3, 4, 2}}, {{4, 2, 3}, {3, 4, 4}}, {{4, 2, 4}, {3, 4, 2}},
 {{4, 3, 2}, {3, 2, 3}}, {{4, 3, 2}, {3, 2, 4}}, {{4, 3, 3}, {3, 2, 4}}, {{4, 3, 4}, {3, 2, 2}},
 {{4, 3, 2}, {3, 4, 3}}, {{4, 3, 2}, {3, 4, 4}}, {{4, 3, 3}, {3, 4, 2}}, {{4, 3, 4}, {3, 4, 2}},
 {{4, 4, 2}, {3, 2, 3}}, {{4, 4, 2}, {3, 2, 4}}, {{4, 4, 3}, {3, 2, 2}}, {{4, 4, 3}, {3, 2, 4}}}
```

```
Length[%]
```

```
120
```

The lifts look as follows:

```
Map[up3, U25reptriples, {2}]
```

$$\Big\{\Big\{\frac{1-2\alpha+\alpha^2}{1-\alpha+\alpha^2}, \frac{-1+\alpha}{\alpha}\Big\}, \Big\{\frac{1-2\alpha+\alpha^2}{1-\alpha+\alpha^2}, \frac{1-\alpha}{1-\alpha+\alpha^2}\Big\}, \Big\{-\frac{(-1+\alpha)^2}{\alpha}, \frac{-1+\alpha}{\alpha}\Big\}, \Big\{-\frac{(-1+\alpha)^2}{\alpha}, -(-1+\alpha)\,\alpha\Big\},$$

$$\Big\{\frac{1}{1-\alpha+\alpha^2}, \frac{1}{\alpha}\Big\}, \Big\{\frac{1-\alpha+\alpha^2}{\alpha^2}, \frac{1}{\alpha}\Big\}, \Big\{\alpha, \frac{\alpha^2}{1-\alpha+\alpha^2}\Big\}, \{\alpha, 1-\alpha+\alpha^2\}, \Big\{\frac{1}{1-\alpha+\alpha^2}, \frac{1-\alpha}{1-\alpha+\alpha^2}\Big\},$$

$$\Big\{\frac{1-\alpha+\alpha^2}{\alpha^2}, \frac{-1+\alpha}{\alpha}\Big\}, \Big\{\alpha, \frac{-1+\alpha}{\alpha}\Big\}, \{\alpha, -(-1+\alpha)\,\alpha\}, \Big\{\frac{1}{(1-\alpha)\,\alpha}, \frac{1}{\alpha}\Big\}, \Big\{\frac{\alpha}{-1+\alpha}, \frac{\alpha^2}{1-\alpha+\alpha^2}\Big\},$$

$$\Big\{\frac{\alpha}{-1+\alpha}, \frac{1}{\alpha}\Big\}, \Big\{\frac{-1+\alpha-\alpha^2}{-1+\alpha}, 1-\alpha+\alpha^2\Big\}, \Big\{\frac{1}{(1-\alpha)\,\alpha}, -\frac{\alpha}{(-1+\alpha)^2}\Big\}, \Big\{\frac{\alpha}{-1+\alpha}, \frac{1-\alpha+\alpha^2}{(-1+\alpha)^2}\Big\},$$

$$\left\{\frac{\alpha}{-1+\alpha}, -\frac{\alpha}{(-1+\alpha)^2}\right\}, \left\{\frac{-1+\alpha-\alpha^2}{-1+\alpha}, \frac{1-\alpha+\alpha^2}{(-1+\alpha)^2}\right\}, \left\{\frac{1-2\alpha+\alpha^2}{1-\alpha+\alpha^2}, 1-\alpha\right\}, \left\{\frac{1-2\alpha+\alpha^2}{1-\alpha+\alpha^2}, \frac{(-1+\alpha)\alpha}{1-\alpha+\alpha^2}\right\},$$

$$\left\{-\frac{(-1+\alpha)^2}{\alpha}, 1-\alpha\right\}, \left\{-\frac{(-1+\alpha)^2}{\alpha}, \frac{-1+\alpha}{\alpha^2}\right\}, \left\{\frac{1}{1-\alpha+\alpha^2}, \frac{1}{1-\alpha}\right\}, \left\{\frac{1-\alpha+\alpha^2}{\alpha^2}, \frac{1-\alpha+\alpha^2}{(-1+\alpha)\alpha}\right\},$$

$$\left\{\alpha, \frac{\alpha^2}{-1+\alpha}\right\}, \left\{\alpha, \frac{1}{1-\alpha}\right\}, \left\{\frac{1}{1-\alpha+\alpha^2}, \frac{\alpha}{1-\alpha+\alpha^2}\right\}, \left\{\frac{1-\alpha+\alpha^2}{\alpha^2}, \frac{1-\alpha+\alpha^2}{\alpha}\right\}, \left\{\alpha, \frac{1-\alpha+\alpha^2}{\alpha}\right\},$$

$$\left\{\alpha, \frac{\alpha}{1-\alpha+\alpha^2}\right\}, \left\{\frac{1}{(1-\alpha)\alpha}, \frac{1}{1-\alpha}\right\}, \left\{\frac{\alpha}{-1+\alpha}, \frac{\alpha^2}{-1+\alpha}\right\}, \left\{\frac{\alpha}{-1+\alpha}, \frac{1-\alpha+\alpha^2}{(-1+\alpha)\alpha}\right\}, \left\{\frac{-1+\alpha-\alpha^2}{-1+\alpha}, \frac{1}{1-\alpha}\right\},$$

$$\left\{\frac{1}{(1-\alpha)\alpha}, \frac{-1+\alpha}{\alpha^2}\right\}, \left\{\frac{\alpha}{-1+\alpha}, 1-\alpha\right\}, \left\{\frac{\alpha}{-1+\alpha}, \frac{(-1+\alpha)\alpha}{1-\alpha+\alpha^2}\right\}, \left\{\frac{-1+\alpha-\alpha^2}{-1+\alpha}, 1-\alpha\right\},$$

$$\left\{\frac{\alpha^2}{1-\alpha+\alpha^2}, \alpha\right\}, \left\{1-\alpha+\alpha^2, \alpha\right\}, \left\{\frac{1}{\alpha}, \frac{1}{1-\alpha+\alpha^2}\right\}, \left\{\frac{1}{\alpha}, \frac{1-\alpha+\alpha^2}{\alpha^2}\right\}, \left\{\frac{\alpha^2}{1-\alpha+\alpha^2}, \frac{\alpha}{-1+\alpha}\right\},$$

$$\left\{1-\alpha+\alpha^2, \frac{-1+\alpha-\alpha^2}{-1+\alpha}\right\}, \left\{\frac{1}{\alpha}, \frac{1}{(1-\alpha)\alpha}\right\}, \left\{\frac{1}{\alpha}, \frac{\alpha}{-1+\alpha}\right\}, \left\{\frac{1-\alpha+\alpha^2}{(-1+\alpha)^2}, \frac{\alpha}{-1+\alpha}\right\},$$

$$\left\{\frac{1-\alpha+\alpha^2}{(-1+\alpha)^2}, \frac{-1+\alpha-\alpha^2}{-1+\alpha}\right\}, \left\{-\frac{\alpha}{(-1+\alpha)^2}, \frac{1}{(1-\alpha)\alpha}\right\}, \left\{-\frac{\alpha}{(-1+\alpha)^2}, \frac{\alpha}{-1+\alpha}\right\}, \left\{\frac{-1+\alpha}{\alpha}, \frac{1-2\alpha+\alpha^2}{1-\alpha+\alpha^2}\right\},$$

$$\left\{\frac{-1+\alpha}{\alpha}, -\frac{(-1+\alpha)^2}{\alpha}\right\}, \left\{-(-1+\alpha)\alpha, -\frac{(-1+\alpha)^2}{\alpha}\right\}, \left\{\frac{1-\alpha}{1-\alpha+\alpha^2}, \frac{1-2\alpha+\alpha^2}{1-\alpha+\alpha^2}\right\}, \left\{\frac{-1+\alpha}{\alpha}, \frac{1-\alpha+\alpha^2}{\alpha^2}\right\},$$

$$\left\{\frac{-1+\alpha}{\alpha}, \alpha\right\}, \left\{-(-1+\alpha)\alpha, \alpha\right\}, \left\{\frac{1-\alpha}{1-\alpha+\alpha^2}, \frac{1}{1-\alpha+\alpha^2}\right\}, \left\{\frac{\alpha^2}{1-\alpha+\alpha^2}, \frac{(-1+\alpha)\alpha}{1-\alpha+\alpha^2}\right\},$$

$$\left\{1-\alpha+\alpha^2, 1-\alpha\right\}, \left\{\frac{1}{\alpha}, 1-\alpha\right\}, \left\{\frac{1}{\alpha}, \frac{-1+\alpha}{\alpha^2}\right\}, \left\{\frac{\alpha^2}{1-\alpha+\alpha^2}, \frac{\alpha}{1-\alpha+\alpha^2}\right\}, \left\{1-\alpha+\alpha^2, \frac{1-\alpha+\alpha^2}{\alpha}\right\},$$

$$\left\{\frac{1}{\alpha}, \frac{1-\alpha+\alpha^2}{\alpha}\right\}, \left\{\frac{1}{\alpha}, \frac{\alpha}{1-\alpha+\alpha^2}\right\}, \left\{\frac{1-\alpha+\alpha^2}{(-1+\alpha)^2}, \frac{1}{1-\alpha}\right\}, \left\{\frac{1-\alpha+\alpha^2}{(-1+\alpha)^2}, \frac{1-\alpha+\alpha^2}{(-1+\alpha)\alpha}\right\},$$

$$\left\{-\frac{\alpha}{(-1+\alpha)^2}, \frac{\alpha^2}{-1+\alpha}\right\}, \left\{-\frac{\alpha}{(-1+\alpha)^2}, \frac{1}{1-\alpha}\right\}, \left\{\frac{-1+\alpha}{\alpha}, \frac{1}{1-\alpha}\right\}, \left\{\frac{-1+\alpha}{\alpha}, \frac{1-\alpha+\alpha^2}{(-1+\alpha)\alpha}\right\},$$

$$\left\{-(-1+\alpha)\alpha, \frac{\alpha^2}{-1+\alpha}\right\}, \left\{\frac{1-\alpha}{1-\alpha+\alpha^2}, \frac{1}{1-\alpha}\right\}, \left\{\frac{-1+\alpha}{\alpha}, \frac{-1+\alpha}{\alpha^2}\right\}, \left\{\frac{-1+\alpha}{\alpha}, \frac{(-1+\alpha)\alpha}{1-\alpha+\alpha^2}\right\},$$

$$\left\{-(-1+\alpha)\alpha, 1-\alpha\right\}, \left\{\frac{1-\alpha}{1-\alpha+\alpha^2}, 1-\alpha\right\}, \left\{\frac{\alpha^2}{-1+\alpha}, \alpha\right\}, \left\{\frac{1}{1-\alpha}, \frac{1}{1-\alpha+\alpha^2}\right\}, \left\{\frac{1}{1-\alpha}, \alpha\right\},$$

$$\left\{\frac{1-\alpha+\alpha^2}{(-1+\alpha)\alpha}, \frac{1-\alpha+\alpha^2}{\alpha^2}\right\}, \left\{\frac{\alpha^2}{-1+\alpha}, \frac{\alpha}{-1+\alpha}\right\}, \left\{\frac{1}{1-\alpha}, \frac{1}{(1-\alpha)\alpha}\right\}, \left\{\frac{1}{1-\alpha}, \frac{-1+\alpha-\alpha^2}{-1+\alpha}\right\},$$

$$\left\{\frac{1-\alpha+\alpha^2}{(-1+\alpha)\alpha}, \frac{\alpha}{-1+\alpha}\right\}, \left\{1-\alpha, \frac{1-2\alpha+\alpha^2}{1-\alpha+\alpha^2}\right\}, \left\{1-\alpha, -\frac{(-1+\alpha)^2}{\alpha}\right\}, \left\{\frac{-1+\alpha}{\alpha^2}, -\frac{(-1+\alpha)^2}{\alpha}\right\},$$

$$\left\{\frac{(-1+\alpha)\alpha}{1-\alpha+\alpha^2}, \frac{1-2\alpha+\alpha^2}{1-\alpha+\alpha^2}\right\}, \left\{1-\alpha, \frac{\alpha}{-1+\alpha}\right\}, \left\{1-\alpha, \frac{-1+\alpha-\alpha^2}{-1+\alpha}\right\}, \left\{\frac{-1+\alpha}{\alpha^2}, \frac{1}{(1-\alpha)\alpha}\right\},$$

$$\left\{\frac{(-1+\alpha)\alpha}{1-\alpha+\alpha^2}, \frac{\alpha}{-1+\alpha}\right\}, \left\{\frac{1-\alpha+\alpha^2}{\alpha}, \frac{1-\alpha+\alpha^2}{\alpha^2}\right\}, \left\{\frac{1-\alpha+\alpha^2}{\alpha}, \alpha\right\}, \left\{\frac{\alpha}{1-\alpha+\alpha^2}, \frac{1}{1-\alpha+\alpha^2}\right\},$$

$$\left\{\frac{\alpha}{1-\alpha+\alpha^2}, \alpha\right\}, \left\{\frac{\alpha^2}{-1+\alpha}, -\frac{\alpha}{(-1+\alpha)^2}\right\}, \left\{\frac{1}{1-\alpha}, \frac{1-\alpha+\alpha^2}{(-1+\alpha)^2}\right\}, \left\{\frac{1}{1-\alpha}, -\frac{\alpha}{(-1+\alpha)^2}\right\},$$

$$\left\{\frac{1-\alpha+\alpha^2}{(-1+\alpha)\alpha}, \frac{1-\alpha+\alpha^2}{(-1+\alpha)^2}\right\}, \left\{\frac{\alpha^2}{-1+\alpha}, -(-1+\alpha)\alpha\right\}, \left\{\frac{1}{1-\alpha}, \frac{-1+\alpha}{\alpha}\right\}, \left\{\frac{1}{1-\alpha}, \frac{1-\alpha}{1-\alpha+\alpha^2}\right\},$$

$$\left\{\frac{1-\alpha+\alpha^2}{(-1+\alpha)\,\alpha}, \frac{-1+\alpha}{\alpha}\right\}, \{1-\alpha, 1-\alpha+\alpha^2\}, \left\{1-\alpha, \frac{1}{\alpha}\right\}, \left\{\frac{-1+\alpha}{\alpha^2}, \frac{1}{\alpha}\right\}, \left\{\frac{(-1+\alpha)\,\alpha}{1-\alpha+\alpha^2}, \frac{\alpha^2}{1-\alpha+\alpha^2}\right\},$$
$$\{1-\alpha, -(-1+\alpha)\,\alpha\}, \left\{1-\alpha, \frac{1-\alpha}{1-\alpha+\alpha^2}\right\}, \left\{\frac{-1+\alpha}{\alpha^2}, \frac{-1+\alpha}{\alpha}\right\}, \left\{\frac{(-1+\alpha)\,\alpha}{1-\alpha+\alpha^2}, \frac{-1+\alpha}{\alpha}\right\},$$
$$\left\{\frac{1-\alpha+\alpha^2}{\alpha}, 1-\alpha+\alpha^2\right\}, \left\{\frac{1-\alpha+\alpha^2}{\alpha}, \frac{1}{\alpha}\right\}, \left\{\frac{\alpha}{1-\alpha+\alpha^2}, \frac{\alpha^2}{1-\alpha+\alpha^2}\right\}, \left\{\frac{\alpha}{1-\alpha+\alpha^2}, \frac{1}{\alpha}\right\}\}$$

And indeed, for each of these $p^\uparrow / q^\uparrow$ is a fundamental element:

```
Select[%, Not[occursInList[nonzerooneH3funs, #[[1]] / #[[2]]]] &]
```

```
{}
```

It follows that $A$ does have a global lift. This completes the proof of Lemma 22.

# Hydra-4

## Definition and homomorphisms

The Hydra-4 partial field is
$\mathbb{H}_4 = (\mathbb{Q}(\alpha, \beta), \langle -1, \alpha, \beta, 1-\alpha, 1-\beta, \alpha\,\beta-1, \alpha+\beta-2\,\alpha\,\beta\rangle)$,
where $\alpha$, $\beta$ are indeterminates.

The generators of $\mathbb{H}_4$ are

```
H4gens = {-1, α, β, 1 - α, 1 - β, α β - 1, α + β - 2 α β};
```

The fundamental elements of $\mathbb{H}_4$ are

```
H4funs = Together[associates[
    {1, α, β, α β, (α - 1)/(α β - 1), (β - 1)/(α β - 1), -(α (β - 1))/(β (α - 1)), ((α - 1) (β - 1))/(1 - α β), (α (β - 1)^2)/(β (α β - 1)), (β (α - 1)^2)/(α (α β - 1))}]]
```

$$\Big\{0,\ 1,\ \frac{1}{1-\alpha},\ 1-\alpha,\ \frac{1}{\alpha},\ \frac{-1+\alpha}{\alpha},\ \alpha,\ \frac{\alpha}{-1+\alpha},\ \frac{-1+\alpha\beta}{\alpha\beta},\ \frac{1}{1-\beta},\ 1-\beta,\ \frac{1-\alpha}{\alpha(-1+\beta)},\ \frac{1}{\beta},\ \frac{1}{\alpha\beta},\ \frac{1-\beta}{(-1+\alpha)\beta},$$
$$\frac{-1+\beta}{\beta},\ -\frac{\alpha(-1+\beta)}{(-1+\alpha)\beta},\ \beta,\ \alpha\beta,\ \frac{\beta}{-1+\beta},\ -\frac{(-1+\alpha)\beta}{\alpha(-1+\beta)},\ -\frac{\alpha(-1+\beta)^2}{-\alpha-\beta+2\alpha\beta},\ -\frac{(-1+\alpha)^2\beta}{-\alpha-\beta+2\alpha\beta},$$
$$\frac{\alpha+\beta-2\alpha\beta}{\alpha(-1+\beta)^2},\ \frac{-\alpha-\beta+2\alpha\beta}{(-1+\alpha)(-1+\beta)},\ \frac{\alpha+\beta-2\alpha\beta}{(-1+\alpha)^2\beta},\ \frac{1}{1-\alpha\beta},\ -\frac{(-1+\alpha)(-1+\beta)}{-1+\alpha\beta},\ \frac{-\alpha-\beta+2\alpha\beta}{-1+\alpha\beta},$$
$$1-\alpha\beta,\ \frac{1-\alpha\beta}{(-1+\alpha)(-1+\beta)},\ \frac{-1+\alpha\beta}{-\alpha-\beta+2\alpha\beta},\ \frac{-\alpha-\beta+2\alpha\beta}{\alpha(-1+\beta)},\ \frac{-1+\alpha\beta}{\alpha(-1+\beta)},\ \frac{\alpha-\alpha\beta}{-1+\alpha},\ \frac{-\alpha+\alpha\beta}{-\alpha-\beta+2\alpha\beta},$$
$$\frac{-\alpha-\beta+2\alpha\beta}{(-1+\alpha)\beta},\ \frac{-1+\alpha\beta}{(-1+\alpha)\beta},\ \frac{\beta-\alpha\beta}{-1+\beta},\ \frac{-\beta+\alpha\beta}{-\alpha-\beta+2\alpha\beta},\ \frac{-1+\alpha}{-1+\alpha\beta},\ \frac{-1+\beta}{-1+\alpha\beta},\ \frac{\alpha(-1+\beta)}{-1+\alpha\beta},$$
$$\frac{\alpha(-1+\beta)^2}{\beta(-1+\alpha\beta)},\ \frac{(-1+\alpha)\beta}{-1+\alpha\beta},\ \frac{(-1+\alpha)^2\beta}{\alpha(-1+\alpha\beta)},\ \frac{\alpha\beta}{-1+\alpha\beta},\ \frac{-1+\alpha\beta}{-1+\alpha},\ \frac{-1+\alpha\beta}{-1+\beta},\ \frac{\alpha(-1+\alpha\beta)}{(-1+\alpha)^2\beta},$$
$$\frac{\beta(-1+\alpha\beta)}{\alpha(-1+\beta)^2},\ \frac{-\alpha-\beta+2\alpha\beta}{\alpha(-1+\alpha\beta)},\ \frac{-\alpha+\alpha^2\beta}{-\alpha-\beta+2\alpha\beta},\ \frac{-\alpha-\beta+2\alpha\beta}{\beta(-1+\alpha\beta)},\ \frac{-\beta+\alpha\beta^2}{-\alpha-\beta+2\alpha\beta},\ \frac{(-1+\alpha)(-1+\beta)}{-\alpha-\beta+2\alpha\beta}\Big\}$$

```
Length[%]
```

```
56
```

The fundamental elements other than 0, 1 play a special role:

```
nonzerooneH4funs = Complement[H4funs, {0, 1}];
```

There are four homomorphisms $\mathbb{H}_4 \to \mathrm{GF}(5)$. We collect them in the partial-field homomorphism $\phi 4 : \mathbb{H}_4 \to \mathrm{GF}(5) \times \mathrm{GF}(5) \times \mathrm{GF}(5) \times \mathrm{GF}(5)$ defined by $\phi 4(\alpha) = (2, 3, 3, 4)$ and $\phi 4(\beta) = (2, 3, 4, 3)$.

```
ϕ4[0] = {0, 0, 0, 0};
ϕ4[1] = {1, 1, 1, 1};
ϕ4[-1] = {4, 4, 4, 4};
ϕ4[x_] := toGF5[x /. {α → {2, 3, 3, 4}, β → {2, 3, 4, 3}}];
```

This is indeed a partial-field homomorphism:

```
ϕ4 /@ H4gens
```

```
{{4, 4, 4, 4}, {2, 3, 3, 4}, {2, 3, 4, 3}, {4, 3, 3, 2}, {4, 3, 2, 3}, {3, 3, 1, 1}, {1, 3, 3, 3}}
```

## Fundamental elements

Our aim in this section is to verify the following statement:

**Lemma 23.** The fundamental elements of $\mathbb{H}_4$ are

$$\text{Asc}\left\{1, \alpha, \beta, \alpha\beta, \frac{\alpha-1}{\alpha\beta-1}, \frac{\beta-1}{\alpha\beta-1}, -\frac{\alpha(\beta-1)}{\beta(\alpha-1)}, \frac{(\alpha-1)(\beta-1)}{1-\alpha\beta}, \frac{\alpha(\beta-1)^2}{\beta(\alpha\beta-1)}, \frac{\beta(\alpha-1)^2}{\alpha(\alpha\beta-1)}\right\}.$$

Elements of this partial field are of the form
$\pm\alpha^u \beta^v (1-\alpha)^w (1-\beta)^x (\alpha\beta-1)^y (\alpha+\beta-2\alpha\beta)^z$.
First we have to show that the elements in the Lemma are indeed fundamental. We inspect $1-p$ for all $p$ in this list:

```
Together[1 - {1, α, β, α β, (α - 1)/(α β - 1), (β - 1)/(α β - 1), -(α (β - 1))/(β (α - 1)), ((α - 1) (β - 1))/(1 - α β), (α (β - 1)^2)/(β (α β - 1)), (β (α - 1)^2)/(α (α β - 1))}]
```

$$\left\{0,\ 1-\alpha,\ 1-\beta,\ 1-\alpha\beta,\ \frac{-\alpha+\alpha\beta}{-1+\alpha\beta},\ \frac{-\beta+\alpha\beta}{-1+\alpha\beta},\ \frac{-\alpha-\beta+2\alpha\beta}{(-1+\alpha)\beta},\ \frac{-\alpha-\beta+2\alpha\beta}{-1+\alpha\beta},\ \frac{-\alpha-\beta+2\alpha\beta}{\beta(-1+\alpha\beta)},\ \frac{-\alpha-\beta+2\alpha\beta}{\alpha(-1+\alpha\beta)}\right\}$$

These are all powers of the generators, so all elements of
$\left\{1, \alpha, \beta, \alpha\beta, \frac{\alpha-1}{\alpha\beta-1}, \frac{\beta-1}{\alpha\beta-1}, -\frac{\alpha(\beta-1)}{\beta(\alpha-1)}, \frac{(\alpha-1)(\beta-1)}{1-\alpha\beta}, \frac{\alpha(\beta-1)^2}{\beta(\alpha\beta-1)}, \frac{\beta(\alpha-1)^2}{\alpha(\alpha\beta-1)}\right\}$
are indeed fundamental. It follows immediately that all their associates are fundamental.

To show that this list is indeed complete, we first try to bound the exponents of fundamental elements.

First note that exchanging $\alpha$ and $\beta$ is a partial-field automorphism, so any bound on $u$ implies a bound on $v$ and any bound on $w$ implies a bound on $x$. We will need a number of homomorphisms for our bounds. First we consider two homomorphisms to $\mathbb{H}_3$. We list the images of the generators of $\mathbb{H}_4$ in the following table.

```
Grid[{H4gens, Together[H4gens /. {α → 1/(1 - α), β → α^2 - α + 1}], Together[H4gens /. β → α/(α - 1)]},
 Dividers → {False, {False, True}}]
```

| | | | | | | |
|---|---|---|---|---|---|---|
| $-1$ | $\alpha$ | $\beta$ | $1-\alpha$ | $1-\beta$ | $-1+\alpha\beta$ | $\alpha+\beta-2\alpha\beta$ |
| $-1$ | $\frac{1}{1-\alpha}$ | $1-\alpha+\alpha^2$ | $\frac{\alpha}{-1+\alpha}$ | $\alpha-\alpha^2$ | $-\frac{\alpha^2}{-1+\alpha}$ | $\frac{\alpha^3}{-1+\alpha}$ |
| $-1$ | $\alpha$ | $\frac{\alpha}{-1+\alpha}$ | $1-\alpha$ | $-\frac{1}{-1+\alpha}$ | $\frac{1-\alpha+\alpha^2}{-1+\alpha}$ | $-\frac{\alpha^2}{-1+\alpha}$ |

Since $\alpha^2-\alpha+1$ occurs with exponent at most 1 and at least $-1$ in the fundamental elements of $\mathbb{H}_3$, it follows that $-1 \le u \le 1$, $-1 \le v \le 1$, $-1 \le y \le 1$. For bounds on $x$, $w$, $z$ we consider homomorphisms to $\mathbb{H}_2$.

```
αβims = {{-ⅈ, -ⅈ}, {-ⅈ, ⅈ (1 - ⅈ) / 2}, {(1 - ⅈ) / 2, ⅈ},
   {ⅈ (1 - ⅈ) / 2, -ⅈ}, {1 - ⅈ, ⅈ (1 - ⅈ)}, {1 - ⅈ, 1 / 2}, {ⅈ (1 - ⅈ), 1 - ⅈ},
   {ⅈ (1 - ⅈ), 1 / 2}, {ⅈ, (1 - ⅈ) / 2}, {ⅈ, ⅈ}, {1 / 2, 1 - ⅈ}, {1 / 2, ⅈ (1 - ⅈ)}};
```

```
imtable = Table[H4gens /. {α → αβims[[k, 1]], β → αβims[[k, 2]]}, {k, 1, Length[αβims]}];
```

```
Grid[Join[{H4gens}, imtable], Dividers → {False, {False, True}}]
```

| $-1$ | $\alpha$ | $\beta$ | $1-\alpha$ | $1-\beta$ | $-1+\alpha\beta$ | $\alpha+\beta-2\alpha\beta$ |
|---|---|---|---|---|---|---|
| $-1$ | $-i$ | $-i$ | $1+i$ | $1+i$ | $-2$ | $2-2i$ |
| $-1$ | $-i$ | $\frac{1}{2}+\frac{i}{2}$ | $1+i$ | $\frac{1}{2}-\frac{i}{2}$ | $-\frac{1}{2}-\frac{i}{2}$ | $-\frac{1}{2}+\frac{i}{2}$ |
| $-1$ | $\frac{1}{2}-\frac{i}{2}$ | $i$ | $\frac{1}{2}+\frac{i}{2}$ | $1-i$ | $-\frac{1}{2}+\frac{i}{2}$ | $-\frac{1}{2}-\frac{i}{2}$ |
| $-1$ | $\frac{1}{2}+\frac{i}{2}$ | $-i$ | $\frac{1}{2}-\frac{i}{2}$ | $1+i$ | $-\frac{1}{2}-\frac{i}{2}$ | $-\frac{1}{2}+\frac{i}{2}$ |
| $-1$ | $1-i$ | $1+i$ | $i$ | $-i$ | $1$ | $-2$ |
| $-1$ | $1-i$ | $\frac{1}{2}$ | $i$ | $\frac{1}{2}$ | $-\frac{1}{2}-\frac{i}{2}$ | $\frac{1}{2}$ |
| $-1$ | $1+i$ | $1-i$ | $-i$ | $i$ | $1$ | $-2$ |
| $-1$ | $1+i$ | $\frac{1}{2}$ | $-i$ | $\frac{1}{2}$ | $-\frac{1}{2}+\frac{i}{2}$ | $\frac{1}{2}$ |
| $-1$ | $i$ | $\frac{1}{2}-\frac{i}{2}$ | $1-i$ | $\frac{1}{2}+\frac{i}{2}$ | $-\frac{1}{2}+\frac{i}{2}$ | $-\frac{1}{2}-\frac{i}{2}$ |
| $-1$ | $i$ | $i$ | $1-i$ | $1-i$ | $-2$ | $2+2i$ |
| $-1$ | $\frac{1}{2}$ | $1-i$ | $\frac{1}{2}$ | $i$ | $-\frac{1}{2}-\frac{i}{2}$ | $\frac{1}{2}$ |
| $-1$ | $\frac{1}{2}$ | $1+i$ | $\frac{1}{2}$ | $-i$ | $-\frac{1}{2}+\frac{i}{2}$ | $\frac{1}{2}$ |

All elements of this table are elements of $\mathbb{H}_2$, so these are indeed the images of partial-field homomorphisms. Recall that each fundamental element of $\mathbb{H}_2$ has a norm between $\frac{1}{2}$ and 2. By taking the base-2 logarithm of the norms we obtain linear bounds on the exponents:

```
Simplify[Log[2, Abs[imtable]]]
```

$$\left\{\left\{0, 0, 0, \frac{1}{2}, \frac{1}{2}, 1, \frac{3}{2}\right\}, \left\{0, 0, -\frac{1}{2}, \frac{1}{2}, -\frac{1}{2}, -\frac{1}{2}, -\frac{1}{2}\right\}, \left\{0, -\frac{1}{2}, 0, -\frac{1}{2}, \frac{1}{2}, -\frac{1}{2}, -\frac{1}{2}\right\},\right.$$
$$\left\{0, -\frac{1}{2}, 0, -\frac{1}{2}, \frac{1}{2}, -\frac{1}{2}, -\frac{1}{2}\right\}, \left\{0, \frac{1}{2}, \frac{1}{2}, 0, 0, 0, 1\right\}, \left\{0, \frac{1}{2}, -1, 0, -1, -\frac{1}{2}, -1\right\},$$
$$\left\{0, \frac{1}{2}, \frac{1}{2}, 0, 0, 0, 1\right\}, \left\{0, \frac{1}{2}, -1, 0, -1, -\frac{1}{2}, -1\right\}, \left\{0, 0, -\frac{1}{2}, \frac{1}{2}, -\frac{1}{2}, -\frac{1}{2}, -\frac{1}{2}\right\},$$
$$\left.\left\{0, 0, 0, \frac{1}{2}, \frac{1}{2}, 1, \frac{3}{2}\right\}, \left\{0, -1, \frac{1}{2}, -1, 0, -\frac{1}{2}, -1\right\}, \left\{0, -1, \frac{1}{2}, -1, 0, -\frac{1}{2}, -1\right\}\right\}$$

```
constraints = Join[-1 ≤ #.{0, u, v, w, x, y, z} ≤ 1 & /@ %, {-1 ≤ u ≤ 1, -1 ≤ v ≤ 1, -1 ≤ y ≤ 1}]
```

$$\Big\{-1 \le \frac{w}{2} + \frac{x}{2} + y + \frac{3z}{2} \le 1,\ -1 \le -\frac{v}{2} + \frac{w}{2} - \frac{x}{2} - \frac{y}{2} - \frac{z}{2} \le 1,\ -1 \le -\frac{u}{2} - \frac{w}{2} + \frac{x}{2} - \frac{y}{2} - \frac{z}{2} \le 1,$$
$$-1 \le -\frac{u}{2} - \frac{w}{2} + \frac{x}{2} - \frac{y}{2} - \frac{z}{2} \le 1,\ -1 \le \frac{u}{2} + \frac{v}{2} + z \le 1,\ -1 \le \frac{u}{2} - v - x - \frac{y}{2} - z \le 1,\ -1 \le \frac{u}{2} + \frac{v}{2} + z \le 1,$$
$$-1 \le \frac{u}{2} - v - x - \frac{y}{2} - z \le 1,\ -1 \le -\frac{v}{2} + \frac{w}{2} - \frac{x}{2} - \frac{y}{2} - \frac{z}{2} \le 1,\ -1 \le \frac{w}{2} + \frac{x}{2} + y + \frac{3z}{2} \le 1,$$
$$-1 \le -u + \frac{v}{2} - w - \frac{y}{2} - z \le 1,\ -1 \le -u + \frac{v}{2} - w - \frac{y}{2} - z \le 1,\ -1 \le u \le 1,\ -1 \le v \le 1,\ -1 \le y \le 1\Big\}$$

Now we find the maximum and minimum value for $u$, $v$, $w$, $x$, $y$, $z$ subject to the linear constraints that we obtained.

```
G = {u, v, w, x, y, z, -u, -v, -w, -x, -y, -z};
```

```
Table[{G[[i]], Maximize[{G[[i]], And @@ constraints}, {u, v, w, x, y, z}, Reals][[1]]},
 {i, 1, Length[G]}]
```

```
{{u, 1}, {v, 1}, {w, 3}, {x, 3}, {y, 1}, {z, 2},
 {-u, 1}, {-v, 1}, {-w, 3}, {-x, 3}, {-y, 1}, {-z, 2}}
```

In other words, $-3 \le w \le 3$, $-3 \le x \le 3$ and $-2 \le z \le 2$. Now we have reduced the set of possible exponents of fundamental elements to a finite list.

```
candidateH4funExps = Flatten[Table[{s, u, v, w, x, y, z}, {s, 0, 1},
    {u, -1, 1}, {v, -1, 1}, {w, -3, 3}, {x, -3, 3}, {y, -1, 1}, {z, -2, 2}], 6];
```

```
Length[%]
```

```
13 230
```

To find out which elements are fundamental we map the candidates to GF($p$) for some large $p$, and choose the images of $\alpha$, $\beta$ so that this map is injective. Then we select the elements $x$ for which $1 - x$ is also in the list. To avoid dealing with fractions we reconstruct the list using the PowerMod function:

```
p = Prime[10 000 000]
```

```
179 424 673
```

```
aIm = 11; bIm = 19;
```

```
H4genIms = Mod[H4gens /. {α → aIm, β → bIm}, p]
```

```
{179 424 672, 11, 19, 179 424 663, 179 424 655, 208, 179 424 285}
```

```
candidateH4funIms =
  {0} ⋃ (Mod[Times @@ PowerMod[H4genIms, #, p], p] & /@ candidateH4funExps);
```

This map is indeed injective (the union symbol removes duplicates; note that we have added 0 so the list has one element more than the list of exponents)

```
Length[%]
```

```
13 231
```

Select the elements $x$ for which $1 - x$ is also in the list:

```
H4funIms = candidateH4funIms ⋂ Mod[1 - candidateH4funIms, p];
```

```
Length[%]
```

```
56
```

This list of remaining candidates necessarily contains the images of all known fundamental elements. Since the number of remaining elements equals the number of known fundamental elements, we conclude that those are, indeed, all. This completes the proof of Lemma 23. We record the corresponding exponents for use below. First we recompute the images, this time not taking the union to avoid sorting.

```
candidateH4funImsNoZero = Mod[Times @@ PowerMod[H4genIms, #, p], p] & /@ candidateH4funExps;
```

Next we find the positions of the images of fundamental elements in this list, and take the exponents with corresponding indices. We remove the exponent corresponding to 1.

```
H4funExps = Complement[candidateH4funExps[[
    Flatten[Position[candidateH4funImsNoZero, #] & /@ H4funIms]]], {{0, 0, 0, 0, 0, 0, 0}}];
```

```
Length[H4funExps]
```

```
54
```

## The automorphism group

In this section we prove the following result:

**Lemma 24.** The automorphism group of $\mathbb{H}_4$ is isomorphic to $S_4$.

The automorphisms correspond with coordinate permutations in GF(5) × GF(5) × GF(5) × GF(5).

The partial-field automorphism group of $\mathbb{H}_4$ is determined uniquely by the images of $\alpha$ and $\beta$. This is easy to see since each automorphism maps 0 to 0, 1 to 1, and $-1$ to $-1$, so the real number line is fixed. Note that $\alpha$ and $\beta$ have to be mapped to fundamental elements, and those fundamental elements have to be distinct, so the group will be finite. Note that we only have to check if $\left\{1, \alpha, \beta, \alpha\,\beta, \frac{\alpha-1}{\alpha\,\beta-1}, \frac{\beta-1}{\alpha\,\beta-1}, -\frac{\alpha\,(\beta-1)}{\beta\,(\alpha-1)}, \frac{(\alpha-1)\,(\beta-1)}{1-\alpha\,\beta}, \frac{\alpha\,(\beta-1)^2}{\beta\,(\alpha\,\beta-1)}, \frac{\beta\,(\alpha-1)^2}{\alpha\,(\alpha\,\beta-1)}\right\}$ map to fundamental elements: the remainder follows automatically. To speed matters up, we first check if the images of fundamental elements under the homomorphism to GF($p$) from the previous section map to fundamental elements. Only on the remaining sets will we test for equality of the rational expressions.

Note that 0 and 1 are mapped to themselves, so we do not have to test for those.

```
nonzerooneH4funIms = Mod[Times @@ PowerMod[H4genIms, #, p], p] & /@ H4funExps;
```

```
Length[%]
```

```
54
```

We recompute the next line to ensure the ordering is identical to the ordering of the images:

```
nonzerooneH4funs = (Times @@ (H4gens^#)) & /@ H4funExps;
```

```
Length[%]
```

```
54
```

```
validImage[aIm_, bIm_] := Module[{nonzerooneH4funImsp, H4genIms},
   (* Check if the generators get mapped to nonzero numbers *)
   H4genIms = Mod[H4gens /. {α → aIm, β → bIm}, p];
   If[MemberQ[H4genIms, 0], Return[False]];
   (* Construct the list of images of fundamental elements *)
   nonzerooneH4funImsp := Mod[Times @@ PowerMod[H4genIms, #, p], p] & /@ H4funExps;
   (* Check if this list equals the old list of images,
   and return True or False accordingly *)
   Length[Complement[nonzerooneH4funIms, nonzerooneH4funImsp]] == 0
  ];
```

```
Table[validImage[nonzerooneH4funIms[[x]], nonzerooneH4funIms[[y]]],
  {x, 1, Length[nonzerooneH4funIms]}, {y, 1, Length[nonzerooneH4funIms]}];
```

```
Position[%, True]
```

```
{{13, 21}, {13, 37}, {14, 16}, {14, 36}, {15, 20}, {15, 40}, {16, 14}, {16, 36},
 {20, 15}, {20, 40}, {21, 13}, {21, 37}, {36, 14}, {36, 16}, {37, 13}, {37, 21},
 {40, 15}, {40, 20}, {45, 48}, {45, 52}, {48, 45}, {48, 52}, {52, 45}, {52, 48}}
```

```
H4autImages = nonzerooneH4funs[[#]] & /@ %
```

$$\left\{\left\{1-\beta,\ \frac{\alpha\,(1-\beta)}{\alpha+\beta-2\,\alpha\,\beta}\right\},\ \left\{1-\beta,\ -\frac{-1+\alpha\,\beta}{1-\beta}\right\},\ \left\{1-\alpha,\ \frac{(1-\alpha)\,\beta}{\alpha+\beta-2\,\alpha\,\beta}\right\},\ \left\{1-\alpha,\ -\frac{-1+\alpha\,\beta}{1-\alpha}\right\},\ \{\beta,\ \alpha\},\right.$$
$$\left\{\beta,\ -\frac{-1+\alpha\,\beta}{\alpha+\beta-2\,\alpha\,\beta}\right\},\ \left\{\frac{(1-\alpha)\,\beta}{\alpha+\beta-2\,\alpha\,\beta},\ 1-\alpha\right\},\ \left\{\frac{(1-\alpha)\,\beta}{\alpha+\beta-2\,\alpha\,\beta},\ -\frac{-1+\alpha\,\beta}{1-\alpha}\right\},\ \{\alpha,\ \beta\},\ \left\{\alpha,\ -\frac{-1+\alpha\,\beta}{\alpha+\beta-2\,\alpha\,\beta}\right\},$$
$$\left\{\frac{\alpha\,(1-\beta)}{\alpha+\beta-2\,\alpha\,\beta},\ 1-\beta\right\},\ \left\{\frac{\alpha\,(1-\beta)}{\alpha+\beta-2\,\alpha\,\beta},\ -\frac{-1+\alpha\,\beta}{1-\beta}\right\},\ \left\{-\frac{-1+\alpha\,\beta}{1-\alpha},\ 1-\alpha\right\},\ \left\{-\frac{-1+\alpha\,\beta}{1-\alpha},\ \frac{(1-\alpha)\,\beta}{\alpha+\beta-2\,\alpha\,\beta}\right\},$$
$$\left\{-\frac{-1+\alpha\,\beta}{1-\beta},\ 1-\beta\right\},\ \left\{-\frac{-1+\alpha\,\beta}{1-\beta},\ \frac{\alpha\,(1-\beta)}{\alpha+\beta-2\,\alpha\,\beta}\right\},\ \left\{-\frac{-1+\alpha\,\beta}{\alpha+\beta-2\,\alpha\,\beta},\ \beta\right\},\ \left\{-\frac{-1+\alpha\,\beta}{\alpha+\beta-2\,\alpha\,\beta},\ \alpha\right\},$$
$$\left\{-\frac{(1-\alpha)\,(1-\beta)}{\alpha+\beta-2\,\alpha\,\beta},\ -\frac{(1-\alpha)\,\beta}{1-\beta}\right\},\ \left\{-\frac{(1-\alpha)\,(1-\beta)}{\alpha+\beta-2\,\alpha\,\beta},\ -\frac{\alpha\,(1-\beta)}{1-\alpha}\right\},\ \left\{-\frac{(1-\alpha)\,\beta}{1-\beta},\ -\frac{(1-\alpha)\,(1-\beta)}{\alpha+\beta-2\,\alpha\,\beta}\right\},$$
$$\left.\left\{-\frac{(1-\alpha)\,\beta}{1-\beta},\ -\frac{\alpha\,(1-\beta)}{1-\alpha}\right\},\ \left\{-\frac{\alpha\,(1-\beta)}{1-\alpha},\ -\frac{(1-\alpha)\,(1-\beta)}{\alpha+\beta-2\,\alpha\,\beta}\right\},\ \left\{-\frac{\alpha\,(1-\beta)}{1-\alpha},\ -\frac{(1-\alpha)\,\beta}{1-\beta}\right\}\right\}$$

```
Length[%]
```

```
24
```

The images of the representative set of fundamental elements:

```
candidateAuts =
  Table[{α, β, α β, (α - 1)/(α β - 1), (β - 1)/(α β - 1), -(α (β - 1))/(β (α - 1)), ((α - 1) (β - 1))/(1 - α β), (α (β - 1)^2)/(β (α β - 1)), (β (α - 1)^2)/(α (α β - 1))} /.
    {α → H4autImages[[x, 1]], β → H4autImages[[x, 2]]}, {x, 1, Length[H4autImages]}];
```

For these remaining ones we will do the (much slower) full check:

```
isH4fun[x_] := occursInList[nonzerooneH4funs, x];
areH4funs[x_List] := And @@ (isH4fun /@ x);
```

```
H4auts = Select[candidateAuts, areH4funs];
```

```
Length[%]
```

```
24
```

Conclusion: the automorphism group of $\mathbb{H}_4$ has 24 elements and is isomorphic to $S_4$, the symmetric group on 4 symbols. Indeed, the partial-field homomorphism $\phi 4$ defined above induces a $1-1$ correspondence between the automorphisms of $\mathbb{H}_4$ and coordinate permutations of $\mathrm{GF}(5)\times\mathrm{GF}(5)\times\mathrm{GF}(5)\times\mathrm{GF}(5)$.

## Representations of $U_{2,5}$

We prove the following fact:

**Lemma 25.** Let $M$ be a 3-connected, $\mathbb{H}_4$-representable matroid. If $M$ has a $U_{2,5}$- or $U_{3,5}$-minor then $M$ has at least four inequivalent representations over GF(5).

To prove this it suffices to show that each $\mathbb{H}_4$-representation of $U_{2,5}$ gives rise to four inequivalent representations over GF(5). Since any $\mathbb{H}_4$-representation matrix of $M$ must contain one of these as a minor, the result then follows. First we enumerate all such representations. By normalizing and suppressing the identity matrix at the front, we see

```
Clear[p];
```

```
A = {{1, 1, 1}, {1, p, q}}; MatrixForm[A]
```

$$\begin{pmatrix} 1 & 1 & 1 \\ 1 & p & q \end{pmatrix}$$

This is an $\mathbb{H}_4$-matrix representing $U_{2,5}$ if and only if $p$, $q$, $\frac{p}{q}$ are fundamental elements, $p$ and $q$ are not equal to 0 or 1, and $p \neq q$. Moreover, two such matrices are equivalent if and only if they are equal.

Create a table of all candidate pairs $p$, $q$. We suppress the output:

```
candidatepqPairs = Flatten[Table[{nonzerooneH4funs[[x]], nonzerooneH4funs[[y]]},
     {x, 1, Length[nonzerooneH4funs]}, {y, 1, x - 1}], 1] ⋃
   Flatten[Table[{nonzerooneH4funs[[x]], nonzerooneH4funs[[y]]},
     {x, 1, Length[nonzerooneH4funs]}, {y, x + 1, Length[nonzerooneH4funs]}], 1];
```

Filter out those for which $p/q$ is not fundamental or equal to 1:

```
pqPairs = Select[candidatepqPairs, occursInList[nonzerooneH4funs, #[[1]] / #[[2]]] &];
```

```
Length[%]
```

```
360
```

To show that each of these gives rise to four inequivalent representations over GF(5) we consider the partial-field homomorphism $\phi 4$.

Let $\xi_i : \text{GF}(5) \times \text{GF}(5) \times \text{GF}(5) \times \text{GF}(5) \to \text{GF}(5)$ be the projection onto the $i$th coordinate. We have to show that $\xi_i(\phi(A))$ is not equivalent to $\xi_j(\phi(A))$ for all matrices $A$ computed above. We only need to look at the pairs $p$, $q$. These are their images under $\phi 4$:

```
GF5t4pqpairs = Map[ϕ4, pqPairs, {2}];
```

Indeed, each pair of representations is inequivalent. A concise way of checking this is to show that for all $1 \leq i < j \leq 4$, $\left| \xi_i(p) - \xi_j(p) \right| + \left| \xi_i(q) - \xi_j(q) \right| > 0$.

```
Map[{Abs[#[[1, 1]] - #[[1, 2]]] + Abs[#[[2, 1]] - #[[2, 2]]],
    Abs[#[[1, 1]] - #[[1, 3]]] + Abs[#[[2, 1]] - #[[2, 3]]],
    Abs[#[[1, 1]] - #[[1, 4]]] + Abs[#[[2, 1]] - #[[2, 4]]],
    Abs[#[[1, 2]] - #[[1, 3]]] + Abs[#[[2, 2]] - #[[2, 3]]],
    Abs[#[[1, 2]] - #[[1, 4]]] + Abs[#[[2, 2]] - #[[2, 4]]],
    Abs[#[[1, 3]] - #[[1, 4]]] + Abs[#[[2, 3]] - #[[2, 4]]]} &, %, {1}];
```

```
Position[%, 0]
```

```
{}
```

This completes the proof of Lemma 25.

## Lifting

We prove the following fact:

**Lemma 26.** Let $M$ be a 3-connected matroid with at least four inequivalent representations over GF(5). Then $M$ is representable over $\mathbb{H}_4$.

Consider a $\text{GF}(5)\times\text{GF}(5)\times\text{GF}(5)\times\text{GF}(5)$-matrix $A$ representing $M$ such that the four projections $\xi_i(A)$ are pairwise inequivalent. We want to apply the Lift Theorem, so we start by constructing a lifting function for Cr($A$). First we need to find out what Cr($A$) is.

**Claim 26.1.** No element of Cr($A$) is of the form $(x, x, x, y)$ for $x \notin \{0, 1\}$.

Proof. Assume that $(x, x, x, y) \in \text{Cr}(A)$ for some $x \in \{2, 3, 4\}$. Let $\mathbb{P}'$ be the sub-partial field of $\text{GF}(5)\times\text{GF}(5)\times\text{GF}(5)\times\text{GF}(5)$ such that the first three coordinates of all elements are equal. Since $A$ is not a $\mathbb{P}'$-matrix, the Confinement Theorem, in particular Corollary 15, implies that either $A$ or $A^T$ contains a minor

$$D = \begin{pmatrix} (1, 1, 1, 1) & (1, 1, 1, 1) & (1, 1, 1, 1) \\ (1, 1, 1, 1) & (x, x, x, y) & (p, q, r, s) \end{pmatrix},$$

with not all three of $p$, $q$, $r$ equal. There are only two elements $u$ of GF(5) such that $\begin{pmatrix} 1 & 1 & 1 \\ 1 & x & u \end{pmatrix}$ represents $U_{2,5}$. Hence two of the first three projections of $D$ must be equivalent. But then the Stabilizer Theorem, in particular Lemma 8, implies that two of the four projections of $A$ are equivalent, a contradiction. □

By considering all permutations of the coordinates we are left with the following cross ratios:

$F' := \mathcal{F}(\text{GF}(5)\times\text{GF}(5)\times\text{GF}(5)\times\text{GF}(5)) \setminus \{x : \text{at least 3 coordinates of } x \text{ are equal},\ x \notin \{(0, 0, 0, 0), (1, 1, 1, 1)\}\}$.

```
Fp = Complement[{{0, 0, 0, 0}, {1, 1, 1, 1}} ∪
   Flatten[Table[{w, x, y, z}, {w, 2, 4}, {x, 2, 4}, {y, 2, 4}, {z, 2, 4}], 3], Flatten[
   Table[{{x, x, x, y}, {x, x, y, x}, {x, y, x, x}, {y, x, x, x}}, {x, 2, 4}, {y, 2, 4}], 2]]
```

```
{{0, 0, 0, 0}, {1, 1, 1, 1}, {2, 2, 3, 3}, {2, 2, 3, 4}, {2, 2, 4, 3}, {2, 2, 4, 4}, {2, 3, 2, 3},
 {2, 3, 2, 4}, {2, 3, 3, 2}, {2, 3, 3, 4}, {2, 3, 4, 2}, {2, 3, 4, 3}, {2, 3, 4, 4}, {2, 4, 2, 3},
 {2, 4, 2, 4}, {2, 4, 3, 2}, {2, 4, 3, 3}, {2, 4, 3, 4}, {2, 4, 4, 2}, {2, 4, 4, 3},
 {3, 2, 2, 3}, {3, 2, 2, 4}, {3, 2, 3, 2}, {3, 2, 3, 4}, {3, 2, 4, 2}, {3, 2, 4, 3},
 {3, 2, 4, 4}, {3, 3, 2, 2}, {3, 3, 2, 4}, {3, 3, 4, 2}, {3, 3, 4, 4}, {3, 4, 2, 2},
 {3, 4, 2, 3}, {3, 4, 2, 4}, {3, 4, 3, 2}, {3, 4, 3, 4}, {3, 4, 4, 2}, {3, 4, 4, 3},
 {4, 2, 2, 3}, {4, 2, 2, 4}, {4, 2, 3, 2}, {4, 2, 3, 3}, {4, 2, 3, 4}, {4, 2, 4, 2},
 {4, 2, 4, 3}, {4, 3, 2, 2}, {4, 3, 2, 3}, {4, 3, 2, 4}, {4, 3, 3, 2}, {4, 3, 3, 4},
 {4, 3, 4, 2}, {4, 3, 4, 3}, {4, 4, 2, 2}, {4, 4, 2, 3}, {4, 4, 3, 2}, {4, 4, 3, 3}}
```

```
Length[%]
```

```
56
```

The restriction of $\phi 4$ to $\mathcal{F}(\mathbb{H}_4)$ is a bijection between $\mathcal{F}(\mathbb{H}_4)$ and $F'$:

```
Sort[ϕ4 /@ H4funs]
```

```
{{0, 0, 0, 0}, {1, 1, 1, 1}, {2, 2, 3, 3}, {2, 2, 3, 4}, {2, 2, 4, 3}, {2, 2, 4, 4}, {2, 3, 2, 3},
 {2, 3, 2, 4}, {2, 3, 3, 2}, {2, 3, 3, 4}, {2, 3, 4, 2}, {2, 3, 4, 3}, {2, 3, 4, 4}, {2, 4, 2, 3},
 {2, 4, 2, 4}, {2, 4, 3, 2}, {2, 4, 3, 3}, {2, 4, 3, 4}, {2, 4, 4, 2}, {2, 4, 4, 3},
 {3, 2, 2, 3}, {3, 2, 2, 4}, {3, 2, 3, 2}, {3, 2, 3, 4}, {3, 2, 4, 2}, {3, 2, 4, 3},
 {3, 2, 4, 4}, {3, 3, 2, 2}, {3, 3, 2, 4}, {3, 3, 4, 2}, {3, 3, 4, 4}, {3, 4, 2, 2},
 {3, 4, 2, 3}, {3, 4, 2, 4}, {3, 4, 3, 2}, {3, 4, 3, 4}, {3, 4, 4, 2}, {3, 4, 4, 3},
 {4, 2, 2, 3}, {4, 2, 2, 4}, {4, 2, 3, 2}, {4, 2, 3, 3}, {4, 2, 3, 4}, {4, 2, 4, 2},
 {4, 2, 4, 3}, {4, 3, 2, 2}, {4, 3, 2, 3}, {4, 3, 2, 4}, {4, 3, 3, 2}, {4, 3, 3, 4},
 {4, 3, 4, 2}, {4, 3, 4, 3}, {4, 4, 2, 2}, {4, 4, 2, 3}, {4, 4, 3, 2}, {4, 4, 3, 3}}
```

It follows immediately that the function $F' \to \mathbb{H}_4$ that is the inverse of $\phi$ is a lifting function.

```
up4[w_, x_, y_, z_] := H4funs[[Position[ϕ4 /@ H4funs, {w, x, y, z}][[1, 1]]]];
up4[{w_, x_, y_, z_}] := up4[w, x, y, z];
```

Some examples:

```
up4[1, 1, 1, 1]
```

```
1
```

```
up4[2, 4, 3, 4]
```

$$\frac{\beta}{-1+\beta}$$

Now suppose that $A$ has no global lift. Then the Lift Theorem (Theorem 12) tells us that $A$ or $A^T$ has a minor of the form

$$\begin{pmatrix} 1 & 1 & 0 & 1 \\ 1 & 0 & 1 & 1 \\ 0 & 1 & 1 & 1 \end{pmatrix} \text{ or } \begin{pmatrix} 1 & 1 & 1 \\ 1 & p & q \end{pmatrix}, \text{ for some } p, q \in \mathrm{Cr}(A).$$

without a local lift. The first of these has a cross ratio (2, 2, 2, 2), so it does not occur as a minor of $A$, by Claim 26.1. Therefore we have to ensure that all $\mathrm{GF}(5)\times\mathrm{GF}(5)\times\mathrm{GF}(5)\times\mathrm{GF}(5)$-representations of $U_{2,5}$ that can occur in $A$ have a local lift. Let $D$ be such a minor. Suppose that $\xi_i(D)$ and $\xi_j(D)$ are equivalent. By the Stabilizer Theorem, and in particular Lemma 8, this implies that $\xi_i(A)$ and $\xi_j(A)$ are equivalent, a contradiction to our choice of $A$. Hence $\xi_i(D)$ and $\xi_j(D)$ are inequivalent. This gives us $6\times 5\times 4\times 3 = 360$ representations of $U_{2,5}$. Again we only list the pairs $(p, q)$ in

$$D = \begin{pmatrix} 1 & 1 & 1 \\ 1 & p & q \end{pmatrix}.$$

```
U25reps = { {2, 3}, {2, 4}, {3, 2}, {3, 4}, {4, 2}, {4, 3}}
```

```
{{2, 3}, {2, 4}, {3, 2}, {3, 4}, {4, 2}, {4, 3}}
```

```
U25repquads = Transpose /@ Permutations[U25reps, {4}];
```

```
Length[%]
```

```
360
```

These are the lifts (we suppress the output):

```
Map[up4, U25repquads, {2}];
```

```
Length[%]
```

```
360
```

And indeed, for each of these $p^{\uparrow} / q^{\uparrow}$ is a fundamental element:

```
Select[%%, Not[occursInList[nonzerooneH4funs, #[[1]] / #[[2]]]] &]
```

```
{}
```

It follows that $A$ does have a global lift. This completes the proof of Lemma 26.

# Hydra-5

## Definition and homomorphisms

The Hydra-5 partial field is defined as
$\mathbb{H}_5 = (\mathbb{Q}(\alpha, \beta, \gamma), \langle -1, \alpha, \beta, \gamma, 1-\alpha, 1-\beta, 1-\gamma, \alpha-\gamma, \gamma-\alpha\beta, (1-\gamma)-(1-\alpha)\beta \rangle)$,
where $\alpha$, $\beta$, $\gamma$ are indeterminates.

The generators of $\mathbb{H}_5$ are

```
H5gens = {-1, α, β, γ, 1 - α, 1 - β, 1 - γ, α - γ, γ - α β, (1 - γ) - (1 - α) β};
```

The fundamental elements of $\mathbb{H}_5$ are

$$\text{H5funs} = \text{Together}\left[\text{associates}\left[\left\{1, \alpha, \beta, \gamma, \frac{\alpha\beta}{\gamma}, \frac{\alpha}{\gamma}, \frac{(1-\alpha)\gamma}{\gamma-\alpha}, \frac{(\alpha-1)\beta}{\gamma-1}, \frac{\alpha-1}{\gamma-1}, \frac{\gamma-\alpha}{\gamma-\alpha\beta}, \frac{(\beta-1)(\gamma-1)}{\beta(\gamma-\alpha)}, \frac{\beta(\gamma-\alpha)}{\gamma-\alpha\beta}, \frac{(\alpha-1)(\beta-1)}{\gamma-\alpha}, \frac{\beta(\gamma-\alpha)}{(1-\gamma)(\gamma-\alpha\beta)}, \frac{(1-\alpha)(\gamma-\alpha\beta)}{\gamma-\alpha}, \frac{1-\beta}{\gamma-\alpha\beta}\right\}\right]\right];$$

```
Length[%]
```

```
92
```

The fundamental elements other than 0, 1 play a special role:

```
nonzerooneH5funs = Complement[H5funs, {0, 1}];
```

There are six homomorphisms $\mathbb{H}_5 \to \mathrm{GF}(5)$. We collect them in the partial-field homomorphism $\phi 5 : \mathbb{H}_5 \to \mathrm{GF}(5) \times \mathrm{GF}(5) \times \mathrm{GF}(5) \times \mathrm{GF}(5) \times \mathrm{GF}(5)$ defined by $\phi 5(\alpha) = (4, 3, 3, 4, 2, 2)$, $\phi 5(\beta) = (3, 2, 2, 4, 3, 4)$, and $\phi 5(\gamma) = (3, 2, 4, 2, 3, 4)$.

```
ϕ5[0] = {0, 0, 0, 0, 0, 0};
ϕ5[1] = {1, 1, 1, 1, 1, 1};
ϕ5[-1] = {4, 4, 4, 4, 4, 4};
ϕ5[x_] := toGF5[x /. {α → {4, 3, 3, 4, 2, 2}, β → {3, 2, 2, 4, 3, 4}, γ → {3, 2, 4, 2, 3, 4}}];
```

This is indeed a partial-field homomorphism:

```
ϕ5 /@ H5gens
```

```
{{4, 4, 4, 4, 4, 4}, {4, 3, 3, 4, 2, 2}, {3, 2, 2, 4, 3, 4},
 {3, 2, 4, 2, 3, 4}, {2, 3, 3, 2, 4, 4}, {3, 4, 4, 2, 3, 2},
 {3, 4, 2, 4, 3, 2}, {1, 1, 4, 2, 4, 3}, {1, 1, 3, 1, 2, 1}, {2, 3, 1, 1, 1, 1}}
```

## Fundamental elements

Our aim in this section is to verify the following statement:

**Lemma 27.** The fundamental elements of $\mathbb{H}_5$ are
$$\mathrm{Asc}\left\{1, \alpha, \beta, \gamma, \frac{\alpha\beta}{\gamma}, \frac{\alpha}{\gamma}, \frac{(1-\alpha)\gamma}{\gamma-\alpha}, \frac{(\alpha-1)\beta}{\gamma-1}, \frac{\alpha-1}{\gamma-1}, \frac{\gamma-\alpha}{\gamma-\alpha\beta}, \frac{(\beta-1)(\gamma-1)}{\beta(\gamma-\alpha)}, \frac{\beta(\gamma-\alpha)}{\gamma-\alpha\beta}, \frac{(\alpha-1)(\beta-1)}{\gamma-\alpha}, \frac{\beta(\gamma-\alpha)}{(1-\gamma)(\gamma-\alpha\beta)}, \frac{(1-\alpha)(\gamma-\alpha\beta)}{\gamma-\alpha}, \frac{1-\beta}{\gamma-\alpha\beta}\right\}.$$

Elements of this partial field are of the form
$$\pm\alpha^r \beta^s \gamma^t (1-\alpha)^u (1-\beta)^v (1-\gamma)^w (\alpha-\gamma)^x (\gamma-\alpha\beta)^y ((1-\gamma)-(1-\alpha)\beta)^z.$$
First we have to show that the elements in the Lemma are indeed fundamental. We inspect $1-p$ for all $p$ in this list:

```
Together[1 - {1, α, β, γ, αβ/γ, α/γ, (1 - α) γ/(γ - α), (α - 1) β/(γ - 1), (α - 1)/(γ - 1), (γ - α)/(γ - α β),
    (β - 1) (γ - 1)/(β (γ - α)), β (γ - α)/(γ - α β), (α - 1) (β - 1)/(γ - α), β (γ - α)/((1 - γ) (γ - α β)), (1 - α) (γ - α β)/(γ - α), (1 - β)/(γ - α β)}]
```

$$\left\{0,\ 1-\alpha,\ 1-\beta,\ 1-\gamma,\ \frac{-\alpha\beta+\gamma}{\gamma},\ \frac{-\alpha+\gamma}{\gamma},\ \frac{\alpha-\alpha\gamma}{\alpha-\gamma},\ \frac{-1+\beta-\alpha\beta+\gamma}{-1+\gamma},\ \frac{-\alpha+\gamma}{-1+\gamma},\ \frac{-\alpha+\alpha\beta}{\alpha\beta-\gamma},\ \frac{-1+\beta-\alpha\beta+\gamma}{\beta\,(-\alpha+\gamma)},\ \frac{-\gamma+\beta\gamma}{\alpha\beta-\gamma},\ \frac{1-\beta+\alpha\beta-\gamma}{\alpha-\gamma},\ -\frac{\gamma-\beta\gamma+\alpha\beta\gamma-\gamma^2}{(-1+\gamma)\ (-\alpha\beta+\gamma)},\ \frac{\alpha-\alpha\beta+\alpha^2\beta-\alpha\gamma}{\alpha-\gamma},\ \frac{1-\beta+\alpha\beta-\gamma}{\alpha\beta-\gamma}\right\}$$

These are all powers of the generators, so all elements of
$\left\{1,\ \alpha,\ \beta,\ \gamma,\ \frac{\alpha\beta}{\gamma},\ \frac{\alpha}{\gamma},\ \frac{(1-\alpha)\gamma}{\gamma-\alpha},\ \frac{(\alpha-1)\beta}{\gamma-1},\ \frac{\alpha-1}{\gamma-1},\ \frac{\gamma-\alpha}{\gamma-\alpha\beta},\ \frac{(\beta-1)(\gamma-1)}{\beta(\gamma-\alpha)},\ \frac{\beta(\gamma-\alpha)}{\gamma-\alpha\beta},\ \frac{(\alpha-1)(\beta-1)}{\gamma-\alpha},\ \frac{\beta(\gamma-\alpha)}{(1-\gamma)(\gamma-\alpha\beta)},\ \frac{(1-\alpha)(\gamma-\alpha\beta)}{\gamma-\alpha},\ \frac{1-\beta}{\gamma-\alpha\beta}\right\}$
are indeed fundamental. It follows immediately that all their associates are fundamental.

To show that this list is indeed complete, we first try to bound the exponents of fundamental elements.

We consider the following homomorphisms to $\mathbb{H}_2$.

```
αβγims = {{-1, 1 - i, i}, {-1, 1 + i, -i}, {-i, 1/2 - i/2, -1}, {-i, 1/2 - i/2, 1/2 - i/2}, {-i, i, 1 - i},
    {-i, i, i}, {1/2 - i/2, -1, 1/2 + i/2}, {1/2 - i/2, 1 - i, 1 - i}, {1/2 - i/2, 1 + i, -i}, {1/2 - i/2, 1/2, 1/2},
    {1/2 + i/2, -1, 1/2 - i/2}, {1/2 + i/2, 1 - i, i}, {1/2 + i/2, 1 + i, 1 + i}, {1/2 + i/2, 1/2, 1/2}, {1 - i, -i, -i},
    {1 - i, -i, 1 + i}, {1 - i, 1/2 + i/2, 1/2 - i/2}, {1 - i, 1/2 + i/2, 2}, {1 + i, 1/2 - i/2, 1/2 + i/2},
    {1 + i, 1/2 - i/2, 2}, {1 + i, i, 1 - i}, {1 + i, i, i}, {i, -i, -i}, {i, -i, 1 + i}, {i, 1/2 + i/2, -1},
    {i, 1/2 + i/2, 1/2 + i/2}, {1/2, 2, 1/2 - i/2}, {1/2, 2, 1/2 + i/2}, {2, 1 - i, 1 - i}, {2, 1 + i, 1 + i}};
```

Indeed, all generators get mapped to elements of $\mathbb{H}_2$:

```
imtable = Table[H5gens /. {α → αβγims[[k, 1]], β → αβγims[[k, 2]], γ → αβγims[[k, 3]]},
    {k, 1, Length[αβγims]}];
```

```
Grid[Join[{H5gens}, imtable], Dividers → {False, {False, True}}]
```

| $-1$ | $\alpha$ | $\beta$ | $\gamma$ | $1-\alpha$ | $1-\beta$ | $1-\gamma$ | $\alpha-\gamma$ | $-\alpha\beta+\gamma$ | $1-(1-\alpha)\beta-\gamma$ |
|---|---|---|---|---|---|---|---|---|---|
| $-1$ | $-1$ | $1-i$ | $i$ | $2$ | $i$ | $1-i$ | $-1-i$ | $1$ | $-1+i$ |
| $-1$ | $-1$ | $1+i$ | $-i$ | $2$ | $-i$ | $1+i$ | $-1+i$ | $1$ | $-1-i$ |
| $-1$ | $-i$ | $\frac{1}{2}-\frac{i}{2}$ | $-1$ | $1+i$ | $\frac{1}{2}+\frac{i}{2}$ | $2$ | $1-i$ | $-\frac{1}{2}+\frac{i}{2}$ | $1$ |
| $-1$ | $-i$ | $\frac{1}{2}-\frac{i}{2}$ | $\frac{1}{2}-\frac{i}{2}$ | $1+i$ | $\frac{1}{2}+\frac{i}{2}$ | $\frac{1}{2}+\frac{i}{2}$ | $-\frac{1}{2}-\frac{i}{2}$ | $1$ | $-\frac{1}{2}+\frac{i}{2}$ |
| $-1$ | $-i$ | $i$ | $1-i$ | $1+i$ | $1-i$ | $i$ | $-1$ | $-i$ | $1$ |
| $-1$ | $-i$ | $i$ | $i$ | $1+i$ | $1-i$ | $1-i$ | $-2i$ | $-1+i$ | $2-2i$ |
| $-1$ | $\frac{1}{2}-\frac{i}{2}$ | $-1$ | $\frac{1}{2}+\frac{i}{2}$ | $\frac{1}{2}+\frac{i}{2}$ | $2$ | $\frac{1}{2}-\frac{i}{2}$ | $-i$ | $1$ | $1$ |
| $-1$ | $\frac{1}{2}-\frac{i}{2}$ | $1-i$ | $1-i$ | $\frac{1}{2}+\frac{i}{2}$ | $i$ | $i$ | $-\frac{1}{2}+\frac{i}{2}$ | $1$ | $-1+i$ |
| $-1$ | $\frac{1}{2}-\frac{i}{2}$ | $1+i$ | $-i$ | $\frac{1}{2}+\frac{i}{2}$ | $-i$ | $1+i$ | $\frac{1}{2}+\frac{i}{2}$ | $-1-i$ | $1$ |
| $-1$ | $\frac{1}{2}-\frac{i}{2}$ | $\frac{1}{2}$ | $\frac{1}{2}$ | $\frac{1}{2}+\frac{i}{2}$ | $\frac{1}{2}$ | $\frac{1}{2}$ | $-\frac{i}{2}$ | $\frac{1}{4}+\frac{i}{4}$ | $\frac{1}{4}-\frac{i}{4}$ |
| $-1$ | $\frac{1}{2}+\frac{i}{2}$ | $-1$ | $\frac{1}{2}-\frac{i}{2}$ | $\frac{1}{2}-\frac{i}{2}$ | $2$ | $\frac{1}{2}+\frac{i}{2}$ | $i$ | $1$ | $1$ |
| $-1$ | $\frac{1}{2}+\frac{i}{2}$ | $1-i$ | $i$ | $\frac{1}{2}-\frac{i}{2}$ | $i$ | $1-i$ | $\frac{1}{2}-\frac{i}{2}$ | $-1+i$ | $1$ |
| $-1$ | $\frac{1}{2}+\frac{i}{2}$ | $1+i$ | $1+i$ | $\frac{1}{2}-\frac{i}{2}$ | $-i$ | $-i$ | $-\frac{1}{2}-\frac{i}{2}$ | $1$ | $-1-i$ |
| $-1$ | $\frac{1}{2}+\frac{i}{2}$ | $\frac{1}{2}$ | $\frac{1}{2}$ | $\frac{1}{2}-\frac{i}{2}$ | $\frac{1}{2}$ | $\frac{1}{2}$ | $\frac{i}{2}$ | $\frac{1}{4}-\frac{i}{4}$ | $\frac{1}{4}+\frac{i}{4}$ |
| $-1$ | $1-i$ | $-i$ | $-i$ | $i$ | $1+i$ | $1+i$ | $1$ | $1$ | $i$ |
| $-1$ | $1-i$ | $-i$ | $1+i$ | $i$ | $1+i$ | $-i$ | $-2i$ | $2+2i$ | $-1-i$ |
| $-1$ | $1-i$ | $\frac{1}{2}+\frac{i}{2}$ | $\frac{1}{2}-\frac{i}{2}$ | $i$ | $\frac{1}{2}-\frac{i}{2}$ | $\frac{1}{2}+\frac{i}{2}$ | $\frac{1}{2}-\frac{i}{2}$ | $-\frac{1}{2}-\frac{i}{2}$ | $1$ |
| $-1$ | $1-i$ | $\frac{1}{2}+\frac{i}{2}$ | $2$ | $i$ | $\frac{1}{2}-\frac{i}{2}$ | $-1$ | $-1-i$ | $1$ | $-\frac{1}{2}-\frac{i}{2}$ |
| $-1$ | $1+i$ | $\frac{1}{2}-\frac{i}{2}$ | $\frac{1}{2}+\frac{i}{2}$ | $-i$ | $\frac{1}{2}+\frac{i}{2}$ | $\frac{1}{2}-\frac{i}{2}$ | $\frac{1}{2}+\frac{i}{2}$ | $-\frac{1}{2}+\frac{i}{2}$ | $1$ |
| $-1$ | $1+i$ | $\frac{1}{2}-\frac{i}{2}$ | $2$ | $-i$ | $\frac{1}{2}+\frac{i}{2}$ | $-1$ | $-1+i$ | $1$ | $-\frac{1}{2}+\frac{i}{2}$ |
| $-1$ | $1+i$ | $i$ | $1-i$ | $-i$ | $1-i$ | $i$ | $2i$ | $2-2i$ | $-1+i$ |
| $-1$ | $1+i$ | $i$ | $i$ | $-i$ | $1-i$ | $1-i$ | $1$ | $1$ | $-i$ |
| $-1$ | $i$ | $-i$ | $-i$ | $1-i$ | $1+i$ | $1+i$ | $2i$ | $-1-i$ | $2+2i$ |
| $-1$ | $i$ | $-i$ | $1+i$ | $1-i$ | $1+i$ | $-i$ | $-1$ | $i$ | $1$ |
| $-1$ | $i$ | $\frac{1}{2}+\frac{i}{2}$ | $-1$ | $1-i$ | $\frac{1}{2}-\frac{i}{2}$ | $2$ | $1+i$ | $-\frac{1}{2}-\frac{i}{2}$ | $1$ |
| $-1$ | $i$ | $\frac{1}{2}+\frac{i}{2}$ | $\frac{1}{2}+\frac{i}{2}$ | $1-i$ | $\frac{1}{2}-\frac{i}{2}$ | $\frac{1}{2}-\frac{i}{2}$ | $-\frac{1}{2}+\frac{i}{2}$ | $1$ | $-\frac{1}{2}-\frac{i}{2}$ |
| $-1$ | $\frac{1}{2}$ | $2$ | $\frac{1}{2}-\frac{i}{2}$ | $\frac{1}{2}$ | $-1$ | $\frac{1}{2}+\frac{i}{2}$ | $\frac{i}{2}$ | $-\frac{1}{2}-\frac{i}{2}$ | $-\frac{1}{2}+\frac{i}{2}$ |
| $-1$ | $\frac{1}{2}$ | $2$ | $\frac{1}{2}+\frac{i}{2}$ | $\frac{1}{2}$ | $-1$ | $\frac{1}{2}-\frac{i}{2}$ | $-\frac{i}{2}$ | $-\frac{1}{2}+\frac{i}{2}$ | $-\frac{1}{2}-\frac{i}{2}$ |
| $-1$ | $2$ | $1-i$ | $1-i$ | $-1$ | $i$ | $i$ | $1+i$ | $-1+i$ | $1$ |
| $-1$ | $2$ | $1+i$ | $1+i$ | $-1$ | $-i$ | $-i$ | $1-i$ | $-1-i$ | $1$ |

Recall that each fundamental element of $\mathbb{H}_2$ has a norm between $\frac{1}{2}$ and 2. By taking the base-2 logarithm of the norms we obtain linear bounds on the exponents:

```
Simplify[Log[2, Abs[imtable]]];
```

```
constraints = -1 ≤ #.{0, r, s, t, u, v, w, x, y, z} ≤ 1 & /@ %
```

$$\Big\{-1 \le \frac{s}{2} + u + \frac{w}{2} + \frac{x}{2} + \frac{z}{2} \le 1,\ -1 \le \frac{s}{2} + u + \frac{w}{2} + \frac{x}{2} + \frac{z}{2} \le 1,$$
$$-1 \le -\frac{s}{2} + \frac{u}{2} - \frac{v}{2} + w + \frac{x}{2} - \frac{y}{2} \le 1,\ -1 \le -\frac{s}{2} - \frac{t}{2} + \frac{u}{2} - \frac{v}{2} - \frac{w}{2} - \frac{x}{2} - \frac{z}{2} \le 1,\ -1 \le \frac{t}{2} + \frac{u}{2} + \frac{v}{2} \le 1,$$
$$-1 \le \frac{u}{2} + \frac{v}{2} + \frac{w}{2} + x + \frac{y}{2} + \frac{3z}{2} \le 1,\ -1 \le -\frac{r}{2} - \frac{t}{2} - \frac{u}{2} + v - \frac{w}{2} \le 1,\ -1 \le -\frac{r}{2} + \frac{s}{2} + \frac{t}{2} - \frac{u}{2} - \frac{x}{2} + \frac{z}{2} \le 1,$$
$$-1 \le -\frac{r}{2} + \frac{s}{2} - \frac{u}{2} + \frac{w}{2} - \frac{x}{2} + \frac{y}{2} \le 1,\ -1 \le -\frac{r}{2} - s - t - \frac{u}{2} - v - w - x - \frac{3y}{2} - \frac{3z}{2} \le 1,$$
$$-1 \le -\frac{r}{2} - \frac{t}{2} - \frac{u}{2} + v - \frac{w}{2} \le 1,\ -1 \le -\frac{r}{2} + \frac{s}{2} - \frac{u}{2} + \frac{w}{2} - \frac{x}{2} + \frac{y}{2} \le 1,$$
$$-1 \le -\frac{r}{2} + \frac{s}{2} + \frac{t}{2} - \frac{u}{2} - \frac{x}{2} + \frac{z}{2} \le 1,\ -1 \le -\frac{r}{2} - s - t - \frac{u}{2} - v - w - x - \frac{3y}{2} - \frac{3z}{2} \le 1,$$
$$-1 \le \frac{r}{2} + \frac{v}{2} + \frac{w}{2} \le 1,\ -1 \le \frac{r}{2} + \frac{t}{2} + \frac{v}{2} + x + \frac{3y}{2} + \frac{z}{2} \le 1,\ -1 \le \frac{r}{2} - \frac{s}{2} - \frac{t}{2} - \frac{v}{2} - \frac{w}{2} - \frac{x}{2} - \frac{y}{2} \le 1,$$
$$-1 \le \frac{r}{2} - \frac{s}{2} + t - \frac{v}{2} + \frac{x}{2} - \frac{z}{2} \le 1,\ -1 \le \frac{r}{2} - \frac{s}{2} - \frac{t}{2} - \frac{v}{2} - \frac{w}{2} - \frac{x}{2} - \frac{y}{2} \le 1,$$
$$-1 \le \frac{r}{2} - \frac{s}{2} + t - \frac{v}{2} + \frac{x}{2} - \frac{z}{2} \le 1,\ -1 \le \frac{r}{2} + \frac{t}{2} + \frac{v}{2} + x + \frac{3y}{2} + \frac{z}{2} \le 1,\ -1 \le \frac{r}{2} + \frac{v}{2} + \frac{w}{2} \le 1,$$
$$-1 \le \frac{u}{2} + \frac{v}{2} + \frac{w}{2} + x + \frac{y}{2} + \frac{3z}{2} \le 1,\ -1 \le \frac{t}{2} + \frac{u}{2} + \frac{v}{2} \le 1,\ -1 \le -\frac{s}{2} + \frac{u}{2} - \frac{v}{2} + w + \frac{x}{2} - \frac{y}{2} \le 1,$$
$$-1 \le -\frac{s}{2} - \frac{t}{2} + \frac{u}{2} - \frac{v}{2} - \frac{w}{2} - \frac{x}{2} - \frac{z}{2} \le 1,\ -1 \le -r + s - \frac{t}{2} - u - \frac{w}{2} - x - \frac{y}{2} - \frac{z}{2} \le 1,$$
$$-1 \le -r + s - \frac{t}{2} - u - \frac{w}{2} - x - \frac{y}{2} - \frac{z}{2} \le 1,\ -1 \le r + \frac{s}{2} + \frac{t}{2} + \frac{x}{2} + \frac{y}{2} \le 1,\ -1 \le r + \frac{s}{2} + \frac{t}{2} + \frac{x}{2} + \frac{y}{2} \le 1\Big\}$$

Now we find the maximum and minimum value for $u, v, w, x, y, z$ subject to the linear constraints that we obtained.

```
G = {r, s, t, u, v, w, x, y, z, -r, -s, -t, -u, -v, -w, -x, -y, -z};
```

```
Table[{G[[i]],
  Maximize[{G[[i]], And @@ constraints}, {r, s, t, u, v, w, x, y, z}, Integers][[1]]}, {i,
  1, Length[G]}]
```

```
{{r, 1}, {s, 1}, {t, 1}, {u, 1}, {v, 1}, {w, 1}, {x, 2}, {y, 1}, {z, 1},
 {-r, 1}, {-s, 1}, {-t, 1}, {-u, 1}, {-v, 1}, {-w, 1}, {-x, 2}, {-y, 1}, {-z, 1}}
```

In other words, $-1 \le r, s, t, u, v, w, y, z \le 1$, $-2 \le x \le 2$. Now we have reduced the set of possible exponents of fundamental elements to a finite list.

```
candidateH5funExps =
  Flatten[Table[{q, r, s, t, u, v, w, x, y, z}, {q, 0, 1}, {r, -1, 1}, {s, -1, 1}, {t, -1, 1},
    {u, -1, 1}, {v, -1, 1}, {w, -1, 1}, {x, -2, 2}, {y, -1, 1}, {z, -1, 1}], 9];
```

```
Length[%]
```

```
65 610
```

To find out which elements are fundamental we map the candidates to GF($p$) for some large $p$, and choose the images of $\alpha$, $\beta$ so that this map is injective. Then we select the elements $x$ for which $1-x$ is also in the list. To avoid dealing with fractions we reconstruct the list using the PowerMod function:

```
p = Prime[1 000 000 000]
```

```
22 801 763 489
```

```
aIm = 17; bIm = 47; cIm = 53;
```

```
H5genIms = Mod[H5gens /. {α → aIm, β → bIm, γ → cIm}, p]
```

```
{22 801 763 488, 17, 47, 53, 22 801 763 473,
 22 801 763 443, 22 801 763 437, 22 801 763 453, 22 801 762 743, 700}
```

```
candidateH5funIms =
  {0} ⋃ (Mod[Times @@ PowerMod[H5genIms, #, p], p] & /@ candidateH5funExps);
```

This map is indeed injective (the union symbol removes duplicates; note that we have added 0 so the list has one element more than the list of exponents)

```
Length[%]
```

```
65 611
```

Select the elements $x$ for which $1-x$ is also in the list:

```
H5funIms = candidateH5funIms ⋂ Mod[1 - candidateH5funIms, p];
```

```
Length[%]
```

```
92
```

This list of remaining candidates necessarily contains the images of all known fundamental elements. Since the number of remaining elements equals the number of known fundamental elements, we conclude that those are, indeed, all. This completes the proof of Lemma 27. We record the corresponding exponents for use below. First we recompute the images, this time not taking the union to avoid sorting.

```
candidateH5funImsNoZero = Mod[Times @@ PowerMod[H5genIms, #, p], p] & /@ candidateH5funExps;
```

Next we find the positions of the images of fundamental elements in this list, and take the exponents with corresponding indices. We remove the exponent corresponding to 1.

```
H5funExps = Complement[
   candidateH5funExps[[Flatten[Position[candidateH5funImsNoZero, #] & /@ H5funIms]]],
   {{0, 0, 0, 0, 0, 0, 0, 0, 0, 0}}];
```

```
Length[H5funExps]
```

```
90
```

## The automorphism group

In this section we prove the following result:

**Lemma 28.** The automorphism group of $\mathbb{H}_5$ is isomorphic to $S_6$.

The automorphisms correspond with coordinate permutations in $GF(5)\times GF(5)\times GF(5)\times GF(5)\times GF(5)\times GF(5)$.

The partial-field automorphism group of $\mathbb{H}_5$ is determined uniquely by the images of $\alpha$, $\beta$, and $\gamma$. This is easy to see since each automorphism maps 0 to 0, 1 to 1, and $-1$ to $-1$, so the real number line is fixed. Note that $\alpha$, $\beta$, and $\gamma$ have to be mapped to fundamental elements, and those fundamental elements have to be distinct, so the group will be finite. Note that we only have to check if $\left\{\alpha, \beta, \gamma, \frac{\alpha\beta}{\gamma}, \frac{\alpha}{\gamma}, \frac{(1-\alpha)\gamma}{\gamma-\alpha}, \frac{(\alpha-1)\beta}{\gamma-1}, \frac{\alpha-1}{\gamma-1}, \frac{\gamma-\alpha}{\gamma-\alpha\beta}, \frac{(\beta-1)(\gamma-1)}{\beta(\gamma-\alpha)}, \frac{\beta(\gamma-\alpha)}{\gamma-\alpha\beta}, \frac{(\alpha-1)(\beta-1)}{\gamma-\alpha}, \frac{\beta(\gamma-\alpha)}{(1-\gamma)(\gamma-\alpha\beta)}, \frac{(1-\alpha)(\gamma-\alpha\beta)}{\gamma-\alpha}, \frac{1-\beta}{\gamma-\alpha\beta}\right\}$ map to fundamental elements: the remainder follows automatically. To speed matters up, we first check if the images of fundamental elements under the homomorphism to $GF(p)$ from the previous section map to fundamental elements. Only on the remaining sets will we test for equality of the rational expressions.

Note that 0 and 1 are mapped to themselves, so we do not have to test for those.

```
nonzerooneH5funIms = Mod[Times @@ PowerMod[H5genIms, #, p], p] & /@ H5funExps;
```

```
Length[%]
```

```
90
```

We recompute the next line to ensure the ordering is identical to the ordering of the images:

```
nonzerooneH5funs = (Times @@ (H5gens^#)) & /@ H5funExps;
```

```
validImage[aIm_, bIm_, cIm_] := Module[{nonzerooneH5funImsp, H5genIms},
   (* Check if the generators get mapped to nonzero numbers *)
   H5genIms = Mod[H5gens /. {α → aIm, β → bIm, γ → cIm}, p];
   If[MemberQ[H5genIms, 0], Return[False]];
   (* Construct the list of images of fundamental elements *)
   nonzerooneH5funImsp := Mod[Times @@ PowerMod[H5genIms, #, p], p] & /@ H5funExps;
   (* Check if this list equals the old list of images,
   and return True or False accordingly *)
   Length[Complement[nonzerooneH5funIms, nonzerooneH5funImsp]] == 0
  ];
```

Computing the following table with True/False entries took about 30 minutes on an average laptop from 2008. If you are using *Mathematica* 7, parallelization should help speed this up.

```
autTable = Table[validImage[nonzerooneH5funIms[[x]], nonzerooneH5funIms[[y]],
    nonzerooneH5funIms[[z]]], {x, 1, Length[nonzerooneH5funIms]},
   {y, 1, Length[nonzerooneH5funIms]}, {z, 1, Length[nonzerooneH5funIms]}];
```

```
autPositions = Position[%, True];
```

```
H5autImages = nonzerooneH5funs[[#]] & /@ %;
```

```
Length[%]
```

```
720
```

The images of the representative set of fundamental elements:

```
candidateAuts =
  Table[{α, β, γ, αβ/γ, α/γ, (1 - α) γ/(γ - α), (α - 1) β/(γ - 1), (α - 1)/(γ - 1), (γ - α)/(γ - α β), (β - 1) (γ - 1)/(β (γ - α)), β (γ - α)/(γ - α β),
     (α - 1) (β - 1)/(γ - α), β (γ - α)/((1 - γ) (γ - α β)), (1 - α) (γ - α β)/(γ - α), (1 - β)/(γ - α β)} /. {α → H5autImages[[x, 1]],
     β → H5autImages[[x, 2]], γ → H5autImages[[x, 3]]}, {x, 1, Length[H5autImages]}];
```

For these remaining ones we will do the (much slower) full check. This took about an hour on an average laptop from 2008.

```
isH5fun[x_] := occursInList[nonzerooneH5funs, x];
areH5funs[x_List] := And @@ (isH5fun /@ x);
```

```
H5auts = Select[candidateAuts, areH5funs];
```

```
Length[H5auts]
```

```
720
```

Conclusion: the automorphism group of $\mathbb{H}_5$ has 720 elements and is isomorphic to $S_6$, the symmetric group on 6 symbols. Indeed, the partial-field homomorphism $\phi 6$ defined above induces a $1-1$ correspondence between the automorphisms of $\mathbb{H}_2$ and coordinate permutations of $\mathrm{GF}(5)\times\mathrm{GF}(5)\times\mathrm{GF}(5)\times\mathrm{GF}(5)\times\mathrm{GF}(5)\times\mathrm{GF}(5)$.

## Representations of $U_{2,5}$

We prove the following fact:

**Lemma 29.** Let $M$ be a 3-connected, $\mathbb{H}_5$-representable matroid. If $M$ has a $U_{2,5}$- or $U_{3,5}$-minor then $M$ has at least six inequivalent representations over GF(5).

To prove this it suffices to show that each $\mathbb{H}_5$-representation of $U_{2,5}$ gives rise to six inequivalent representations over GF(5). Since any $\mathbb{H}_5$-representation matrix of $M$ must contain one of these as a minor, the result then follows. First we enumerate all such representations. By normalizing and suppressing the identity matrix at the front, we see

```
Clear[p];
```

```
A = {{1, 1, 1}, {1, p, q}}; MatrixForm[A]
```

$$\begin{pmatrix} 1 & 1 & 1 \\ 1 & p & q \end{pmatrix}$$

This is an $\mathbb{H}_5$-matrix representing $U_{2,5}$ if and only if $p$, $q$, $\frac{p}{q}$ are fundamental elements, $p$ and $q$ are not equal to 0 or 1, and $p \neq q$. Moreover, two such matrices are equivalent if and only if they are equal.

Create a table of all candidate pairs $p$, $q$. We suppress the output:

```
candidatepqPairs = Flatten[Table[{nonzerooneH5funs[[x]], nonzerooneH5funs[[y]]},
     {x, 1, Length[nonzerooneH5funs]}, {y, 1, x - 1}], 1] ⋃
   Flatten[Table[{nonzerooneH5funs[[x]], nonzerooneH5funs[[y]]},
     {x, 1, Length[nonzerooneH5funs]}, {y, x + 1, Length[nonzerooneH5funs]}], 1];
```

Filter out those for which $p/q$ is not fundamental or equal to 1:

```
pqPairs = Select[candidatepqPairs, occursInList[nonzerooneH5funs, #[[1]] / #[[2]]] &];
```

```
Length[%]
```

```
720
```

To show that each of these gives rise to six inequivalent representations over GF(5) we consider the partial-field homomorphism $\phi 5$.

Let $\xi_i : \text{GF}(5)\times\text{GF}(5)\times\text{GF}(5)\times\text{GF}(5)\times\text{GF}(5)\times\text{GF}(5) \to \text{GF}(5)$ be the projection onto the $i$th coordinate. We have to show that $\xi_i(\phi(A))$ is not equivalent to $\xi_j(\phi(A))$ for all matrices $A$ computed above. We only need to look at the pairs $p$, $q$. These are their images under $\phi 5$:

```
GF5t6pqpairs = Map[ϕ5, pqPairs, {2}];
```

Indeed, each pair of representations is inequivalent. A concise way of checking this is to show that for all $1 \le i < j \le 4$, $|\xi_i(p) - \xi_j(p)| + |\xi_i(q) - \xi_j(q)| > 0$.

```
Map[{Abs[#[[1, 1]] - #[[1, 2]]] + Abs[#[[2, 1]] - #[[2, 2]]],
    Abs[#[[1, 1]] - #[[1, 3]]] + Abs[#[[2, 1]] - #[[2, 3]]],
    Abs[#[[1, 1]] - #[[1, 4]]] + Abs[#[[2, 1]] - #[[2, 4]]],
    Abs[#[[1, 1]] - #[[1, 5]]] + Abs[#[[2, 1]] - #[[2, 5]]],
    Abs[#[[1, 1]] - #[[1, 6]]] + Abs[#[[2, 1]] - #[[2, 6]]],
    Abs[#[[1, 2]] - #[[1, 3]]] + Abs[#[[2, 2]] - #[[2, 3]]],
    Abs[#[[1, 2]] - #[[1, 4]]] + Abs[#[[2, 2]] - #[[2, 4]]], Abs[#[[1, 2]] - #[[1, 5]]] +
     Abs[#[[2, 2]] - #[[2, 5]]], Abs[#[[1, 2]] - #[[1, 6]]] + Abs[#[[2, 2]] - #[[2, 6]]],
    Abs[#[[1, 3]] - #[[1, 4]]] + Abs[#[[2, 3]] - #[[2, 4]]],
    Abs[#[[1, 3]] - #[[1, 5]]] + Abs[#[[2, 3]] - #[[2, 5]]],
    Abs[#[[1, 3]] - #[[1, 6]]] + Abs[#[[2, 3]] - #[[2, 6]]], Abs[#[[1, 4]] - #[[1, 5]]] +
     Abs[#[[2, 4]] - #[[2, 5]]], Abs[#[[1, 4]] - #[[1, 6]]] + Abs[#[[2, 4]] - #[[2, 6]]],
    Abs[#[[1, 5]] - #[[1, 6]]] + Abs[#[[2, 5]] - #[[2, 6]]]} &, %, {1}];
```

```
Position[%, 0]
```

```
{}
```

This completes the proof of Lemma 28.

## Lifting

We prove the following fact:

**Lemma 29.** Let $M$ be a 3-connected matroid with at least five inequivalent representations over GF(5). Then $M$ is representable over $\mathbb{H}_5$.

Consider a $\mathrm{GF}(5)\times\mathrm{GF}(5)\times\mathrm{GF}(5)\times\mathrm{GF}(5)\times\mathrm{GF}(5)$-matrix $A$ representing $M$ such that the five projections $\xi_i(A)$ are pairwise inequivalent. We want to apply the Lift Theorem, so we start by constructing a lifting function for $\mathrm{Cr}(A)$. First we need to find out what $\mathrm{Cr}(A)$ is.

**Claim 29.1.** No element of $\mathrm{Cr}(A)$ is of the form $(x, x, x, y, z)$ for $x \notin \{0, 1\}$.

*Proof.* Assume that $(x, x, x, y, z) \in \mathrm{Cr}(A)$ for some $x \in \{2, 3, 4\}$. Let $\mathbb{P}'$ be the sub-partial field of $\mathrm{GF}(5)\times\mathrm{GF}(5)\times\mathrm{GF}(5)\times\mathrm{GF}(5)\times\mathrm{GF}(5)$ in which the first three coordinates of all elements are equal. Since $A$ is not a $\mathbb{P}'$-matrix, the Confinement Theorem, in particular Corollary 15, implies that either $A$ or $A^T$ contains a minor

$$D = \begin{pmatrix} (1, 1, 1, 1, 1) & (1, 1, 1, 1, 1) & (1, 1, 1, 1, 1) \\ (1, 1, 1, 1, 1) & (x, x, x, y, z) & (p, q, r, s, t) \end{pmatrix},$$

with not all three of $p$, $q$, $r$ equal. There are only two elements $u$ of GF(5) such that $\begin{pmatrix} 1 & 1 & 1 \\ 1 & x & u \end{pmatrix}$ represents $U_{2,5}$. Hence two of the first three projections of $D$ must be equivalent. But then the Stabilizer Theorem, in particular Lemma 8, implies that two of the five projections of $A$ are equivalent, a contradiction. □

By considering all permutations of the coordinates we are left with the following cross ratios:

$F' := \mathcal{F}(\mathrm{GF}(5)\times\mathrm{GF}(5)\times\mathrm{GF}(5)\times\mathrm{GF}(5)\times\mathrm{GF}(5)) \setminus \{x : \text{at least 3 coordinates of } x \text{ are equal},\ x \notin \{(0, 0, 0, 0, 0), (1, 1, 1, 1, 1)\}\}$.

```
Fp = Complement[{{0, 0, 0, 0, 0}, {1, 1, 1, 1, 1}} ⋃
    Flatten[Table[{v, w, x, y, z}, {v, 2, 4}, {w, 2, 4}, {x, 2, 4}, {y, 2, 4}, {z, 2, 4}], 4],
   Flatten[Table[{{x, x, x, y, z}, {x, x, x, z, y}, {x, x, y, x, z}, {x, x, y, z, x},
      {x, x, z, x, y}, {x, x, z, y, x}, {x, y, x, x, z}, {x, y, x, z, x},
      {x, y, z, x, x}, {x, z, x, x, y}, {x, z, x, y, x}, {x, z, y, x, x}, {y, x, x, x, z},
      {y, x, x, z, x}, {y, x, z, x, x}, {y, z, x, x, x}, {z, x, x, x, y}, {z, x, x, y, x},
      {z, x, y, x, x}, {z, y, x, x, x}}, {x, 2, 4}, {y, 2, 4}, {z, 2, 4}], 3]];
```

```
Length[%]
```

```
92
```

Consider the homomorphism $\psi$, obtained from the homomorphism $\phi5$ by dropping the last coordinate. The restriction of $\psi$ to $\mathcal{F}(\mathbb{H}_5)$ is a bijection between $\mathcal{F}(\mathbb{H}_5)$ and $F'$:

```
ψ[0] = {0, 0, 0, 0, 0};
ψ[1] = {1, 1, 1, 1, 1};
ψ[-1] = {4, 4, 4, 4, 4};
ψ[x_] := toGF5[x /. {α → {4, 3, 3, 4, 2}, β → {3, 2, 2, 4, 3}, γ → {3, 2, 4, 2, 3}}];
```

```
Complement[ψ /@ H5funs, Fp]
```

```
{}
```

It follows immediately that the function $F' \to \mathbb{H}_5$ that is the inverse of $\psi$ is a lifting function.

```
up5[v_, w_, x_, y_, z_] := H5funs[[Position[ψ /@ H5funs, {v, w, x, y, z}][[1, 1]]]];
up5[{v_, w_, x_, y_, z_}] := up5[v, w, x, y, z];
```

Some examples:

```
up5[1, 1, 1, 1, 1]
```

```
1
```

```
up5[2, 4, 3, 3, 4]
```

$$\frac{\gamma}{\alpha}$$

Now suppose that $A$ has no global lift. Then the Lift Theorem (Theorem 12) tells us that $A$ or $A^T$ has a minor of the form

$$\begin{pmatrix} 1 & 1 & 0 & 1 \\ 1 & 0 & 1 & 1 \\ 0 & 1 & 1 & 1 \end{pmatrix} \text{ or } \begin{pmatrix} 1 & 1 & 1 \\ 1 & p & q \end{pmatrix}, \text{ for some } p, q \in \mathrm{Cr}(A).$$

without a local lift. The first of these has a cross ratio (2, 2, 2, 2, 2), so it does not occur as a minor of $A$, by Claim 29.1. Therefore we have to ensure that all $\mathrm{GF}(5) \times \mathrm{GF}(5) \times \mathrm{GF}(5) \times \mathrm{GF}(5)$-representations of $U_{2,5}$ that can occur in $A$ have a local lift. Let $D$ be such a minor. Suppose that $\xi_i(D)$ and $\xi_j(D)$ are equivalent. By the Stabilizer Theorem, and in particular Lemma 8, this implies that $\xi_i(A)$ and $\xi_j(A)$ are equivalent, a contradiction to our choice of $A$. Hence $\xi_i(D)$ and $\xi_j(D)$ are inequivalent. This gives us $6 \times 5 \times 4 \times 3 = 360$ representations of $U_{2,5}$. Again we only list the pairs $(p, q)$ in

$$D = \begin{pmatrix} 1 & 1 & 1 \\ 1 & p & q \end{pmatrix}.$$

```
U25reps = { {2, 3}, {2, 4}, {3, 2}, {3, 4}, {4, 2}, {4, 3}}
```

```
{{2, 3}, {2, 4}, {3, 2}, {3, 4}, {4, 2}, {4, 3}}
```

```
U25repquints = Transpose /@ Permutations[U25reps, {5}];
```

```
Length[%]
```

```
720
```

We compute the lifts:

```
H5lifts = Map[up5, U25repquints, {2}];
```

And indeed, for each of these $p^{\uparrow} / q^{\uparrow}$ is a fundamental element:

```
Select[H5lifts, Not[occursInList[nonzerooneH5funs, #[[1]] / #[[2]]]] &]
```

```
{}
```

It follows that $A$ does have a global lift. This completes the proof of Lemma 29.

# Appendix: generating Hydra-3

So far this document has focused on the verification of our results. However, the Hydra-$k$ partial fields did not appear out of thin air. Admittedly, the Gaussian partial field, $\mathbb{H}_2$ was invented by the first author, and in fact set this research in motion, but the others were actually computed as the lift partial field of $\text{GF}(5) \times \ldots \times \text{GF}(5)$, using Definition 5.6 from **[lift]**. As an example we carry out this computation for $\mathbb{H}_3$.

We will not use all generators, but take those that generate all cross ratios:

```
cratgens = {{2, 2, 3}, {2, 3, 2}, {2, 3, 3}, {2, 3, 4}};
```

Indeed these generate all cross ratios other than 0,1 (cf. Lemma 22).

By restricting ourselves to these cross ratios we ensure that 5.6(i)--5.6(iii) hold. For 5.6(iv), we need to ensure that, for every triple $p$, $q$, $r$ such that $p\,q\,r = 1$, the lifted relation $\tilde{p}\,\tilde{q}\,\tilde{r} = 1$ holds. We do this only for three of the triples $p$, $q$, $r$. It turns out that those suffice to determine the partial field completely.

```
pqrtriples = {{toGF5[1 - {2, 3, 4}], toGF5[({2, 2, 3} - 1)/{2, 2, 3}], toGF5[({2, 3, 4} - 1)/{2, 3, 4}]},
   {toGF5[({2, 3, 2} - 1)/{2, 3, 2}], toGF5[1/{2, 3, 4}], toGF5[1/(1 - {2, 3, 4})]},
   {toGF5[({2, 3, 4} - 1)/{2, 3, 4}], toGF5[1/{2, 3, 4}], toGF5[1/(1 - {2, 3, 3})]}}
```

```
{{{4, 3, 2}, {3, 3, 4}, {3, 4, 2}},
 {{3, 4, 3}, {3, 2, 4}, {4, 2, 3}}, {{3, 4, 2}, {3, 2, 4}, {4, 2, 2}}}
```

Indeed, the product of each of these triples equals (1,1,1):

```
toGF5[#[[1]] #[[2]] #[[3]]] & /@ pqrtriples
```

```
{{1, 1, 1}, {1, 1, 1}, {1, 1, 1}}
```

The definition shows how to write these as rational function of the generators above. By replacing $(x, y, z)$ by a variable $s\,x\,y\,z$ and multiplying out the denominators in the expressions $1 - p\,q\,r$ we get the following three polynomials:

```
polrels = {s223 s234 - (1 - s234) (s223 - 1) (s234 - 1),
(1 - s234) s234 s232 - (s232 - 1),
s234 s234 (1 - s233) - (s234 - 1)};
```

We compute a Groebner basis over the integers for this ideal:

```
GroebnerBasis[polrels, {s223, s232, s233, s234}, CoefficientDomain → Integers]
```

$$\{-1 + s234 - s234^2 + s233\ s234^2,\ -1 + s232 - s232\ s234 + s232\ s234^2,\ -s233 + s232\ s233 + s234 - s233\ s234,\ -1 + s223 + s232\ s234\}$$

We choose $\alpha$ := s234, and express the other generators in terms of $\alpha$:

```
Solve[(% /. s234 → α) == 0, {s223, s232, s233}]
```

$$\left\{\left\{s223 \to 1 - \frac{\alpha}{1 - \alpha + \alpha^2},\ s233 \to \frac{1 - \alpha + \alpha^2}{\alpha^2},\ s232 \to \frac{1}{1 - \alpha + \alpha^2}\right\}\right\}$$

...and $\mathbb{H}_3$ is born!

# Appendix: a partial Mathematica command legend

A Mathematica notebook is an interactive environment, consisting of units called cells. A cell can contain text (such as this paragraph), code (input), or output. If a cell containing code is executed, the output appears directly below it. In this document, a sequence of such executions is recorded.

Among other things, the following commands were used:

*l*[[2]] is the second element of a list
*l*[[2, 4]] is the element with coordinates (2,4) in a 2-dimensional array or matrix
*l*[[{2, 4}]] is the sublist of *l* consisting of the 2nd and 4th element
*f*[*x*] is the function *f* applied to *x*.
*f* @@ {*a*, *b*, *c*} = *f*[*a*, *b*, *c*]
*f* /@ {*a*, *b*, *c*} = {*f*[*a*], *f*[*b*], *f*[*c*]}
... ♯ ... & denotes a (pure) function. When invoked, the single argument is put in all places where the ♯ symbol appears.
_? *f* is a pattern. It matches all expressions *x* for which *f*[*x*] yields True.
*f* /. *x* → *y* replace all occurrences of *x* in *f* by *y*.
% previous line of output
%% second to last line of output
A semicolon at the end of a statement suppresses the output that would otherwise be printed directly below its execution.

For the remaining functionality that we used we refer to http://reference.wolfram.com/ .